\documentclass[preprint]{elsarticle}

\usepackage[margin=3cm]{geometry}% by courtesy of Mico

\usepackage{lineno,hyperref}
\modulolinenumbers[5]
\usepackage{amssymb}

\usepackage{subfigure}

\usepackage{amsmath} %Need this for e.g. for hdots
\usepackage{color} %To color comments. Can be removed later.
\usepackage{mathtools} %To use colonequals

\usepackage[]{algorithm2e} %To write algorithms

\newcommand{\mb}[1] {\mathbf{#1}}                           % Short hand for bold-face in mathenv.

\newcommand{\drop}[1] {\Omega_{#1}}                         % Drop interior
\newcommand{\dinterf}[1] {\Gamma_{#1}}                      % Drop interface
\newcommand{\nDrop} {N_\Gamma}                              % Number of drops
\newcommand{\sinterf}[1] {\gamma_{#1}}                      % Solid boundary
\newcommand{\nSolid} {N_\gamma}                             % Number of solids

\newcommand{\norm}[1] {\mb{\hat{n}}_{#1}}                   % Normal vector (real)
\newcommand{\normC}[1] {\hat{n}_{#1}}                       % Normal vector (complex)
\newcommand{\SPC}[1] {\beta(\mb{x})} % Solid particle contribution
\newcommand{\SPCk}[2] {\beta_{k}[\mb{{#1}}](\mb{{#2}})} % Solid particle contribution
\newcommand{\SLe} {\mathbb{S}} % Single layer
\newcommand{\DLe} {\mathbb{D}} % Double layer
\newcommand{\SL}[3] {\SLe_{#1}[\mb{#2}](\mb{#3})}
\newcommand{\DL}[3] {\DLe_{#1}[\mb{#2}](\mb{#3})}
\newcommand{\SLp}[3] {\SLe^{\,P}_{#1}[\mb{#2}](\mb{#3})}
\newcommand{\DLp}[3] {\DLe^{\,P}_{#1}[\mb{#2}](\mb{#3})}
\newcommand{\SLc}[3] {\SLe_{#1}[{#2}]({#3})}
\newcommand{\DLc}[3] {\DLe_{#1}[{#2}]({#3})}
\newcommand{\SLpc}[3] {\SLe^{\,P}_{#1}[{#2}]({#3})}
\newcommand{\DLpc}[3] {\DLe^{\,P}_{#1}[{#2}]({#3})}

\newcommand{\pSolid}[1] {M^{\sinterf{}}_{{#1}}}
\newcommand{\pDrop}[1] {M^{\dinterf{}}_{{#1}}}
\newcommand{\ptot} {M}
\newcommand{\pGen} {\ptot^{\Lambda}}

\newcommand{\rspace} {``real space"}
\newcommand{\kspace} {``$\mb{k}$-space"}
\newcommand{\rspaces} {\rspace ~sum}
\newcommand{\kspaces} {\kspace ~sum}

\newcommand{\biharm} {\mathcal{B}}
\newcommand{\lapl} {\mathcal{L}}
\newcommand{\dd}[0] {\textrm{d}}
\newcommand{\DD}[0] {\textrm{D}}
\newcommand{\OO}[0] {\mathcal{O}}

\usepackage{xcolor}
\definecolor{dodgerblue}{rgb}{0.12, 0.56, 1.0}
\definecolor{deepsaffron}{rgb}{1.0, 0.6, 0.2}
\definecolor{pastelviolet}{rgb}{0.8, 0.6, 0.79}
\definecolor{orchid}{rgb}{0.85, 0.44, 0.84}

\newcommand{\myChange}[1]{\textcolor{black}{#1}}
\newcommand{\revOne}[1]{\textcolor{black}{#1}}
\newcommand{\revTwo}[1]{\textcolor{black}{#1}}

\makeatletter
\def\ps@pprintTitle{%
 \let\@oddhead\@empty
 \let\@evenhead\@empty
 \def\@oddfoot{}%
 \let\@evenfoot\@oddfoot}
\makeatother

\begin{document}

\begin{frontmatter}

\title{An integral equation method for closely interacting surfactant-covered droplets in wall-confined Stokes flow}
\author[add1]{Sara P\aa lsson \corref{cor1}}
\ead{sarapal@kth.se}
\author[add1]{Anna-Karin Tornberg}

\cortext[cor1]{Corresponding author}
\address[add1]{Numerical Analysis, Department of Mathematics, KTH Royal Institute of Technology, 100 44 Stockholm, SE}

\begin{abstract}
A highly accurate method for simulating surfactant-covered droplets in two-dimensional Stokes flow with solid boundaries is presented. The method handles both periodic channel flows of arbitrary shape and stationary solid constrictions. A boundary integral method together with a special quadrature scheme is applied to solve the Stokes equations to high accuracy, also \myChange{for closely interacting droplets}. The problem is considered in a periodic setting and \myChange{Ewald decompositions} for the Stokeslet and stresslet \myChange{are} derived. Computations are \myChange{accelerated} using the spectral Ewald method. The time evolution is handled with a fourth order, adaptive, implicit-explicit time-stepping scheme. The numerical method is tested through several convergence studies and other challenging examples and is shown to handle drops in close proximity both to other drops and solid objects to high accuracy.
\end{abstract}


\end{frontmatter}

%\definecolor{mycolor1}{rgb}{0.1171875,0.54296875,0.76171875}% %
%\definecolor{greycolor1}{rgb}{0.45313,0.45313,0.45313}%
% \definecolor{mycolor1_matlabblue}{rgb}{0.00000,0.44700,0.74100}%
% \definecolor{mycolor2_matlabred}{rgb}{0.85000,0.32500,0.09800}%

% \definecolor{mycolor1}{rgb}{0.11719,0.54297,0.76172}%
% \definecolor{mycolor2}{rgb}{0.45313,0.45313,0.45313}%

%% INTRODUCTION
% !TEX root = ms.tex

% ******************************************************************************
\section{Introduction}
\label{sec:intro}
% ******************************************************************************
The study of deforming droplets on the micro scale is motivated by several applications, one being the design of lab-on-a-chip-devices. In many cases, the study of deformable droplets in a confined flow is especially important. For example, the channel geometry can be used to control the behaviour of the droplets; regarding transport, splitting and fusing of the droplets \cite{droplet_reviewTan}. It is also of interest for the study of flow through porous media, which is relevant to many industries such as e.g. oil recovery. A review of the physics of the problem is given by \citeauthor{Zhang2018} \cite{Zhang2018}.

On \myChange{the} micro scale, the flow can be modelled by the Stokes equations. The surface area to volume ratio is typically very high and interfacial forces are important for the flow dynamics. Surfactants are molecules that alter the surface tension of a drop, which changes the interfacial dynamics and thus the behaviour of the whole system. The inclusion of surfactants is an important tool in drop creation and coalescence prevention \cite{Anna2016}. A review of drops and bubbles in shear flow without constrictions is given by \citeauthor{Rallison1984} \cite{Rallison1984}. Regarding the flow through straight capillaries, \citeauthor{Olbricht1992} \cite{Olbricht1992} performed an extensive experimental study of the shape of a drop as a function of several physical parameters, such as Capillary number and viscosity ratio.
Furthermore, \citeauthor{Shapira1988} \cite{Shapira1988} studied how the shape of a droplet between two parallel plates is affected by the ratio of droplet diameter and channel height using small deformation analysis and reflections. The study of physical parameters is clearly of interest, however this paper focuses on the development of a highly accurate numerical method for simulating deforming droplets in wall-bounded Stokes flow with the inclusion of stationary, solid constrictions.

% Overview of BIE for drops and constrictions/walls
The numerical method described in this paper contains a boundary integral equation method. An overview of some studies of deforming droplets with different kinds of constrictions using boundary integral methods follows: In 3D, \citeauthor{Zinchenko2006} \cite{Zinchenko2006} developed a method for simulating \myChange{a single} deforming drop squeezing through solid particle constrictions. This method is based on a Hebeker representation for the solid particle contribution and formulates a system of Fredholm integral equations of the second kind. They considered a uniform flow pushing the drop through constrictions and studied the effect of Capillary number and viscosity ratio between drop and bulk on the deformation. To handle interactions between droplets and particles, a special desingularization technique was used. \revTwo{Whilst allowing for simulations of droplets squeezing through narrow constrictions, this method demands special consideration for all combinations of near drop and solid evaluations, such as drop-drop, solid-drop, drop-solid etc. and was restricted to solid particles of spheroidal shapes.}
Furthermore, in \cite{Griggs2007} the authors considered the motion of a single drop between two parallel plane walls, but for this paper the solid walls were handled by modifying the Green's functions. \citeauthor{janssen2007} \cite{janssen2007} regarded the deformation of drops with unity viscosity ratio between parallel plates, and how the degree of confinement affects the behaviour of the drop. \myChange{Their} method is based on a boundary integral formulation in 3D, where the walls are taken into account by using Green's functions associated with the walls. They have since extended their method to handle also non-unity viscosity ratios \cite{janssen2008} and also regarded unity viscosity ratios of drops together with insoluble surfactants in
\cite{janssen2008a}. \myChange{Their} method is however restricted to flat parallel walls. \citeauthor{Tsai1994} \cite{Tsai1994} studied the dynamics of a 3D axisymmetric drop in straight and capillary tubes as a function of viscosity ratio and Capillary number.

In 2D, \citeauthor{Zhou1993} \cite{Zhou1993} considered suspensions of drops in channels, using a periodic suspension of viscous drops. Here the drops were ordered in a single file, but studies were also made for random suspensions \cite{Zhou1993a}. A Fredholm integral equation of the second kind was obtained using periodic Green's functions that represented the flow due to a periodic array of 2D point forces in a channel. In the above references, the flow was driven by the relative translation of the two walls. In \cite{Zhou1994}, the flow was driven by a constant pressure drop.
\citeauthor{Li2000} \cite{Li2000} considered wall-bounded channel flow of a suspension of many droplets, for different Capillary numbers and viscosity ratios. They used a boundary integral equation with modified periodic Green's functions to take the walls into account. \citeauthor{DeBisschop2002} \cite{DeBisschop2002} considered the motion of a two-dimensional bubble rising in an inclined channel in Stokes flow. They considered both clean and surfactant-covered bubbles. The fluid velocity was computed using a periodic Green's function (in the $x$-direction). The authors compared their results with that of experiments regarding inclined walls, showing good agreement.

%Other methods for drops and walls/constrictions
Other numerical methods to simulate deformable droplets in wall-bounded flow include \citeauthor{Lee2006} \cite{Lee2006}, who considered the effect of surfactants on the deformation of drops and bubbles in flow with non-zero Reynolds number. They used finite-differences for the Navier-Stokes equations, finite volumes for the insoluble surfactants and Peskin's immersed interface method for the interface tracking. In 3D, \citeauthor{Wang2012} \cite{Wang2012} studied deformable drops in a square channel using a  boundary element method. \citeauthor{MORTAZAVI2000} \cite{MORTAZAVI2000} studied three dimensional deforming drops in a tube using a finite difference/front-tracking scheme.
In \myChange{2D}, \citeauthor{Claus2018} \cite{Claus2018} used a cut finite element method to simulate deforming bubbles in Navier-Stokes flow. Their results included those of a bubble squeezing through a $5:1:5$ contraction/expansion micro-channel. \citeauthor{Chung2009} \cite{Chung2009} investigated the effect of viscosity ratio and Capillary number on a similar construction, using the finite element front-tracking method.

% Overview some 2D things for drops no walls and also vesicles with walls (a bit newer)
Surfactant-covered droplets in two-dimensional Stokes flow without the presence of walls and solid constrictions have been previously simulated with boundary integral equation methods, e.g. by \citeauthor{kropinskilushi2011} \cite{kropinskilushi2011} and \myChange{by the current authors}  \citeauthor{Palsson2018} \cite{Palsson2018}. The results in \cite{Palsson2018} were thoroughly validated using exact and semi-analytical solutions obtained by conformal mapping theory. The validation tests showed the ability of the method to obtain \myChange{high accuracy} (e.g. 8 correct digits) in solutions also after a long time with significant droplet deformation and close interactions between droplets.

Resolving the interactions of droplets in close proximity is a challenge for all numerical methods; grid based methods face the need for fine meshes and remeshing, whilst boundary integral equation methods necessitates the handling of nearly-singular integrals. The work in \cite{Palsson2018} utilised a special quadrature scheme \cite{Ojala2015} in order to resolve these interactions and achieve very accurate solutions as discussed above. The interactions of other objects in flow, such as vesicles, yield similar challenges.
\citeauthor{Rahimian2010} \cite{Rahimian2010} used a boundary integral method in 2D extending that by \citeauthor{kropinski2001} \cite{kropinski2001} to study how vesicles deform over time in confined flows, for example when squeezing through constrictions.
\citeauthor{Quaife2014} \cite{Quaife2014} studied vesicles suspended in a viscous Stokesian fluid, including channel constrictions and other solid geometries. This method was revisited in \cite{Quaife2016} with the inclusion of an adaptive time-stepping scheme.
\revTwo{To handle closely interacting vesicles, interpolation was used between an on-surface evaluation and evaluation points at a sufficient distance from the interface. The interpolation method generalises to 3D easily, but introduces several parameter selections for optimal use such as the spacing of interpolation points allowing for the use of regular upsampled quadrature.}
\revTwo{\citeauthor{Marple2016} \cite{Marple2016} simulated vesicles in periodic channel flows of arbitrary shape, where periodicity was imposed through an extra linear condition. An advantage of this approach is that it allows for the use of already existing fast solvers for the free space Green's functions. To handle vesicles in close proximity, a globally compensated trapezoidal rule was used. However, the close evaluation scheme was applied only to vesicle-vesicle interactions and not including the channel walls.}

In this paper, the fluid flow problem is considered in a periodic domain. \myChange{As a consequence,} all periodic images need to be considered when evaluating the fluid velocity. To make this \myChange{approach}  computationally viable, a fast method is needed. For periodic systems, fast methods utilising an Ewald summation approach \cite{Ewald1921} are especially suitable. With this approach, the periodic sums are split into two parts: one real space sum and one Fourier space sum. The computation of the Fourier space sum is then accelerated using FFTs. In 3D, Ewald decompositions using different ``screening functions" have been performed, see e.g. \cite{Beenakker1986,Hasimoto1959,Lindbo2010}.
In 2D, \citeauthor{VanDeVorst1996} \cite{VanDeVorst1996} derived the formulation to split the Stokeslet and the stresslet to compute the flowfield of a fluid in a domain with pores in the Stokes regime. However, this specific derivation yielded a non-symmetric expression for the stresslet decomposition, similar to that obtained by \cite{Marin2012b} in an alternative derivation for the three-dimensional case. The computations of the decomposed expressions can be sped up using the spectral Ewald method \cite{AfKlinteberg2014,Lindbo2011a}. When considering instead a problem without periodicity, either the spectral Ewald method for free-space \cite{AfKlinteberg2017a} or a Fast Multipole Method (FMM) can be utilised. The method in this paper can easily be modified to the free-space case, similar to that in \cite{Palsson2018}.

This paper presents a \myChange{new,} highly accurate numerical method for the simulation of deforming droplets in a periodic two-dimensional Stokes flow in the presence of \myChange{solid} constrictions and channel \myChange{walls}. The method handles close interactions between drops as well as drops and solid objects without an increase in error. The drops may be clean or covered by insoluble surfactants. The method is an extension of that in \cite{Palsson2018}. In this paper, the tools described in that paper are extended upon to include also solid walls and stationary objects. The method in this paper is general, i.e. considers both drops and solids in any configuration, and makes no distinction between channel walls and other solid objects, \myChange{contrary to other existing methods}. Moreover, this paper also includes a new derivation of the split of both the Stokeslet and the stresslet, using the Hasimoto screening function. This gives a decomposition of the Stokeslet equal to that of \citeauthor{VanDeVorst1996} \cite{VanDeVorst1996}, but a symmetric expression for the split of the stresslet. \myChange{Furthermore, this paper
\myChange{derives estimates of} the truncation errors that arise in from the Ewald decomposition.}

The problem setting in this paper is limited to two dimensions. Certain physical effects will be lost through this simplification, however it has been noted that a substantial degree of physical relevance remains \cite{kropinski2001,Zhou1993}. Furthermore, the reduction in dimension allows for larger simulations with an increased number of close interactions, due to the substantially reduced computational cost. In terms of numerical methods, this work shows the advantage of boundary integral formulations for highly accurate treatment of interface dynamics and close interactions, and there is ongoing work to develop the same abilities in three dimensions, see e.g. \cite{Sorgentone2018} and the references therein.

The paper is organised as follows: in \S\ref{sec:form} the governing equations, nondimensionalisation and boundary integral formulation are introduced. The numerical method is described in \S\ref{sec:meth}, and in \S\ref{sec:ewald} the \myChange{spectral} Ewald method in two dimensions to handle periodicity is described, together with the decomposition of the Stokeslet and the stresslet and truncation error estimates to facilitate parameter selection. In \S\ref{sec:results} the capabilities of the method are demonstrated through numerical tests.
% ------------------------------------------------------------------------------

%% FORMULATION
% !TEX root = ms.tex

% ******************************************************************************
\section{Problem formulation}
\label{sec:form}
% ******************************************************************************
The equations and mathematical tools needed to simulate surfactant-covered drops in a wall-bounded flow with solid constrictions are described in this section. First, the Stokes equations which govern the flow of the problem and the convection-diffusion equation for the surfactant concentration are stated both in dimensional and nondimensional form. Then follows a description of how to reformulate the Stokes equations for deformable drops in the presence of stationary solid objects and walls into an integral equation. Finally, the periodic extension of the problem is described.
% ------------------------------------------------------------------------------

% ******************************************************************************
\subsection{Governing equations}
\begin{figure}[h!]
  \centering
  \includegraphics[width=0.3\textwidth]{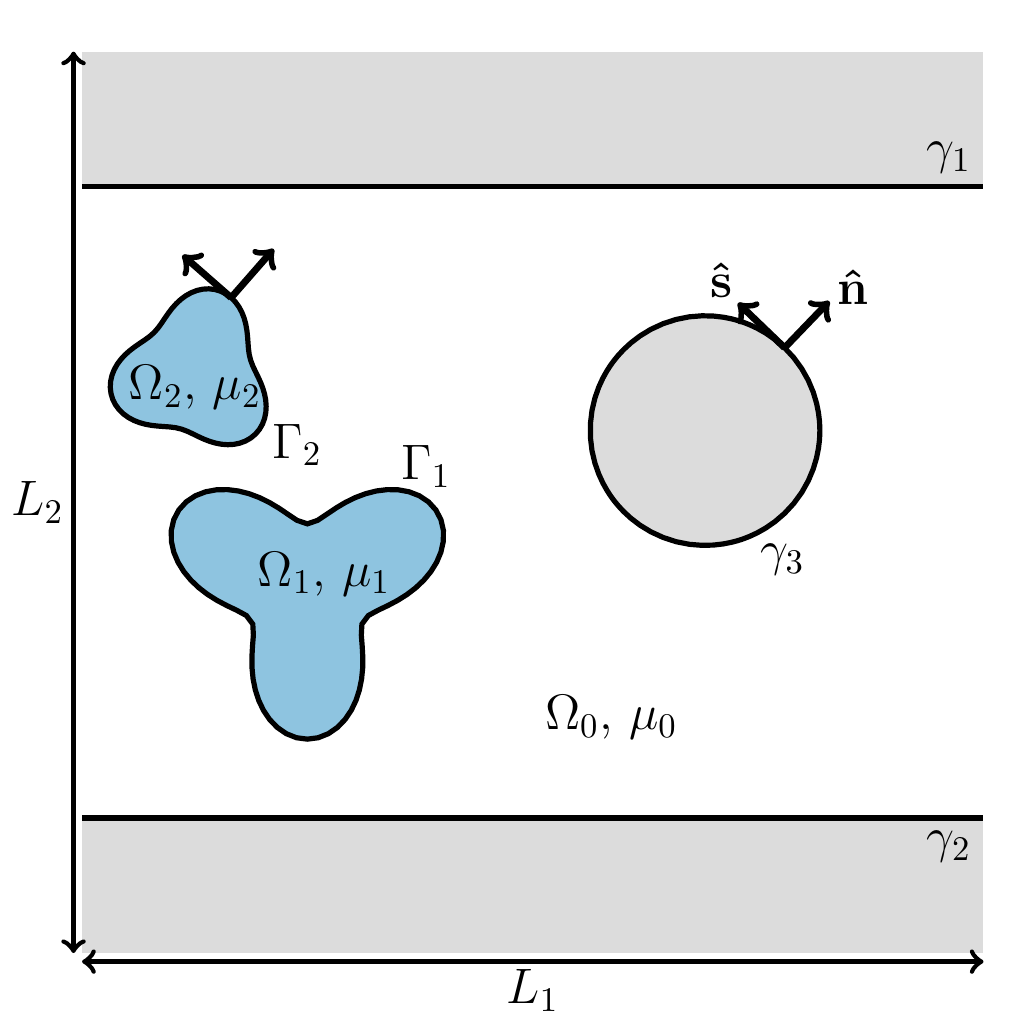}
  \caption{Example of domain configuration.}
  \label{fig:form_domain}
\end{figure}

% ******************************************************************************
The equations governing the physical problem are the incompressible Stokes equations, which in their dimensional form are
\begin{align}
\begin{split}
  \mu_0 \Delta \mb{u}_0 &= \nabla p_0, \; \nabla\cdot\mb{u}_0 = 0, \; \mb{x}\in\drop{0}, \\
  \mu_k \Delta \mb{u}_k &= \nabla p_k, \; \nabla\cdot\mb{u}_k = 0, \; \mb{x}\in\drop{k}, \; k=1,\hdots,\nDrop.
\end{split}
\label{eq:form_Dstokes}
\end{align}
Here, $\drop{0}$ is the bulk fluid surrounding the drops and solid objects and $\drop{k}$ is the interior of the drop $k$. There are in total $\nDrop$ drops. Furthermore, $p_k$ is the pressure and $\mb{u}_k$ the velocity, for $k=0,1,\hdots,\nDrop$. The drops and the bulk fluid are separated by the interfaces $\dinterf{k}$, on which the normal stress balance
\begin{align}
  -(p_0-p_k)\norm{k} + 2\left( \mu_0\mb{e}_0 -\mu_k\mb{e}_k\right)\cdot \norm{k} = \sigma_k\kappa_k\norm{k} - \nabla_s\sigma_k,
  \label{eq:form_Dinterf}
\end{align}
holds, where $\norm{k}$ is the outward-facing normal, $\mb{e}_k$ the strain tensor for the bulk and the interior of the drops, $\kappa_k=\nabla_s\cdot\norm{k}$ the curvature and $\nabla_s$ the surface gradient for $s$ traversing the drop $k$ in an anti-clockwise direction. Moreover, $\sigma_k:=\sigma_k(s,t)$ is the surface-tension coefficient of drop $k$ at time $t$. The fluid velocity is continuous on the drop boundaries, i.e. $\mb{u}_k=\mb{u}_0$ on $\dinterf{k}$. The interfaces are discretised anti-clockwise with a parameter $s\in[0,L_k(t)]$ where $L_k(t)$ is the length of the interface of drop $k$ at time $t$.
An example of a domain configuration can be found in Figure~\ref{fig:form_domain}. The drops translate and deform according to the ODE
\begin{align}
  \dfrac{\dd\mb{x}_k}{\dd t} = \mb{u}_k(\mb{x}_k,\sigma_k,t),
  \label{eq:form_dxdt}
\end{align}
for all points $\mb{x}_k\in\dinterf{k}$.
The boundaries of the solids are denoted by $\sinterf{k}$, for all solids $k=1,\hdots,\nSolid$. All solid boundaries have a no-slip boundary condition, i.e. the fluid velocity relative to the solid boundaries is always zero. The flow problem to be considered is the case where deformable drops are moving in channels or close to solid objects in an added flow field $\mb{u}_\infty$.

In addition, insoluble surfactants are considered. Their concentration is described by $\rho_k(s,t)$ and governed by a convection-diffusion equation on each interface; for $\mb{x}_k(s,t)\in\dinterf{k}$,
\begin{align}
  \dfrac{\DD\rho_k}{\DD t} - \rho_k\left(\nabla_s\cdot\mb{u}_k\right) = D_{\dinterf{}}\nabla_s^2\rho_k,
  \label{eq:form_Dsurf}
\end{align}
where $\frac{\DD}{\DD t}$ is the material derivative and $D_{\dinterf{}}$ is the diffusion coefficient along the interface \cite{stone}. As the surfactants are insoluble, the mass of surfactants is conserved along each interface,
\begin{align}
  \dfrac{\dd}{\dd t}\int\limits_{\dinterf{k}(t)} \rho_k(t)\,\dd S = 0,
\end{align}
for each drop $k=1,\hdots,\nDrop$. The surfactant concentration and the surface-tension coefficient are coupled through an equation of state. Here a linear equation of state is considered \myChange{\cite{Pawar1996}},
\begin{align}
  \sigma_k(s,t) = \sigma_0 - RT\rho_k(s,t),
  \label{eq:form_Deqstate}
\end{align}
for each drop $k$, where $\sigma_0$ is the surface-tension coefficient of a clean drop, $R$ is the universal gas constant and $T$ the temperature. This equation of state can be trivially exchanged to others.
% ------------------------------------------------------------------------------

% ******************************************************************************
\subsection{Nondimensionalisation}
% ******************************************************************************
All lengths are nondimensionalised using a characteristic length $r_0$, which unless otherwise stated is defined as the radius of the largest drop. The velocity is nondimensionalised by a characteristic velocity $U$ which is chosen from the imposed far-field flow as $U=\max(|\mb{u}_\infty(\mb{x})|)$. Furthermore, the surface-tension coefficient is nondimensionalised by the surface-tension coefficient of a clean drop, $\sigma_0$. \myChange{Consequently, the} characteristic pressure \myChange{is} $\frac{\mu U}{r_0}$ and \myChange{the} characteristic time \myChange{is} $\frac{r_0}{U}$. Also, the surfactant concentration is nondimensionalised by the initial surfactant concentration on the largest drop of the problem.

For the rest of this paper, all quantities are considered in their nondimensional form. The Stokes equations \eqref{eq:form_Dstokes} then read
\begin{align}
  \begin{split}
    \Delta \mb{u}_0 &= \nabla p_0, \; \nabla\cdot\mb{u}_0 = 0, \; \mb{x}\in\Omega_0, \\
    \lambda_k \Delta \mb{u}_k &= \nabla p_k, \; \nabla\cdot\mb{u}_k = 0, \; \mb{x}\in\Omega_k, \; k=1,\hdots,\nDrop,
  \end{split}
  \label{eq:form_stokes}
\end{align}
where $\lambda_k:=\frac{\mu_k}{\mu_0}$ is the viscosity ratio between the fluid of drop $k$ and the bulk. An inviscid bubble corresponds to the limit where $\lambda_k=0$. The no-slip condition on the solid boundaries is $\mb{u}_k = 0$ for all $\mb{x}\in\sinterf{k}$, $k=1,\hdots,\nSolid$. The stress balance over each interface $\dinterf{k}$ \eqref{eq:form_Dinterf} is rewritten as
\begin{align}
  -(p_0-p_k)\norm{k} + 2\left(\mb{e}_0 -\lambda_k\mb{e}_k\right)\cdot \norm{k} = \dfrac{1}{\text{Ca}} \left( \sigma_k\kappa_k\norm{k} - \nabla_s\sigma_k\right),
  \label{eq:form_interf}
\end{align}
where the Capillary number $\text{Ca}$ is defined as $\text{Ca}=\frac{U\mu_0}{\sigma_0}$.

The convection-diffusion equation governing the surfactant concentration \eqref{eq:form_Dsurf} becomes
\begin{align}
  \dfrac{\DD \rho_k}{\DD t} - \rho_k\left(\nabla_s\cdot\mb{u}_k\right) = \dfrac{1}{\text{Pe}_{\dinterf{}}}\nabla_s^2\rho_k, \; \mb{x}(s,t)\in\dinterf{k}, \; k=1,\hdots,\nDrop,
  \label{eq:form_surf}
\end{align}
where $\text{Pe}_{\dinterf{}} = \frac{r_0\mu_0}{\sigma_0 D_{\dinterf{}}}$ is the Peclet number. Furthermore, \eqref{eq:form_Deqstate} becomes
\begin{align}
  \sigma_k(s,t) = 1 - E\rho_k(s,t),
  \label{eq:form_eqstate}
\end{align}
where $E=\frac{RT\rho_0}{\sigma_0}$ is the so-called elasticity number.
% ------------------------------------------------------------------------------

% ******************************************************************************
\subsection{Boundary integral formulation}
% ******************************************************************************
A thorough derivation of the formulation for drops and solid particles in 3D can be found in \cite{Zinchenko2006}. Here, the same approach is followed but the formulation is rewritten for the 2D case.

For any point $\mb{x}$ in the bulk fluid $\Omega_0$, the velocity can be written as
\begin{align}
  \mb{u}(\mb{x}) = \sum_{k=1}^{\nDrop}\SL{\dinterf{k}}{f_k}{x}+ \mb{u}_\infty(\mb{x}) + \sum_{k=1}^{\nDrop}(\lambda_k-1) \DL{\dinterf{k}}{u}{x} + \SPC{x},
  \label{eq:form_ubulk}
\end{align}
where $\mb{f_k} :=\frac{1}{\text{Ca}}\left(\sigma_k\kappa_k\norm{k}-\nabla_s\sigma_k\right)$ from \eqref{eq:form_Dinterf}, $\SLe$ and $\DLe$ stand for the single-layer and double-layer contributions respectively and $\SPC{}$ stands for the solid-particle contribution as discussed below. The single- and double-layer potentials are defined as
\begin{align}
  \SL{\Lambda}{g}{x} &= \int\limits_{\Lambda} \mb{g}(\mb{y})\cdot\mb{G}(\mb{x}-\mb{y})\,\dd S_y. \label{eq:form_singleL} \\
  \DL{\Lambda}{g}{x} &= \int\limits_{\Lambda} \mb{g}(\mb{y})\cdot\mb{T}(\mb{x}-\mb{y})\cdot\norm{}(\mb{y})\,\dd S_y, \label{eq:form_doubleL}
\end{align}
for a drop interface or solid boundary, where $\Lambda=\dinterf{k}$ for a drop $k$ or $\Lambda=\sinterf{k}$ for a solid $k$. Here, $\mb{G}$ and $\mb{T}$ are the Stokeslet and stresslet respectively, which in 2D are defined as
\begin{align}
  G_{jl}(\mb{r}) &= -\dfrac{1}{4\pi}\left(-\delta_{jl}\log(|\mb{r}|)+\dfrac{\mb{r}_j\mb{r}_l}{|\mb{r}|^2}\right), \label{eq:form_stokeslet} \\
  T_{jlm}(\mb{r}) &= -\dfrac{1}{4\pi}\left(-4\dfrac{\mb{r}_j\mb{r}_l\mb{r}_m}{|\mb{r}|^4}\right), \label{eq:form_stresslet}
\end{align}
where $\mb{r}=\mb{x}-\mb{y}$. The solid-particle contribution can be defined as a single-layer potential over the solid boundaries, however, this generates an ill-conditioned system. In order to obtain a well-conditioned system the same approach as by \citeauthor{Zinchenko2006} \cite{Zinchenko2006} is taken, where a Hebeker representation \cite{Hebeker1986} is used to represent the solid-particle contribution. In this representation, the flow exterior to the solid particles is represented as a combination of single- and double-layer potentials,
\begin{align}
  \SPC{x} = \sum_{k=1}^{\nSolid} \SPCk{q}{x} \coloneqq \sum_{k=1}^{\nSolid} 2\DL{\sinterf{k}}{q}{x} + \eta \SL{\sinterf{k}}{q}{x}, \label{eq:form_hebeker}
\end{align}
where $\mb{q}$ is the so-called Hebeker density and $\eta$ is a proportionality factor which is set to $\eta=1$. Another option to achieve a well-conditioned system would be to use a completion flow as in \cite{Power1987}.

Finally, taking the limit of \eqref{eq:form_ubulk} as the point $\mb{x}\rightarrow \dinterf{\ell}$ and $\mb{x}\rightarrow \sinterf{\ell}$, gives a system of Fredholm integral equations of the second-kind \cite{Zinchenko2006},
\begin{align}
  \begin{split}
  \mb{u}(\mb{x}) - 2\sum_{k=1}^{\nDrop}\left(\dfrac{\lambda_k-1}{\myChange{\lambda_\ell}+1}\right) \DL{\dinterf{k}}{u}{x} &- \sum_{k=1}^{\nSolid}\dfrac{2}{\myChange{\lambda_\ell}+1}\SPCk{q}{x} =
    \sum_{k=1}^{\nDrop} \dfrac{2}{\myChange{\lambda_\ell}+1}\SL{\dinterf{k}}{f_k}{x} + \dfrac{2}{\lambda_\ell+1}\mb{u}_\infty(\mb{x}),
 \end{split}
 \label{eq:form_fullu}
\end{align}
where $\mb{x}\in \dinterf{\ell}$ for all fluid interfaces $\dinterf{\ell}\in\dinterf{}:=\bigcup\limits_{k=1}^{\nDrop}\dinterf{k}$ and
\begin{align}
  \begin{split}
  \mb{q}(\mb{x}) - \sum_{k=1}^{\nDrop}(\lambda_k-1)\DL{\dinterf{k}}{u}{x} &- \sum_{k=1}^{\nSolid}\SPCk{q}{x} = \sum_{k=1}^{\nDrop}\SL{\dinterf{k}}{f_k}{x} + \mb{u}_\infty(\mb{x}),
\end{split}
 \label{eq:form_fullq}
\end{align}
for $\mb{x}\in \sinterf{}:= \bigcup\limits_{k=1}^{\nSolid} \sinterf{k}$ (all solid boundaries). This is a system of Fredholm integral equations of the second kind which needs to be solved for each time step to obtain the velocity $\mathbf{u}$ with which the drops are moving.
% ------------------------------------------------------------------------------

% ******************************************************************************
\subsection{Periodicity}
\label{sec:form_per}
% ******************************************************************************
In this paper, the flow problem is considered with periodic boundary conditions in both the $x$- and $y$-direction. When computing the flow $\mb{u}$ at any point $\mb{x}\in \drop{0} \cup \sinterf{} \cup \dinterf{}$ through \eqref{eq:form_fullu} and \eqref{eq:form_fullq}, this means that
$\DL{\Lambda}{u}{x}$ and
 $\SL{\Lambda}{u}{x}$ contain the integrals over all periodic images over surfaces $\Lambda$. With their periodic replicas, they become
  \begin{align}
      \SLp{\Lambda}{g}{x} &= \sum\limits_{\mb{p}\in\mathbb{Z}^2}\int\limits_{\Lambda} \mb{g}(\mb{y})\cdot\mb{G}(\mb{x}-\mb{y}-\tau(\mb{p}))\,\dd S_y, \label{eq:SLp} \\
      \DLp{\Lambda}{g}{x} &= \sum\limits_{\mb{p}\in\mathbb{Z}^2}\int_\Lambda \mb{g}(\mb{y})\cdot\mb{T}(\mb{x}-\mb{y}-\tau(\mb{p}))\cdot\mb{n}(\mb{y})\,\dd S_y, \label{eq:DLp}
  \end{align}
  where $\tau(\mb{p}) = (p_1 L_1, p_2L_2)^T$ for a periodic box of size $L_1 \times L_2$ and $\mb{p}=(p_1,p_2)^T$, $p_1,\,p_2\in\mathbb{Z}$. \myChange{The interpretation and evaluation of these sums will be discussed in \S\ref{sec:ewald}.}

% ------------------------------------------------------------------------------

%% NUMERICAL METHODS
% !TEX root = ms.tex

% ******************************************************************************
\section{Numerical method}
\label{sec:meth}
% ******************************************************************************
To compute the evolution of deforming surfactant-covered drops, \eqref{eq:form_dxdt} and \eqref{eq:form_surf} need to be solved, generating the following system
\begin{align}
  \dfrac{\dd \mb{x}}{\dd t} &= \mb{u}_k(\mb{x},\sigma,t), \label{eq:meth_dxdt}\\
  \dfrac{\DD\rho_k}{\DD t} &- \rho_k\left(\nabla_s \cdot \mb{u}_k\right) = \dfrac{1}{\text{Pe}_{\dinterf{}}}\nabla_s^2\rho_k, \label{eq:meth_surfdt}
\end{align}
for all $\mb{x}\in\dinterf{k}$, for $k=1,\hdots,\nDrop$. The velocity $\mb{u}_k$ is determined by solving the boundary integral formulation described in \S\ref{sec:form}. Several components are needed to obtain an accurate solution to this system, most of them are described in detail in \cite{Palsson2018}. In this section, an overview of the method will be given.

Moving the drop interfaces by \eqref{eq:meth_dxdt} can result in a clustering of discretisation points on the interfaces. This is not ideal, as it necessitates remeshing. One can instead modify the velocity, as was done by \citeauthor{Hou1994} \cite{Hou1994} for elasticity problems and \citeauthor{kropinski2001} \cite{kropinski2001} for drops and bubbles. This approach will instead move the drops with velocity $\tilde{\mb{u}}_k$, i.e.
\begin{align*}
  \dfrac{\dd\mb{x}}{\dd t} = \tilde{\mb{u}}_k(\mb{x},\sigma,t), \; \mb{x}\in\dinterf{k},
\end{align*}
where the normal component of $\tilde{\mb{u}}_k$, $u_n$, is the same as for the fluid velocity $\mb{u}_k$ and the tangential velocity is modified as described in \S\ref{sec:meth_ut}. \myChange{It is possible to modify the tangential velocity since the normal velocity alone governs the deformation.}
Inserting \myChange{the} new velocity into \eqref{eq:meth_dxdt} and \eqref{eq:meth_surfdt} and expanding the material derivative gives the new system
\begin{align}
  \dfrac{\dd\mb{x}}{\dd t} &= \tilde{\mb{u}}_k(\mb{x},\sigma_k,t), \label{eq:meth_dzdt2} \\
  \dfrac{\partial \rho_k}{\partial t} &= \dfrac{\tilde{u}_t}{s_\alpha^k(t)}\dfrac{\partial\rho_k}{\partial \alpha} -\dfrac{1}{s_\alpha^k(t)}\dfrac{\partial(\rho_k u_t)}{\partial \alpha} - \rho_k u_n \kappa +\dfrac{1}{\text{Pe}_{\dinterf{}}
  s_\alpha^k(t)^2}\dfrac{\partial^2\rho_k}{\partial \alpha^2}, \; \alpha\in[0,2\pi], \label{eq:meth_surfdt2}
\end{align}
for all $\mb{x}\in\bigcup_{k=1}^{\nDrop}\dinterf{k}$. Here, $s_\alpha^k(\alpha,t) = \frac{1}{2\pi}L_k(t)$ where $L_k(t)$ is the length of drop interface $k$ at time $t$. \myChange{The droplet interfaces are parametrised with $\alpha$.}

\begin{figure}[h!]
\centering
\includegraphics[width=0.8\textwidth]{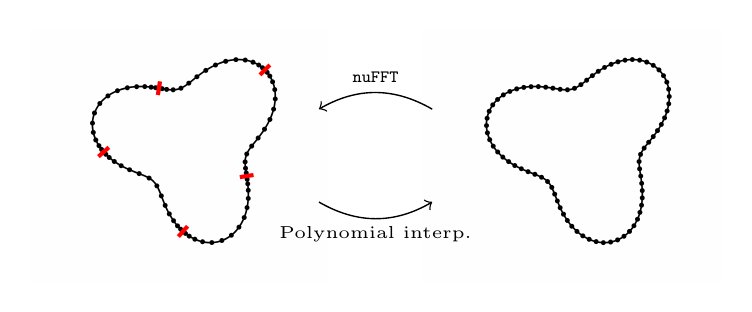}
\caption{Hybrid method using two discretisation for the drop interfaces: Grid 1 (left), composite Gauss-Legendre discretisation, and Grid 2 (right), uniform in arc length. Both grids have the same number of points, in this case $N=80$, i.e. $5$ panels. Red markers in the plot to the left mark panel divisions.}
\label{fig:meth_domainA}
\end{figure}

A hybrid method for discretising the equations in \eqref{eq:meth_dxdt} and \eqref{eq:meth_surfdt} is used.  The hybrid method consists of two discretisations: ``Grid 1'' is a panel-based composite 16-point Gauss-Legendre discretisation used on both drops and solids, and ``Grid 2'' is a uniform discretisation in arc-length, used only for the drops. Grid 1 will be used to determine the velocity from the boundary integral formulation in \eqref{eq:form_fullu}, $\mb{u}_k(\mb{x},\sigma,t)$, which is used to compute the normal component $u_n$ used to move the drops. Grid 2 will be used for determining the appropriate tangential velocity, for solving the surfactant equation \eqref{eq:meth_surfdt2} and for updating both
$\rho$ and $\mb{x}$ in time. To go from the equidistant discretisation to the Gauss-Legendre one a non-uniform \texttt{FFT} is used, see \citeauthor{Greengard2004} \cite{Greengard2004}. To go the opposite way, a 16-point polynomial interpolation on each panel is used. A schematic is shown in Figure~\ref{fig:meth_domainA}.

The problem setting is as follows; there are $\nDrop$ drops and $\nSolid$ solids. A solid $k$ is discretised by $\pSolid{k}$ Gauss-Legendre points (Grid 1), giving a total of discretisation points on the solids as $\pSolid{} := \sum_{k=1}^{\nSolid} \pSolid{k}$. A drop is discretised by $\pDrop{k}$ discretisation points uniform in arc-length (Grid 2), but also by  $\pDrop{k}$ Gauss-Legendre points (Grid 1). The total number of drop discretisation points is $\pDrop{} := \sum_{k=1}^{\nDrop} \pDrop{k}$ \myChange{which} makes the total number of discretisation points of the system $\ptot = \pSolid{} + \pDrop{}$,
and the number of unknowns $2\ptot$.
% ------------------------------------------------------------------------------

% ******************************************************************************
\subsection{Complex variable notation}
\label{sec:meth_complex}
% ******************************************************************************
When regarding the problem of deforming drops in 2D, it is beneficial to consider the formulation in complex variable notation where a point $\mb{x}$ corresponds to $z=x+iy$. Considering $z,\tau\in\mathbb{C}$, the complex counterparts of $\SL{\Lambda}{g}{x}$ and $\DL{\Lambda}{g}{x}$ are denoted $\SLc{\Lambda}{g}{z}$ and  $\DLc{\Lambda}{g}{z}$,
for $\Lambda=\dinterf{}$ (layer potential over drop interfaces) and $\Lambda=\sinterf{}$ (layer potential over solid boundaries). Furthermore, in complex notation
\begin{align}
  \SLc{\Lambda}{g}{z} = -\dfrac{1}{8\pi}\int\limits_\Lambda g(\tau)|\dd\tau| + \dfrac{1}{4\pi}\int\limits_\Lambda g(\tau)\log(|\tau-z|)|\dd\tau| - \dfrac{1}{8\pi}\int\limits_\Lambda \overline{g(\tau)}
  \dfrac{(\tau-z)}{(\overline{\tau}-\overline{z})}|\dd\tau|,
  \label{eq:form_singleC}
\end{align}
and
\begin{align}
  \DLc{\Lambda}{g}{z} = \dfrac{1}{2\pi}\int\limits_\Lambda g(\tau)\Im\left\{\dfrac{\dd\tau}{\tau-z}\right\} + \dfrac{1}{2\pi}\int\limits_\Lambda \overline{g(\tau)}\dfrac{\Im\left\{
  ( \overline{\tau}-\overline{z})\dd\tau\right\}}{(\overline{\tau}-\overline{z})^2}.
  \label{eq:form_doubleC}
\end{align}
For the remainder of this section, the problem will be treated in this complex setting.
% ------------------------------------------------------------------------------

% ******************************************************************************
\subsection{Computing the velocity $u(z,\sigma,t)$}
\label{sec:meth_u}
% ******************************************************************************
To compute the velocity the composite Gauss-Legendre discretisations (``Grid 1'') of drop and solid interfaces are used. The system to discretise is given by \eqref{eq:form_fullu} together with \eqref{eq:form_fullq}, using the periodic expressions for $\SLe$ and $\DLe$ in \eqref{eq:SLp} and \eqref{eq:DLp}. This system of Fredholm integral equations of the second kind will be solved by a Nystr\"{o}m method. In its discretised version \myChange{it} becomes, for
 $u_\ell \approx u(z_\ell)$ and $q_\ell \approx q(z_\ell)$ where $z_{\ell}$ are the Gauss-Legendre discretisation points with associated weights $w_\ell$, on either a drop interface or solid boundary,
\begin{align}
  \begin{split}
  u_\ell - & \sum_{k=1}^{\nDrop}2\left(\dfrac{\lambda_k-1}{\myChange{\lambda_\ell}+1}\right) \DLpc{\dinterf{k}}{u}{z_\ell} - \sum_{k=1}^{\nSolid}\dfrac{2}{\lambda_\ell+1}\left(2\DLpc{\sinterf{k}}{q}{z_\ell} +
   \SLpc{\sinterf{k}}{q}{z_\ell}\right) \\
        &= \sum_{k=1}^{\nDrop}\dfrac{2}{\myChange{\lambda_\ell}+1} \SLpc{\dinterf{k}}{f}{z_\ell} + \dfrac{2}{\lambda_\ell + 1}u_\infty(z_\ell), \; \forall \ell=1,\hdots,\pDrop,
\\
  q_\ell -& \sum_{k=1}^{\nDrop}(\lambda_k-1)\DLpc{\dinterf{k}}{u}{z_\ell} - \sum_{k=1}^{\nSolid}\left(2 \DLpc{\sinterf{k}}{q}{z_\ell} +  \SLpc{\sinterf{k}}{q}{z_\ell}\right) \\
  &= \sum_{k=1}^{\nDrop} \SLpc{\dinterf{k}}{f}{z_\ell} + u_\infty(z_\ell), \; \forall \ell=1,\hdots \pSolid.
\end{split}
\label{eq:meth_discSystem}
\end{align}
Here, $f_\ell:= \dfrac{1}{\text{Ca}}\left(\sigma_\ell \kappa_\ell \normC{\ell} - \nabla_s \sigma_\ell\right)$ from \eqref{eq:form_interf}, $\DLpc{\Lambda}{g}{z_\ell}$ and $\SLpc{\Lambda}{g}{z_\ell}$ are the discretised versions of $\DLc{\Lambda}{f}{z}$
in \eqref{eq:form_doubleC} and $\SLc{\Lambda}{f}{z}$ in \eqref{eq:form_singleC} respectively, including the periodic extension as mentioned in \S\ref{sec:form_per}. In \S\ref{sec:ewald} the computation of these periodic sums is described. Define $\tau(\mb{p})=p_1L_1+ip_2L_2$
 for $\mb{p}=(p_1,p_2)$ with $p_1,\,p_2\in\mathbb{Z}$
and $L_1$, $L_2$ the sides of the periodic box. The periodic double and single layer potentials are then defined as
\begin{align}
  \DLpc{\Lambda}{g}{z} &= \sum_{\mb{p}} \dfrac{1}{2\pi} \sum_{m=1}^{\pGen} \left[ g_m \mathbb{M}_{m}^{(1)}(z+\tau(\mb{p})) + \overline{g_m} \mathbb{M}_{m}^{(2)}(z+\tau(\mb{p})) \right], \text{ and }
   \label{eq:meth_DLcDisc} \\
  \SLpc{\Lambda}{g}{z} &=  \sum_{\mb{p}} \left\{ -\dfrac{1}{8\pi} \sum_{m=1}^{\pGen} \left[  w_m g_m|z_m^\prime|  + \overline{g_m} \mathbb{M}_{m}^{(3)}(z+\tau(\mb{p})) \right] + \dfrac{1}{4\pi}
   \mathop{\sum_{m=1}}_{z_m \neq z}^{\pGen} g_m \mathbb{M}_{m}^{(4)}(z+\tau(\mb{p})) \right\}, \label{eq:meth_SLcDisc}
\end{align}
for $\Lambda=\dinterf{},\,\sinterf{}$ and where (for $z\neq z_m$)
\begin{align}
  \begin{cases}
    \mathbb{M}_{m}^{(1)}(z) &= w_m \Im\left\{\dfrac{z_m^\prime}{z_m-z}\right\}, \\
    \mathbb{M}_{m}^{(2)}(z) &= w_m \dfrac{ \Im\left\{ z_m^\prime (\overline{z_m}- \overline{z})\right\}}{(\overline{z_m}-\overline{z})^2}, \\
    \mathbb{M}_{m}^{(3)}(z) &= w_m |z_m^\prime|\left(\dfrac{z_m-z}{\overline{z}_m-\overline{z}}\right), \\
    \mathbb{M}_{m}^{(4)}(z) &=  w_m |z_m^\prime| \log\left(|z_m-z|\right).
  \end{cases}
  \label{eq:meth_M}
\end{align}
Observe that the limits when $z_m=z$ are finite for $\mathbb{M}^{(1)}$, $\mathbb{M}^{(2)}$ and $\mathbb{M}^{(3)}$ and given by
\begin{align*}
  \mathbb{M}_{m}^{(1)}(z_m) = w_m \Im\left\{\dfrac{z_m^{\prime\prime}}{2z_m^\prime}\right\}, \quad \mathbb{M}_{m}^{(2)}(z_m) = w_m \dfrac{\Im\left\{ z_m^{\prime\prime}\overline{z_m^\prime}\right\}}{2(\overline{z_m^\prime})^2}, \quad \mathbb{M}_{m}^{(3)}(z_m) = w_m |z_m^\prime| \left(\dfrac{z_m^\prime}{\overline{z_m^\prime}}\right).
\end{align*}
\myChange{Again, the interpretation of both $\DLpc{\Lambda}{g}{z}$ and $\SLpc{\Lambda}{g}{z}$ is discussed in \S\ref{sec:ewald}.}
The second integral of \eqref{eq:form_singleC} corresponding to the integral of  $\mathbb{M}_{m}^{(4)}$ is more complicated as the kernel is non-smooth. Also, note that special quadrature is needed for the nearly-singular case when drops get close to each other or solid objects/walls. How to deal with both these issues is described in \S\ref{sec:meth_SQ}. In this subsection, it is assumed these problems can be dealt with efficiently and a highly accurate solution obtained for all points on all drops.

The discretised system in \eqref{eq:meth_discSystem} is then solved using \texttt{gmres}. The periodic sums are all computed with the spectral Ewald method, as described in \S\ref{sec:ewald}. The discretised system in \eqref{eq:meth_discSystem} has a unique solution by the Fredholm Alternative, and as it is a Fredholm integral equation of the second kind has spectral properties which enable \texttt{gmres} to converge in few iterations. The authors observe that the number of iterations of \texttt{gmres} vary with viscosity ratio, but have in none of the test simulations felt the need for a preconditioner.

% ******************************************************************************
\subsection{Special quadrature}
\label{sec:meth_SQ}
% ******************************************************************************
All integrals of $\mathbb{M}_{m}^{(k)}(z)$, $k\in[1,4]$, in \S\ref{sec:meth_u} become near-singular, i.e. when evaluating them at a point $z$ such that $\|z_m-z\| \ll 1$, for some $m$, yielding large numerical errors when the regular Gauss-Legendre quadrature rule is applied. This is the case for example when drops get close to each other or solids. How these errors behave and can be estimated was studied in \cite{Palsson2018}. To obtain accurate approximations of the integrals at any distance from the interfaces, a special quadrature scheme will be employed. The main idea of the specialised quadrature has been described in
\cite{Ojala2012,Ojala2015,Palsson2018}. This quadrature scheme will be employed for near-interactions in the case of integrals of
$\mathbb{M}_{m}^{(k)}$ for $k\in[1,4]$. In the case of $\mathbb{M}_{m}^{(4)}$ special treatment is needed also for the on-surface evaluations. For all kernels in \eqref{eq:meth_M}, the main idea of the special quadrature is similar. In short, for an integral of the form
$\int_\Gamma f(\tau)K(\tau,z)\dd\tau $, the idea is to express the function $f(\tau)$ as a polynomial in $\tau$, where the coefficients can be computed using a Vandermonde system. One can then use recursive formulas to compute the integrals as needed analytically. Following the notation in \cite{Ojala2012}, all the integrals of $\mathbb{M}_{m}^{(k)}$, $k\in[1,4]$ can be written on one of the following forms
\begin{align}
  I_1(z) = \int\limits_\Lambda \dfrac{h(\tau)}{\tau-z}\dd\tau, \quad I_2(z) = \int\limits_\Lambda \dfrac{h(\tau)}{(\tau-z)^2}\dd\tau, \quad I_3(z) = \int\limits_\Lambda h(\tau)\log(|\tau-z|)|\dd\tau|,
  \label{eq:meth_I123}
\end{align}
for any smooth boundary $\Lambda$. How to handle $I_1$ and $I_2$ is described in \cite{Ojala2015,Palsson2018}. The special quadrature to deal with $I_3$ can be found in \cite{Helsing2009,Ojala2012}. A brief overview of all three cases is given here. Note that the special quadrature treatment is strictly short-ranged, so there is no need to involve any periodicity in the calculations, except when considering drops and solids close to the edge of the periodic box. Furthermore, note that the third integral $I_3$ can be rewritten in the following way
\begin{align*}
  I_3(z) = \int\limits_\Lambda h(\tau)\log(|\tau-z|)|\dd\tau| = \Im\left\{ \int\limits_\Lambda \dfrac{h(\tau)}{\normC{\tau}} \log(\tau-z)\dd\tau \right\},
\end{align*}
where $\normC{\tau}$ is the normal of $\Lambda$ at point $\tau$.

For all three integrals, the approach is the same. Consider $\Lambda$ as a panel on either $\dinterf{}$ or $\sinterf{}$, rotated and scaled such that its endpoints are at $-1$ and $1$ in the complex plane. The evaluation point $z$ is rotated and scaled along with $\Lambda$. Now, expanding $h(\tau)$ as a monomial,  coefficients $c_\ell$ can be computed such that
\begin{align*}
  h(\tau) \approx \sum_{\ell=0}^{n-1} c_\ell \tau^\ell,
\end{align*}
where $n$ is the number of Gauss-Legendre points on each panel, here set to $n=16$. Inserting the interpolation into the integrals $I_1$, $I_2$ and $I_3$ the following expressions are obtained
\begin{align}
  I_1(z) \approx \sum_{\ell=0}^{n-1}c_\ell p_\ell, \quad I_2(z) \approx \sum_{\ell=0}^{n-1} c_\ell q_\ell \; \text{ and } \; I_3(z) \approx \sum_{\ell=0}^{n-1} d_\ell \Im\left\{ r_\ell\right\},
  \label{eq:meth_SQ123}
\end{align}
where
\begin{align}
  p_\ell := \int\limits_{-1}^1 \dfrac{\tau^\ell \dd\tau}{\tau-z}, \quad q_\ell := \int\limits_{-1}^1 \dfrac{\tau^\ell \dd\tau}{(\tau-z)^2} \; \text{ and } \; r_\ell = \int\limits_{-1}^1 \tau^\ell \log(\tau-z)\dd\tau.
  \label{eq:meth_pqr}
\end{align}
The analytical integrals $p_\ell$, $q_\ell$ and $r_\ell$ can be computed through recursion formulas, given in \ref{sec:app_sq}. The polynomial coefficients $d_\ell$ in \eqref{eq:meth_SQ123} are given by a polynomial expansion of the function $f(\tau)~=~h(\tau)/\normC{\tau}$, similar to that of $c_\ell$.

Regarding the on-surface evaluation of $I_3(z)$, the same approach is used. It can be seen as computing the quadrature weights for the log-kernel with a particular target point $z$. Using this approach, the quadrature weights can be precomputed and saved for all quadrature points on a panel in a $n \times n$ matrix. This can be  extended to a rectangular matrix to include special quadrature treatment also for target points on the neighbouring panels.
\\ \\
\myChange{The following \S\ref{sec:meth_ut}-\S\ref{sec:meth_Sapad} are kept brief and discussed in further detail in \cite{Palsson2018}.}
 % ******************************************************************************
 \subsection{Modifying the tangential velocity \myChange{on the droplet interfaces}}
 \label{sec:meth_ut}
 % ******************************************************************************
As previously stated, the tangential velocity needs to be modified in order to avoid clustering of discretisation points on the interfaces. Denote the velocity computed in \S\ref{sec:meth_u} by $u(z)$ and the modified velocity by $\tilde{u}(z)$. Both velocities are considered on the uniform grid (``Grid 2''). The normal components of the velocities are the same, i.e.
\begin{align*}
  \tilde{u}_n:=\Re\left\{\tilde{u}(s)\overline{\normC{}}(s)\right\} = \Re\left\{u(s)\overline{\normC{}}(s)\right\} =: u_n.
\end{align*}
The tangential component is modified according to (see \cite{Ojala2015})
 \begin{align}
   \tilde{u}_t(s) = \dfrac{s}{2\pi}\int\limits_0^{2\pi} \Im\left\{\dfrac{z^{\prime\prime}(q)}{z^\prime(q)}\right\}u_n(q)\dd q - \int\limits_0^s\Im\left\{ \dfrac{z^{\prime\prime}(q)}{z^\prime(q)}\right\}u_n(q)\dd q.
   \label{eq:meth_ut}
 \end{align}
 Note that the modified velocity $\tilde{u}(s) = \left[ u_n(s) + i\tilde{u}_t(s)\right]\normC{}(s)$ still fulfills the kinematic condition.

% ******************************************************************************
\subsection{Solving the surfactant equations}
\label{sec:meth_surf}
% ******************************************************************************
\noindent The equation for insoluble surfactants is solved as described in \cite{Palsson2018}. The equation to compute the surfactant concentration is given by \eqref{eq:meth_surfdt2}. This equation can be solved using a pseudo-spectral method which generates a system of ODEs to solve, one for each Fourier coefficient of $\rho_k(\alpha,t)$; $\widehat{\rho}^k_j(t)$:
\begin{align}
  \dfrac{\dd\widehat{\rho}^k_j}{\dd t} &= \widehat{\left(f_E^k\right)}_j + \widehat{\left(f_I^k\right)}_j, \; j=-\dfrac{\pDrop{k}}{2},\hdots\dfrac{\pDrop{k}}{2}-1,
\end{align}
for all drops $k$. Here, $f_I$ corresponds to the part of \eqref{eq:meth_surfdt2} that needs to be treated implicitly due to stiffness, and $f_E$ corresponds to everything else which is handled explicitly, i.e.
\begin{align*}
  f_E^k =  \dfrac{\tilde{u}_t}{s_\alpha^k(t)}\dfrac{\partial\rho_k}{\partial \alpha} -\dfrac{1}{s_\alpha^k(t)}\dfrac{\partial(\rho_k u_t)}{\partial \alpha} - \rho_k u_n \kappa_k, \text{ and } f_I^k = \dfrac{1}{\text{Pe}_{\dinterf{}}
  s_\alpha^k(t)^2}\dfrac{\partial^2\rho_k}{\partial \alpha^2}.
\end{align*}
Also, $\widehat{\left(f_E^k\right)}_j$ and $\widehat{\left(f_I^k\right)}_j$ correspond to the $j$th Fourier coefficient of $f_E^k$ and $f_I^k$ respectively. Since this is a pseudo-spectral method, the surfactant concentration is computed on the uniform grid (``Grid 2'').

% ******************************************************************************
\subsection{Adjusting the grid spacing \myChange{on the droplet interfaces}}
\label{sec:meth_Sapad}
% ******************************************************************************
\revOne{As the drop deforms over time the uniform discretisation will keep its points equidistant in arc-length due to the modified tangential velocity in \eqref{eq:meth_ut}. The spatial accuracy of the method depends on the grid spacing. As the droplet interface will stretch and contract over time, upsampling and downsampling is needed to keep the grid spacing similar in arc-length through the whole simulation. This approach keeps the resolution of the droplet interfaces approximately constant and reduces the cost of the simulation when droplets contract.}
\revOne{To facilitate the interchange of information between ``Grid 1'' and ``Grid 2'' as described  in Figure~\ref{fig:meth_domainA}, the number of discretisation points on each droplet is always kept as a multiple of 16. Using \texttt{FFT}s to upsample and downsample is trivial on the uniform grid.}
% As the drop deforms over time the uniform discretisation will keep its points equidistant with distance $\Delta s$ due to the modified tangential velocity in \eqref{eq:meth_ut}. \revOne{The spatial accuracy of the method depends on this distance. As the droplet interface will stretch and contract over time, the number of discretisation points for a desired level of accuracy can change. This can lead to either under- or over-resolved droplet boundaries as the droplets evolve in time.  To overcome this, a spatial adaptivity scheme to keep $\Delta s$ constant is utilised. The scheme boils down to regular upsampling/downsampling to keep $\Delta s$ approximately constant over time. To facilitate the interchange of information between ``Grid 1'' and ``Grid 2'', the number of discretisation points on each droplet is always kept as a multiple of 16. Using \texttt{FFT}s, this operation is trivial on the uniform grid. Each droplet interface is discretised with a uniform discretisation due to the pseudo-spectral method for solving the surfactant concentration. The adaptive scheme could be further improved by taking e.g. changes in curvature into account, or be expanded with the use of error estimates. The current scheme is however efficient and has proved adequate for most test cases.}

% ******************************************************************************
\subsection{Time-stepping scheme}
\label{sec:meth_time}
% ******************************************************************************
The coupled system to time-step is
\begin{align}
\begin{split}
  \dfrac{\dd z}{\dd t} &= \tilde{u}(z,\sigma,t), \; z\in\dinterf{} \\
  \dfrac{\dd\hat{\rho}^k_j}{\dd t} &= \widehat{\left(f_E^k\right)}_j + \widehat{\left(f_I^k\right)}_j, \; j=-\dfrac{\pDrop{k}}{2},\hdots\dfrac{\pDrop{k}}{2}-1, \; k=1,\hdots,\nDrop,
\end{split}
\end{align}
where $\widehat{\left(f_E^k\right)}_j$ should be treated explicitly and $\widehat{\left(f_I^k\right)}_j$ implicitly. In \cite{Palsson2018}, a second order time-stepping scheme was utilised. This time-stepper was chosen after a comparison of several schemes in \cite{Palsson2017}. It has for this paper been updated to a fourth order time-stepping scheme, which gives a considerable gain in computational cost. The scheme needs to handle adaptivity in time for both surfactant concentration $\rho$ and position $z$, and utilise the same stages for both equations. Therefore, a fourth order adaptive scheme by \citeauthor{Kennedy2003} \cite{Kennedy2003} is used, which uses the ``ARK4(3)6L[2]SA-ERK'' for the explicit parts together with the diagonally implicit ``ARK4(3)6L[2]SA-ESDIRK'' for the implicit part.
These are additive Runge-Kutta methods where adaptivity is acquired by comparing to a lower order scheme. For a Butcher tableau of the scheme, the reader is referred to Appendix C in \cite{Kennedy2003}. The time-step is modified using
\begin{align*}
  dt_{new} = \max\left( dt_{old}\left[s_f \dfrac{tol}{r}\right]^{1/4},\epsilon\right),
\end{align*}
where $s_f=0.8$ is a safety factor, $tol$ the given time-stepping tolerance, $\epsilon$ is machine epsilon and $r=\max(r_z,r_\rho)$, where $r_z$ is the measured error in $z$ and $r_\rho$ the measured error in $\rho$. \revOne{The measured errors in both $z$ and $\rho$ are computed by the comparison between the fourth and third order results in the time-stepping scheme, as the relative difference in two-norm.}

% ******************************************************************************
\subsection{Summary of the numerical method}
\label{sec:meth_sum}
% ******************************************************************************
Above, each step of the numerical method is described. Here follows an overview of how the different parts are put together.

Initially, all drop boundaries are discretised with a discretisation that is uniform in arc-length (``Grid 2''). The solid boundaries are discretised with a composite 16-point Gauss-Legendre scheme (``Grid 1''). The surfactant concentration is initialised on the uniform grid of the drop interfaces. Time-stepping is performed as described in \S\ref{sec:meth_time}. For every stage from $\tilde{t}$ to $\tilde{t}+c\,dt$ in the time-stepping scheme, the following steps are taken:
\begin{enumerate}
  \item The uniform drop discretisation is transformed to the panel-based G-L quadrature through a \texttt{nuFFT}, i.e. "Grid 2" $\rightarrow$ "Grid 1".
  \item The velocity $u$ for time $\tilde{t}$ is computed by solving the integral equation on both the drop interfaces and the solid boundaries, see \S\ref{sec:meth_u}. To compute the velocities for the periodic problem the spectral Ewald method as described in \S\ref{sec:ewald} is utilised. Special quadrature as in \S\ref{sec:meth_SQ} is used to obtain high accuracy for all discretisation points.
  \item Once the velocity $u$ for time $\tilde{t}$ is obtained, this is interpolated back to the uniform grid, "Grid 1" $\rightarrow$ "Grid 2", for the drop discretisation points.
  \item The velocity is then modified to $\tilde{u}$ using \eqref{eq:meth_ut}. Additionally, $f_E(\tilde{t})$ as in \S\ref{sec:meth_surf} is computed.
  \item The new position and surfactant concentration at time $\tilde{t}+c\,dt$ is computed using the method in \S\ref{sec:meth_time}.
  \item The surface-tension coefficient at time $\tilde{t}+c\,dt$ is computed through \eqref{eq:form_eqstate}.
\end{enumerate}
Furthermore, if an interface length has changed sufficiently, the number of discretisation points is modified through \texttt{FFT}s to keep $\Delta s$ constant, see \S\ref{sec:meth_Sapad}. Also, a Krasny filter of level $10^{-12}$ is applied to both position and surfactant concentration. \myChange{These two procedures are done in between complete Runge-Kutta time steps.}

% !TEX root = ms.tex

% ******************************************************************************
\section{Periodicity and \myChange{the spectral Ewald method}}
\label{sec:ewald}
% ******************************************************************************
As mentioned in \S\ref{sec:form_per} the flow problem is considered in a periodic setting in both $x$- and $y$-direction. The integrals to compute are \eqref{eq:SLp} and \eqref{eq:DLp} for the Stokeslet and stresslet respectively. They are discretised using the Gauss-Legendre quadrature described in \S\ref{sec:meth}, with quadrature nodes and weights $\mb{x}_n$, $w_n$, $n=1,\hdots,\pGen$, for $\pGen$ the total number of discretisation points on a boundary,
$\Lambda$. Thus, the approximations of the integrals \eqref{eq:SLp} and \eqref{eq:DLp} are
\begin{align}
  \SLp{\Lambda}{g}{x}_j &\approx \sum\limits_{\mb{p}\in\mathbb{Z}^2} \sum\limits_{n=1}^{\pGen} w_n G_{jl}(\mb{x}-\mb{x}_n-\tau(\mb{p}))g_l(\mb{x}_n) \eqqcolon u^G_j(\mb{x}),
  \label{eq:ewald_UG} \\
  \DLp{\Lambda}{g}{x}_j &\approx \sum\limits_{\mb{p}\in\mathbb{Z}^2} \sum\limits_{n=1}^{\pGen} w_n T_{jlm}(\mb{x}-\mb{x}_n-\tau(\mb{p}))g_l(\mb{x}_n)\normC{m}(\mb{x}_n) \eqqcolon u^T_j(\mb{x}).
   \label{eq:ewald_uT}
\end{align}
As previously in \S\ref{sec:form_per}, $\tau(\mb{p}) = (p_1 L_1, p_2L_2)^T$ for a periodic box of size $L_1 \times L_2$ and $\mb{p}=(p_1,p_2)^T$. The Stokeslet $G_{jl}$ and the stresslet $T_{jlm}$ are defined as in \eqref{eq:form_stokeslet}
 and \eqref{eq:form_stresslet} respectively.
 \myChange{As written, the sums \eqref{eq:ewald_UG} and \eqref{eq:ewald_uT} are not well-defined, since they are not convergent. However, $u^G$ can be made sense of as the velocity due to the point forces at $\mb{x}_n$, $n=1,\hdots,M^\Lambda$, repeated periodically with strengths $w_n g_l(\mb{x}_n)$. A pressure gradient will here be assumed to balance the force acting on the fluid, allowing $u^G$ to be expressed as a (slowly) converging sum in Fourier space \cite{VanDeVorst1996}. Similarly, a well-defined but slowly converging Fourier sum can be formulated for $u^T$.}
 Ewald decomposition \cite{Ewald1921} is used to remedy this slow decay. Each sum is split into two parts: one which contains the singularity and converges rapidly, referred to as the \rspaces, and one which contains a smooth periodic function, and thus converges quickly in Fourier space; the \kspaces. \myChange{The computations are accelerated using the spectral Ewald method \cite{AfKlinteberg2014,Lindbo2010}.} In this section, the Ewald decompositions of the 2D Stokeslet and stresslet are presented, together with estimates of the truncation errors and an overview of the \myChange{spectral} Ewald method.

% ******************************************************************************
\subsection{Ewald decomposition}
% ******************************************************************************
To illustrate the idea of an Ewald decomposition, the split into \rspace and \kspace ~is first computed for the Green's function for the biharmonic equation. In 2D, this Green's function has the following form
\begin{align}
  \biharm(|\mb{x}-\mb{y}|) = -\dfrac{|\mb{x}-\mb{y}|^2}{8\pi}\left( \log(|\mb{x}-\mb{y}|) - \alpha \right),
  \label{eq:ewald_B}
\end{align}
and it is the fundamental solution to $-\Delta^2 \biharm(|\mb{x}-\mb{y}|) = \delta(|\mb{x}-\mb{y}|)$. The  choice of constant $\alpha$ is free, and is here chosen to $\alpha=\frac{3}{2}$ preserve the 3D relation between the Green's function and the Stokeslet, i.e. $G_{jl}(\mb{r}) = \left(\Delta \delta_{jl} - \nabla_j \nabla_l\right)\biharm(|\mb{r}|)$ \cite{AfKlinteberg2017a}.

When considering a sum
\begin{align}
  u^{\biharm}(\mb{x}) = \sum\limits_{\mb{p}\in\mathbb{Z}^2}^{*}\sum\limits_{n=1}^{\pGen} \biharm(|\mb{x}-\mb{x}_n-\tau(\mb{p})|)f(\mb{x}_n),
  \label{eq:ub}
\end{align}
where the asterisk in the first sum corresponds to the exclusion of the term $\mb{x}-\mb{x}_n-\tau(\mb{p})=0$, the aim is to find a split into $\biharm^R(\mb{r},\xi)$ and $\widehat{\biharm}^F(\mb{k},\xi)$ such that
\begin{align}
  u^{\biharm}(\mb{x}) = \sum\limits_{\mb{p}\in\mathbb{Z}^2}^{*}\sum\limits_{n=1}^{\pGen} \biharm^R(|\mb{x}-\mb{x}_n-\tau(\mb{p})|,\xi)f(\mb{x}_n)
  + \dfrac{1}{V}\sum\limits_{\mb{k}\neq
   0}\widehat{\biharm}^F(|\mb{k}|,\xi)\sum\limits_{n=1}^{\pGen} f(\mb{x}_n)e^{-i\mb{k}\cdot(\mb{x}-\mb{x}_n)},
   \label{eq:ub2}
\end{align}
where $V=L_1 L_2$. Here, the first sum corresponds to the \rspaces, and the second to the \kspaces. Since the term where $\mb{x}-\mb{x}_n-\tau(\mb{p})=0$ is excluded from the sum $u^{\biharm}$, it should also be excluded from the \kspaces~when evaluating at a target point $\mb{x} = \mb{x}_n$, where $\mb{x}_n$ is any of the source points. This can be done by adding the limit
\[
\lim\limits_{|\mb{r}|\rightarrow 0} \left( \biharm^R(|\mb{r}|,\xi)-\biharm(|\mb{r}|)\right)f(\mb{x}_n),
\]
to the expression above.
\myChange{The sum $u^{\biharm}$ in \eqref{eq:ub} can be considered as a converging sum in Fourier space in a similar manner to that in \cite{VanDeVorst1996}. For the periodic case, the solution is unique up to a constant, which here will be determined by assuming that $\sum_{n=1}^{M^\Lambda}f(\mb{x}_n)=0$. This makes the split in \eqref{eq:ub2} independent of the splitting parameter $\xi$.}

The split is obtained by convolving $\biharm$ with a screening function $\gamma(r,\xi)$, as follows:
\begin{align*}
  \biharm(|\mb{r}|) = \underbrace{\biharm(|\mb{r}|)-\biharm(|\mb{r}|)\ast \gamma(r,\xi)}_{\eqqcolon \biharm^R(|\mb{r}|,\xi)} + \underbrace{\biharm(|\mb{r}|)\ast\gamma(r,\xi)}_{\eqqcolon \biharm^F(|r|,\xi)}.
\end{align*}
The screening function should be defined such that $\biharm^R(|\mb{r}|,\xi)$ is short-range and $\biharm^F(|\mb{r}|,\xi)$ is smooth and long-range. The Hasimoto screening function, defined as
\begin{align}
  \gamma(r,\xi) = \dfrac{\xi^2}{\pi}(2-\xi^2r^2)e^{-\xi^2r^2} \; \leftrightarrow \; \widehat{\gamma}(k,\xi) = \left( 1+\dfrac{k^2}{4\xi^2} \right)e^{-k^2/4\xi^2},
   \label{eq:ewald_hasimoto}
\end{align}
meets these criteria, where $r = |\mb{r}|$ and $k=|\mb{k}|$. In Fourier space, the convolution $\biharm(|\mb{r}|)\ast\gamma(r,\xi)$ is computed as a multiplication, which gives $\widehat{\biharm}^F(|\mb{k}|,\xi) = \widehat{\gamma}(k,\xi)\widehat{\biharm}(|\mb{k}|,\xi)$. For the biharmonic Green's function, $\widehat{\biharm}(|\mb{k}|)=-1/k^4$.
How to compute $\biharm^R(|\mb{r}|,\xi)$ by convolution is described in \ref{sec:app_ewald_B}. The split obtained is the following,
\begin{align}
  \begin{cases}
    \biharm^R(|\mb{r}|,\xi) &= \dfrac{1}{16\pi\xi^2}\left( \xi^2 r^2 E_1(\xi^2 r^2) -e^{-\xi^2r^2}\right) \\
    \widehat{\biharm}^F(|\mb{k}|,\xi) &= \dfrac{-1}{k^4}\left(1+\dfrac{k^2}{4\xi^2}\right)e^{-k^2/4\xi^2}.
  \end{cases}
  \label{eq:ewald_Bsplit}
\end{align}

% ******************************************************************************
\subsubsection{Stokeslet}
% ******************************************************************************
To find the split of the Stokeslet, first note that the Stokeslet can be expressed as an operator $K_{jl}$ acting on $\biharm(|\mb{r}|)$ \cite{AfKlinteberg2017a}, where $K_{jl} = \Delta \delta_{jl} - \nabla_j \nabla_l$. Thus, to compute the \rspace ~part of the Stokeslet, $G^R_{jl}$, this operator is applied to $\biharm^R$ which gives
\begin{align}
    G^R_{jl}(\mb{r},\xi) = K_{jl}\biharm^R(|\mb{r}|,\xi) = \dfrac{1}{4\pi}\left[ e^{-\xi^2 r^2}\left(\hat{r}_j \hat{r}_l-\delta_{jl}\right) + \dfrac{\delta_{jl}}{2}E_1(\xi^2 r^2)\right],
    \label{eq:ewald_GR}
\end{align}
where $\hat{r}_j = r_j/r$. Similarly, the Fourier space part can be expressed as $\widehat{G}^F_{jl}(\mb{k},\xi)=\widehat{K}_{jl}\widehat{\biharm}^F(|\mb{k}|,\xi)$, where $\widehat{K}_{jl}$ denotes the pre-factor that is produced when $K_{jl}$ is applied to $e^{-i\mb{k}\cdot\mb{r}}$, i.e. $\widehat{K}_{jl}=-\delta_{jl}k^2 + k_j k_l$. This gives
\begin{align}
  \widehat{G}^F_{jl}(\mb{k},\xi) = \widehat{K}_{jl}\widehat{\biharm}^F_{jl}(|\mb{k}|,\xi) = \left(\delta_{jl}-\hat{k}_j \hat{k}_l\right)\dfrac{1}{k^2}\left(1+\dfrac{k^2}{4\xi^2}\right)e^{-k^2/4\xi^2},
  \label{eq:ewald_GF}
\end{align}
where $\hat{k}_j=k_j/k$ and $k=|\mb{k}|$ for $\mb{k}=(k_1,k_2)$. Thus, $u_j^G(\mb{x})$ in \eqref{eq:ewald_UG} can be written as
\begin{align}
  \begin{split}
  u_j^G(\mb{x}) = \sum\limits_{\mb{p}\in\mathbb{Z}^2}^{*}\sum\limits_{n=1}^{\pGen}G^R_{jl}(\mb{x}-\mb{x}_n-\tau(\mb{p}),\xi)f_l(\mb{x}_n) + \dfrac{1}{V}\sum\limits_{\mb{k}\neq 0} \widehat{G}^F_{jl}(\mb{k},\xi)\sum\limits_{n=1}^{\pGen}f_l(\mb{x}_n)e^{-i\mb{k}\cdot(\mb{x}-\mb{x}_n)}.
 \end{split}
 \label{eq:ewald_Gsplit}
\end{align}
In the case of the target point $\mb{x}=\mb{x}_n$ for any source point $\mb{x}_n$, the self-contribution that arises from the Fourier part needs to be removed. Thus,
\begin{align*}
    \lim\limits_{|\mb{r}|\rightarrow 0} \left(G^R_{jl}(\mb{r},\xi)-G_{jl}(\mb{r})\right)f_l(\mb{x}_n) = -\left(\dfrac{1}{2}\gamma +\log(\xi)+1\right)\delta_{jl}f_l(\mb{x}_n),
\end{align*}
needs to be added to the expression, where $\gamma$ is the Euler-Mascheroni constant. \myChange{Note that no assumptions are made on $f(\mb{x}_n)$ for this split.}

% ******************************************************************************
\subsubsection{Stresslet}
% ******************************************************************************
The Ewald decomposition of the stresslet is computed in a similar manner to the Stokeslet, by applying an operator $K$ to the decomposition of the biharmonic Green's function. The operator relating the stresslet and $\biharm(\mb{r})$ is \cite{AfKlinteberg2017a},
\begin{align*}
  K_{jlm} = \left( \delta_{jl}\nabla_m + \delta_{lm}\nabla_j + \delta_{jm}\nabla_l\right)\Delta - 2\nabla_j\nabla_l\nabla_m.
\end{align*}
Thus, the \rspace~part of the stresslet is given by
\begin{align}
  T^R_{jlm}(\mb{r},\xi) = K_{jlm}\biharm^R(|\mb{r}|,\xi) = \dfrac{1}{4\pi}e^{-\xi^2 r^2}\left(\dfrac{-4\hat{r}_j\hat{r}_l\hat{r}_m}{r}(1+\xi^2 r^2) + \myChange{2\xi^2r}(\delta_{jl}\hat{r}_m + \delta_{jm}\hat{r}_l + \delta_{lm}\hat{r}_j)\right).
  \label{eq:ewald_TR}
\end{align}
By applying $\widehat{K}_{jlm} = -i\left[ \left( \delta_{jl}k_m + \delta_{jm}k_l +\delta_{lm}k_j\right)k^2 - 2k_j k_l k_m\right]$ to $\widehat{\biharm}^F(|\mb{k}|,\xi)$, the \kspace~part of the stresslet is computed as
\begin{align}
  \widehat{T}^F_{jlm}(\mb{k},\xi) = \widehat{K}_{jlm}\widehat{\biharm}^F(|\mb{k}|,\xi) = i\left[ \left( \delta_{jl}\hat{k}_m + \delta_{jm}\hat{k}_l +\delta_{lm}\hat{k}_j\right) - 2\hat{k}_j \hat{k}_l \hat{k}_m\right]\dfrac{1}{k}\left(1+\dfrac{k^2}{4\xi^2}\right)e^{-k^2/4\xi^2}.
  \label{eq:ewald_TF}
\end{align}
The complete expression thus reads
\begin{align}
  \begin{split}
  u_j^T(\mb{x}) =& \sum\limits_{\mb{p}\in\mathbb{Z}^2}^{*}\sum\limits_{n=1}^{\pGen} T^R_{jlm}(\mb{x}-\mb{x}_n-\tau(\mb{p}),\xi) f_l(\mb{x}_n) n_m(\mb{x}_n)  \\
  &+ \dfrac{1}{V}\sum\limits_{\mb{k}\neq 0} \widehat{T}^F_{jlm}(\mb{k},\xi)\sum_{n=1}^{\pGen} f_l(\mb{x}_n)n_m(\mb{x}_n)e^{-i\mb{k}\cdot(\mb{x}-\mb{x}_n)} + \dfrac{1}{V}\sum\limits_{n=1}^{\pGen} \hat{T}^{F,(0)}_{jlm}(\mb{x}_n)f_l(\mb{x}_n)n_m(\mb{x}_n),
\end{split}
  \label{eq:ewald_Tsplit}
\end{align}
\myChange{where no assumptions are made on $f(\mb{x}_n)$.} For the stresslet, there is no ``self-interaction term" as is the case for the Stokeslet, as $\lim\limits_{|\mb{r}|\rightarrow 0} T^R_{jlm}(\mb{r},\xi)-T_{jlm}(\mb{r})=0$. The term $\widehat{T}^{F,(0)}$ is chosen to guarantee zero-mean flow through the primary periodic cell, and also ensures that the stresslet identity is met \cite{AfKlinteberg2014}. This corresponds to setting
\begin{align*}
  \widehat{T}^{F,(0)}_{jlm}(\mb{y}) = \delta_{lm}y_j.
\end{align*}

As a side note regarding the stresslet, it can also be computed using derivatives of the Laplace Green's function, $\lapl(r)$, and the Stokeslet, i.e.
\begin{align*}
  T_{jlm}(\mb{x}-\mb{y}) = 2\dfrac{\partial \lapl(|\mb{x}-\mb{y})|}{\partial x_l}\delta_{jm} + \dfrac{\partial G_{jl}(\mb{x}-\mb{y})}{\partial x_m} + \dfrac{\partial G_{lm}(\mb{x}-\mb{y})}{\partial x_j}.
\end{align*}
Computing $T^R(\mb{r},\xi)$ and $\widehat{T}^F(\mb{k},\xi)$ through this relation using the Hasimoto screening of $\lapl(|\mb{r}|)$ into $\lapl^R(|\mb{r}|,\xi)$ and $\widehat{\lapl}^F(|\mb{k}|,\xi)$ and $G^R_{jl}(\mb{r},\xi)$,
$\widehat{G}^F_{jl}(\mb{k},\xi)$ from above, generates the same results as those in \eqref{eq:ewald_TR}, \eqref{eq:ewald_TF}. The split of $\lapl$ is derived in
\ref{sec:app_ewald_L}. Moreover, if instead using the so-called Ewald screening function to split $\lapl$,
\begin{align*}
  \gamma_E(r,\xi) = \dfrac{\xi^2}{\pi}e^{-\xi^2r^2} \; \leftrightarrow \; \widehat{\gamma}_E(k,\xi) = e^{-k^2/4\xi^2},
\end{align*}
while keeping the Hasimoto screening function for the Stokeslet, $G$, this generates the same expression as that found by
\citeauthor{VanDeVorst1996} \cite{VanDeVorst1996}. This expression, however, is not symmetric, and will not be used in this work.

% ******************************************************************************
\subsection{Truncation errors}
% ******************************************************************************
The sums in \eqref{eq:ewald_Gsplit} and \eqref{eq:ewald_Tsplit} converge fast, but cannot be computed numerically without truncation. To decide where to make the truncations, the \rspaces~and the \kspaces~need to be regarded separately and their truncation errors estimated. The error estimates and their derivations are described in detail in \ref{sec:app_trunc} and inspired by the error estimates in \cite{AfKlinteberg2017a}. Here, the case of a square box, i.e. $L_1=L_2=L$ is considered for simplicity.

\subsubsection{Truncation errors for the \rspaces}
\begin{figure}[h!]
  \centering
  \includegraphics[width=0.4\textwidth]{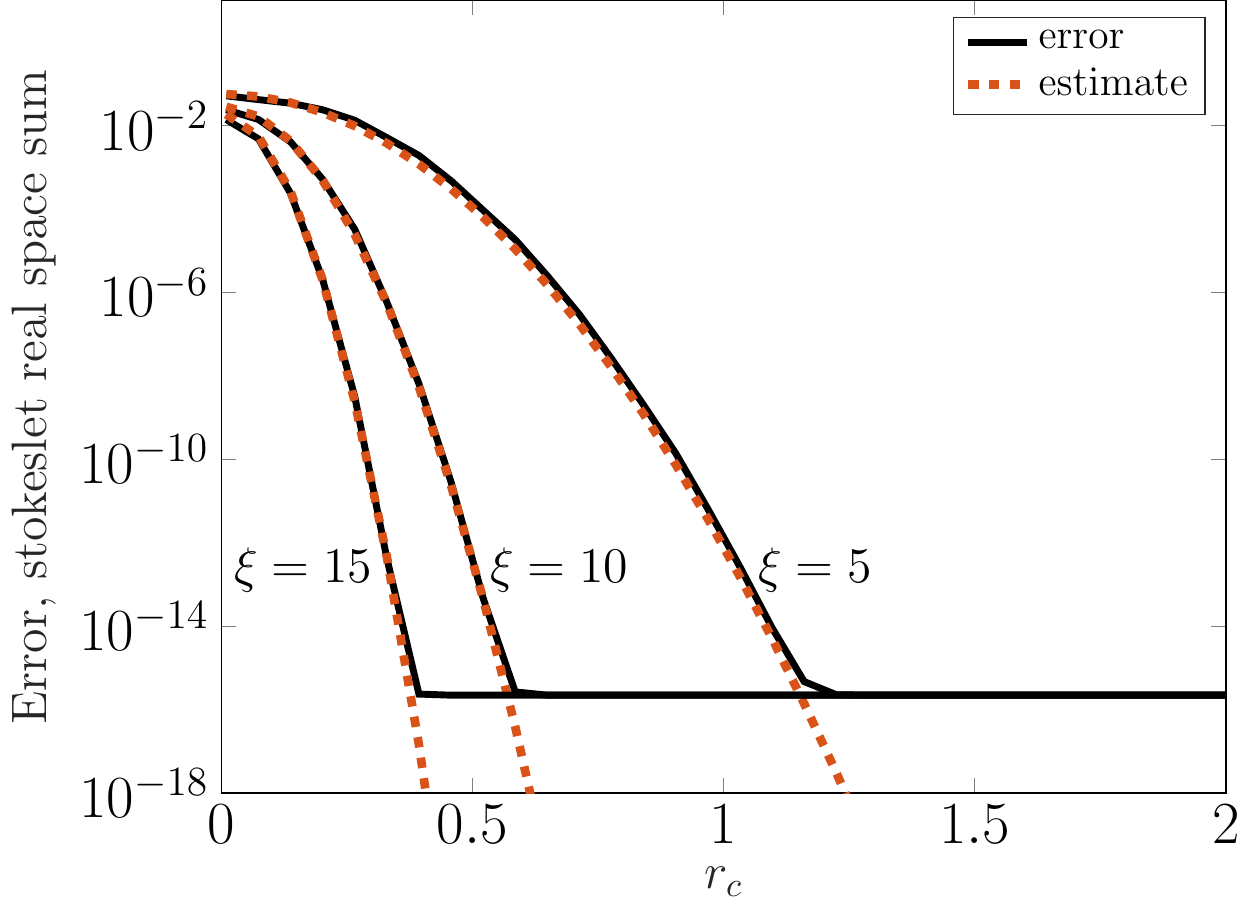}%
  \includegraphics[width=0.4\textwidth]{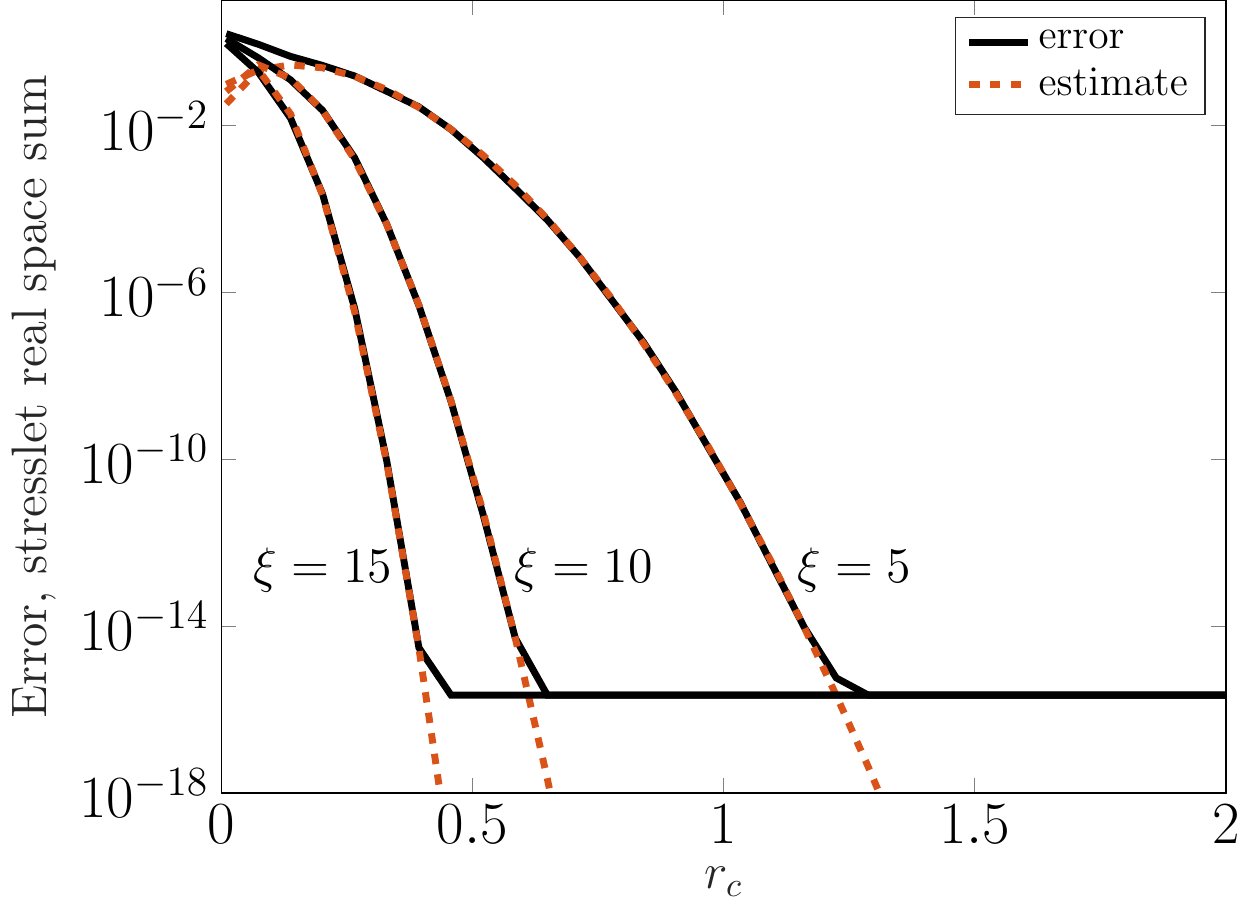}
  \caption{Truncation error estimate for real space sum for Stokeslet (left) and stresslet (right), as a function of $r_c$ for $\xi=5,10,15$. Black lines are the measured errors and dashed lines the estimates of \eqref{eq:ewald_GRest} and \eqref{eq:ewald_TRest}. The system is $N_s=10^3$ randomly distributed point sources within a square of size $L=2\pi$, and $N_t=10^2$ randomly distributed target points in the same square.}
  \label{fig:ewald_realest}
\end{figure}

\noindent The two \rspaces s for the Stokeslet and stresslet respectively are defined as
\begin{align}
  \begin{cases}
  u_j^{G,R}(\mb{x},\xi) &= \sum\limits_{\mb{p}\in\mathbb{Z}^2}^{*}\sum\limits_{n=1}^{\pGen}G^R_{jl}(\mb{x}-\mb{x}_n-\tau(\mb{p}),\xi)f_l(\mb{x}_n), \\
  u_j^{T,R}(\mb{x},\xi) &= \sum\limits_{\mb{p}\in\mathbb{Z}^2}^{*}\sum\limits_{n=1}^{\pGen} T^R_{jlm}(\mb{x}-\mb{x}_n-\tau(\mb{p}),\xi) f_l(\mb{x}_n) \normC{m}(\mb{x}_n)
\end{cases}
\label{eq:ewald_rspaces}
\end{align}
For these sums, only target points within a cut-off radius defined as $r_c$ will be considered. Comparing $u_j^{G,R}$ and $u_j^{T,R}$ with the truncated sums $\tilde{u}_j^{G,R}$ and $\tilde{u}_j^{T,R}$, the RMS of the truncation error is given by
\begin{align*}
  \partial \mb{u}^{G,R} = \sqrt{\dfrac{1}{\pGen}\sum_{n=1}^{\pGen}\left\vert\mb{u}^{G,R}-\tilde{\mb{u}}^{G,R}\right\vert^2}, \quad   \partial \mb{u}^{T,R} = \sqrt{\dfrac{1}{\pGen}\sum_{n=1}^{\pGen}\left\vert\mb{u}^{T,R}-\tilde{\mb{u}}^{T,R}\right\vert^2},
\end{align*}
for the Stokeslet and stresslet respectively. The truncation errors are estimated as
\begin{align}
  \left(\delta \mb{u}^{G,R}\right)^2 \approx \dfrac{Q_G\pi}{4L^2} \dfrac{e^{-2\xi^2 r_c^2}}{\xi^2},
  \label{eq:ewald_GRest}
\end{align}
for the Stokeslet, and
\begin{align}
  \left( \delta\mb{u}^{T,R}\right)^2 \approx \dfrac{2\pi Q_T}{L^2} \xi^2 r_c^2 e^{-2\xi^2 r_c^2},
  \label{eq:ewald_TRest}
\end{align}
for the stresslet, where $Q_G = \sum_{n=1}^{\pGen} \sum_{l=1}^2 f_l^2(\mb{x}_n)$ and $Q_T= \sum_{n=1}^{\pGen} \sum_{l,m=1}^2 f_l^2(\mb{x}_n)\norm{m}^2(\mb{x}_n)$. The derivations of these estimates can be found in \ref{sec:app_truncR}.
The truncation errors and estimates for both Stokeslet and stresslet are shown in Figure~\ref{fig:ewald_realest}. Comparing the \rspace~estimates to those empirically obtained for 3D in \cite{AfKlinteberg2014,Lindbo2010}, the estimates show the same asymptotic behaviour, but vary in the constant factor in front and powers of $(\xi r_c)$ needed.

\subsubsection{Truncation errors for the \kspaces}
\begin{figure}[h!]
  \centering
  \includegraphics[width=0.4\textwidth]{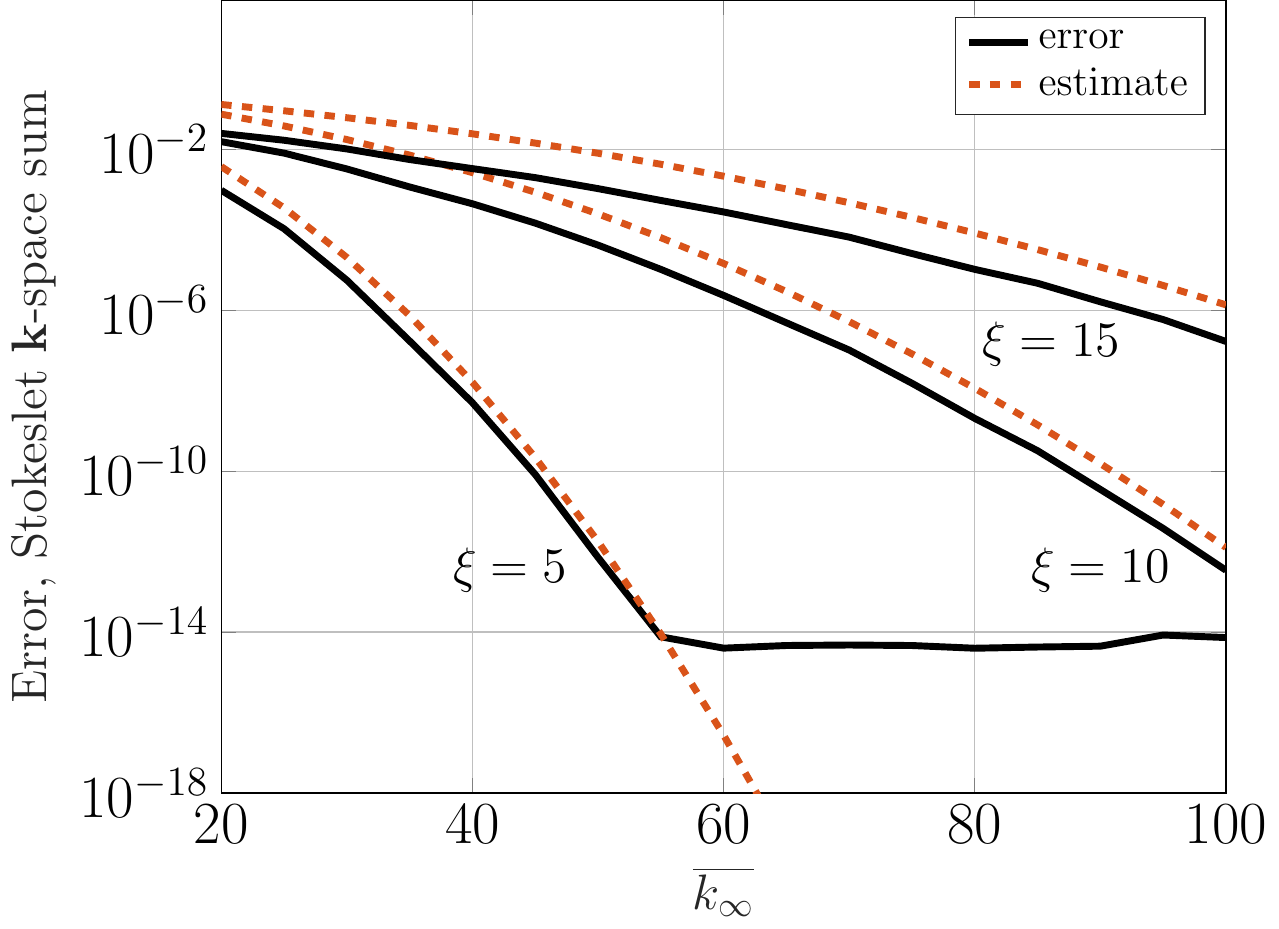}%
  \includegraphics[width=0.4\textwidth]{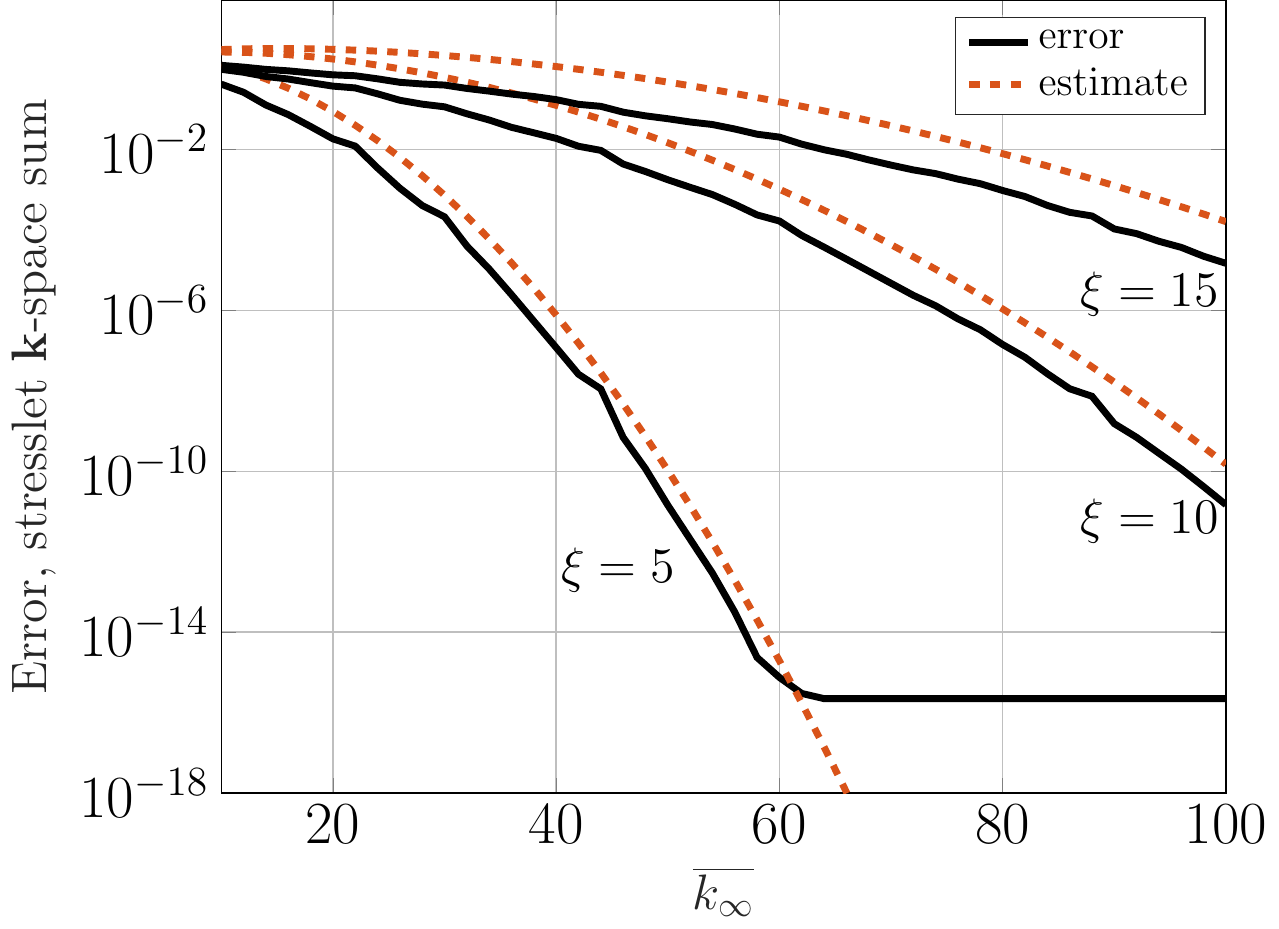}
  \caption{Truncation error estimate for the \kspaces~for the Stokeslet (left) and stresslet (right), as a function of $k_\infty$ for $\xi=5,10,15$. Black lines are the measured errors and the dashed lines the estimates using \eqref{eq:ewald_GFest} and \eqref{eq:ewald_TFest}. The system is the same as that in Figure~\ref{fig:ewald_realest}.}
  \label{fig:ewald_Fest}
\end{figure}

\noindent For the \kspace, the sums to compute are defined as
\begin{align*}
  \begin{cases}
      u_j^{G,F}(\mb{x},\xi)&= \dfrac{1}{V}\sum\limits_{\mb{k}\neq 0} \widehat{G}^F_{jl}(\mb{k},\xi)\sum\limits_{n=1}^{\pGen}f_l(\mb{x}_n)e^{-i\mb{k}\cdot(\mb{x}-\mb{x}_n)}, \\
      u_j^{T,F}(\mb{x},\xi) &= \dfrac{1}{V}\sum\limits_{\mb{k}\neq 0} \widehat{T}^F_{jlm}(\mb{k},\xi)\sum\limits_{n=1}^{\pGen}f_l(\mb{x}_n)\normC{m}(\mb{x}_n)e^{-i\mb{k}\cdot(\mb{x}-\mb{x}_n)}.
  \end{cases}
\end{align*}
and they are truncated in \kspace such that $\mb{k}=(k_1,k_2)$ for $k_1,k_2\in[-k_\infty,k_\infty]$. The RMS of the truncation error is estimated as
\begin{align}
  \left(\delta\mb{u}^{G,F}\right)^2 \approx \dfrac{4Q_G}{L^5\pi k_\infty} e^{-2k_\infty^2/4\xi^2},
  \label{eq:ewald_GFest}
\end{align}
for the Stokeslet and
\begin{align}
  \left(\delta \mb{u}^{T,F}\right)^2 \approx \dfrac{8\pi Q_T}{L^5}k_\infty e^{-2k_\infty^2/4\xi^2},
  \label{eq:ewald_TFest}
\end{align}
for the stresslet, with $Q_G$ and $Q_T$ as defined previously. These estimates are derived in \ref{sec:app_truncF}. The truncation errors together with the estimates are shown in Figure~\ref{fig:ewald_Fest} for both the Stokeslet and the stresslet. Additionally, comparing the \kspace~estimates to those for 3D in \cite{AfKlinteberg2014,Lindbo2010}, the estimates again show the same asymptotic behaviour, but the expressions in front of $e^{-k_\infty^2/4\xi^2}$ differ. The \kspace~truncation estimates are not as precise as their \rspace~counterparts, but always overestimate the errors.

% ******************************************************************************
\subsection{The \myChange{spectral} Ewald method}
\label{sec:semeth}
% ******************************************************************************
Using Ewald decompositions such as those in \eqref{eq:ewald_Gsplit} and \eqref{eq:ewald_Tsplit}, the sums to compute are now rapidly converging. For a system with $N$ discretisation points, computing the sums directly
results in an $\OO(N^2)$ complexity. To \myChange{increase computational efficiency}, the \myChange{spectral} Ewald method \cite{AfKlinteberg2014,AfKlinteberg2017a,Lindbo2010} is used, which makes the computations of the sums $\OO(N\log(N))$ in cost. The method is thoroughly described in the previous references for the three dimensional case, and generalises to 2D easily. In short, the key approach is to evaluate the Fourier space sums on a grid in the \kspaces~which enables the use of FFTs of the size $M^2$ to speed up the computations.
To spread source and target points to the grid truncated Gaussians with $P^2$ points support are used, and the shape parameter of the Gaussians is determined to minimise the approximation error for \myChange{the} given $P$. The real space sums can be computed in $\OO(N)$ time, by constricting the evaluation only to points in a near-neighbour list of each point $\mb{x}_t$, defined as
$NL_t = \left\{ (\mb{x}_s,\mb{p}) : |\mb{x}_t-\mb{x}_s-\tau(\mb{p})|<r_c\right\}$. Under the assumption of a  constant number of those near neighbours, to create such a list is also $\OO(N)$. The parameter $\xi$ from the screening function decide how much work is put into the \rspaces~ and \kspaces~ respectively.

There are several parameters to set in the method. To keep the number of nearest neighbours constant in the real-space sum as the system is scaled up, the cut-off radius $r_c$ is set first. \myChange{The computation of the real-space sums are thus $\OO(N)$ in cost.} Using the estimates in \eqref{eq:ewald_GRest} (Stokeslet) and \eqref{eq:ewald_TRest} (stresslet), for a given tolerance $tol_e$ the splitting parameter $\xi$ is computed. From $\xi$ the corresponding $k_\infty = \frac{M}{2}$ is computed from the estimates in \eqref{eq:ewald_GFest} (Stokeslet) and
\eqref{eq:ewald_TFest} (stresslet). \myChange{If $N$ were to double, this would result in a doubling of the grid size of the FFT, i.e. from $M^2$ to $2M^2$, yielding a computational cost of $\OO(N\log (N))$.}
Moreover, $P$ is set to $P=24$ \myChange{which} keeps the approximation errors of the method close to round off and is in 2D not very costly.

% !TEX root = ms.tex

% ******************************************************************************
\section{Numerical results}
\label{sec:results}
% ******************************************************************************
The numerical method concerning the drop deformation has previously been thoroughly validated using conformal mapping techniques in \cite{Palsson2018}. \revOne{These validation algorithms include analytical solutions for the steady-state of a single droplet in an extensional flow and semi-analytical solutions for the challenging time-dependent problem of deforming bubbles in close proximity.} In this section the extended numerical method described in \S\ref{sec:meth} including solid objects and walls is thoroughly tested, through convergence studies and difficult test cases. Each case is described in detail below.

\subsection{Convergence study - drop squeezing through a constriction}
\begin{figure}[h!]
   \centering
   \includegraphics[width=0.40\textwidth]{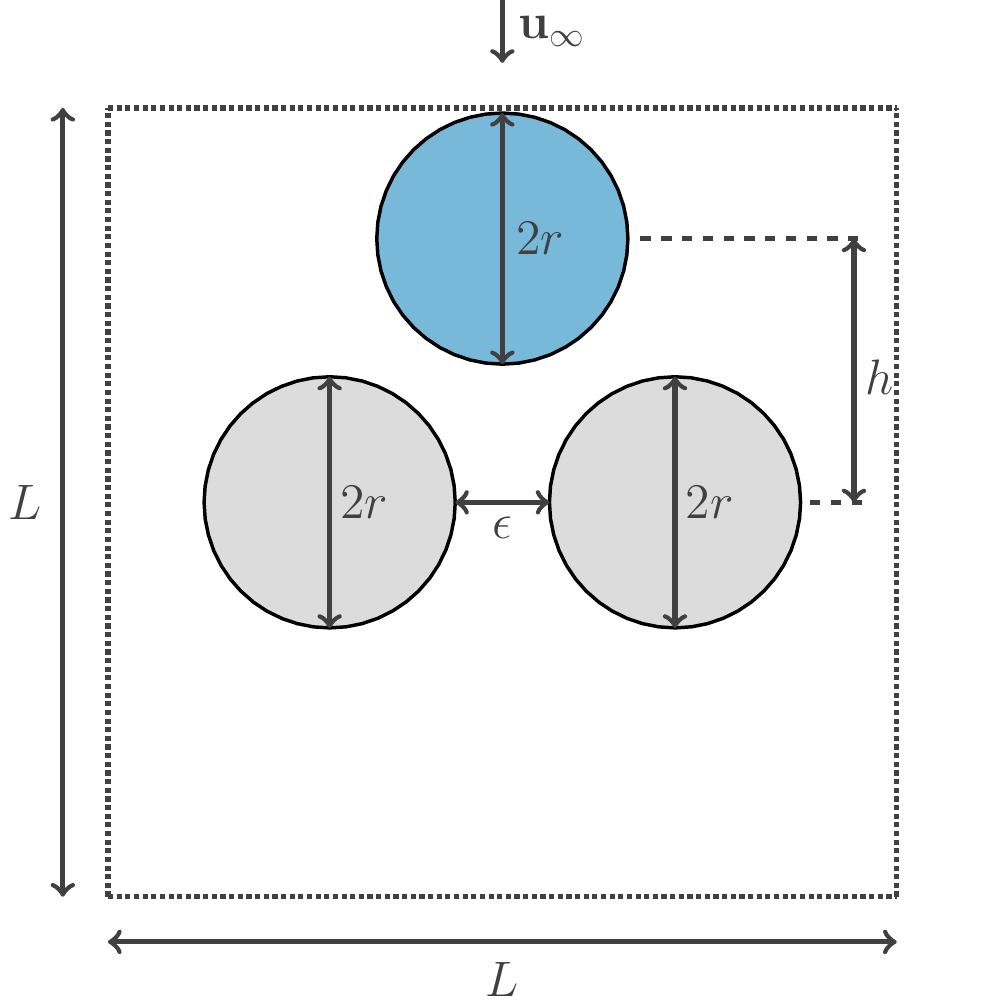}%
   \caption{Initial domain set-up for drop squeezing between two solid discs.}
   \label{fig:P2_T1_t0}
 \end{figure}

\noindent Firstly, a convergence study of drops of different viscosity ratios squeezing through a solid constriction is performed. The domain at time $t=0$ consists of two stationary solid discs of radius $r=1$ that are placed with centre points in $(\pm 1.375,0)$, i.e. their minimum distance is $\epsilon=0.75r$. Also, a drop with radius $r$ and viscosity ratio $\lambda$ is placed with its centre point in $h=(0,2.1)$. The drop is pushed down through the constriction by a imposed flow, $u_\infty=(0,-1)$. The Capillary number is set to $Ca=1$. The viscosity ratios investigated here are $\lambda=0.5$ and $\lambda=2$. The domain set-up is shown in Figure~\ref{fig:P2_T1_t0}. Initially, the drop and solids are discretised with the same number of points $\pDrop{}=\pSolid{}=N$.  Throughout the simulation the arc-length spacing $\Delta s$ is kept constant, and it is computed as the original circumference of the drop over the number of points at time $0$,
 $\Delta s = 2\pi r_{t=0}/N_{t=0}$. The results for varying $N_{t=0}$ is compared against a reference solution computed with time-stepping tolerance $10^{-10}$ and $N_{t=0}=N_{ref} = 1920$ (corresponding to $\Delta s_{ref} \approx 3\cdot 10^{-3}$) discretisation points  at times $t=1,2,3,4,5$. The tests have been performed both for clean and surfactant-covered drops. The problem is considered in a periodic setting with a reference box of size $L^2= 2\pi \times 2\pi$.

\paragraph{Comparing the solutions}
Denoting the reference solution by $z_{ref}$, and a coarser solution using $N$ discretisation points for $\tau$, the aim is to compute the difference between $z_{ref}$ and $\tau$. Previously (see e.g. \cite{Palsson2018}), the coarser solution has been upsampled to the same size as $z_{ref}$ using \texttt{FFT}s. For this convergence study, however, this will not be an optimal approach as will be explained in \S\ref{res:P2_T1_A}. Here, instead a different approach based on a normal projection onto the reference solution will be considered.

\begin{figure}[h!]
  \centering
  \includegraphics[width=0.3\textwidth]{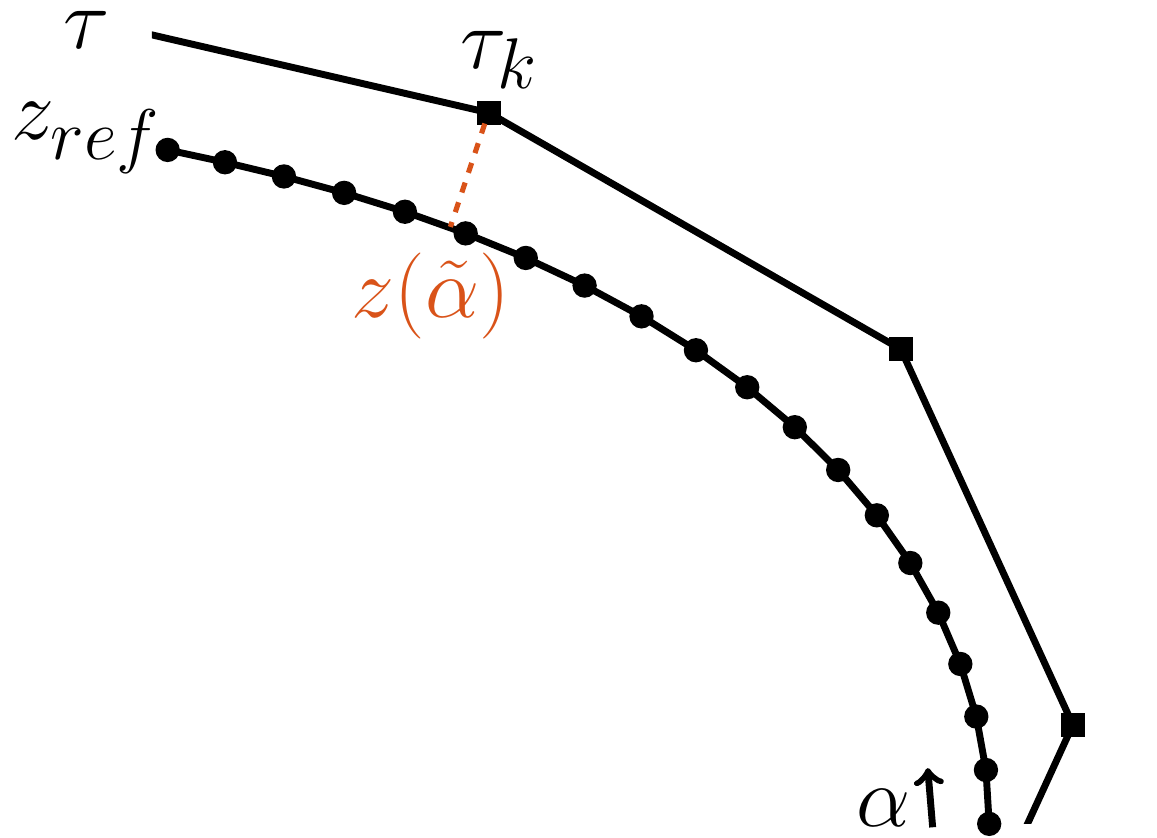}%
  \caption{Schematic of how the difference between reference solution and coarse solution is computed. Normal projection $z(\tilde{\alpha})$ of $\tau_k$ from coarse discretisation onto the reference interface is marked in red.}
  \label{fig:comperr}
\end{figure}

In essence, the closest distance between a point in the coarse discretisation and the reference solution needs to be measured. A potential tangential shift of the point is irrelevant.  A schematic of this procedure can be found in Figure~\ref{fig:comperr}. Considering the $N_{ref}$ discretisation points of $z_{ref}$, they are represented on the equidistant grid and can be seen as discrete points of a periodic function $z(\alpha)$, where $\alpha\in[0,2\pi]$. It is therefore possible to obtain their Fourier coefficients through an \texttt{FFT}. Once these coefficients have been obtained, a normal projection of a coarse discretisation point $\tau_k$ onto the reference interface can be found through a minimisation procedure. This procedure can be formulated as finding the $\tilde{\alpha}$ such that $\|\tau_k - z(\tilde{\alpha})\|_\infty$ is minimised. The difference between the reference solution and the coarse discretisation point $\tau_k$ is defined as this distance. To compute the difference between a coarse solution and the reference  solution, this procedure is repeated for all discretisation points of the coarse solution.

\subsubsection{Clean drops}
\label{res:P2_T1_A}
\begin{figure}[h!]
  \centering
  \includegraphics[width=0.4\textwidth]{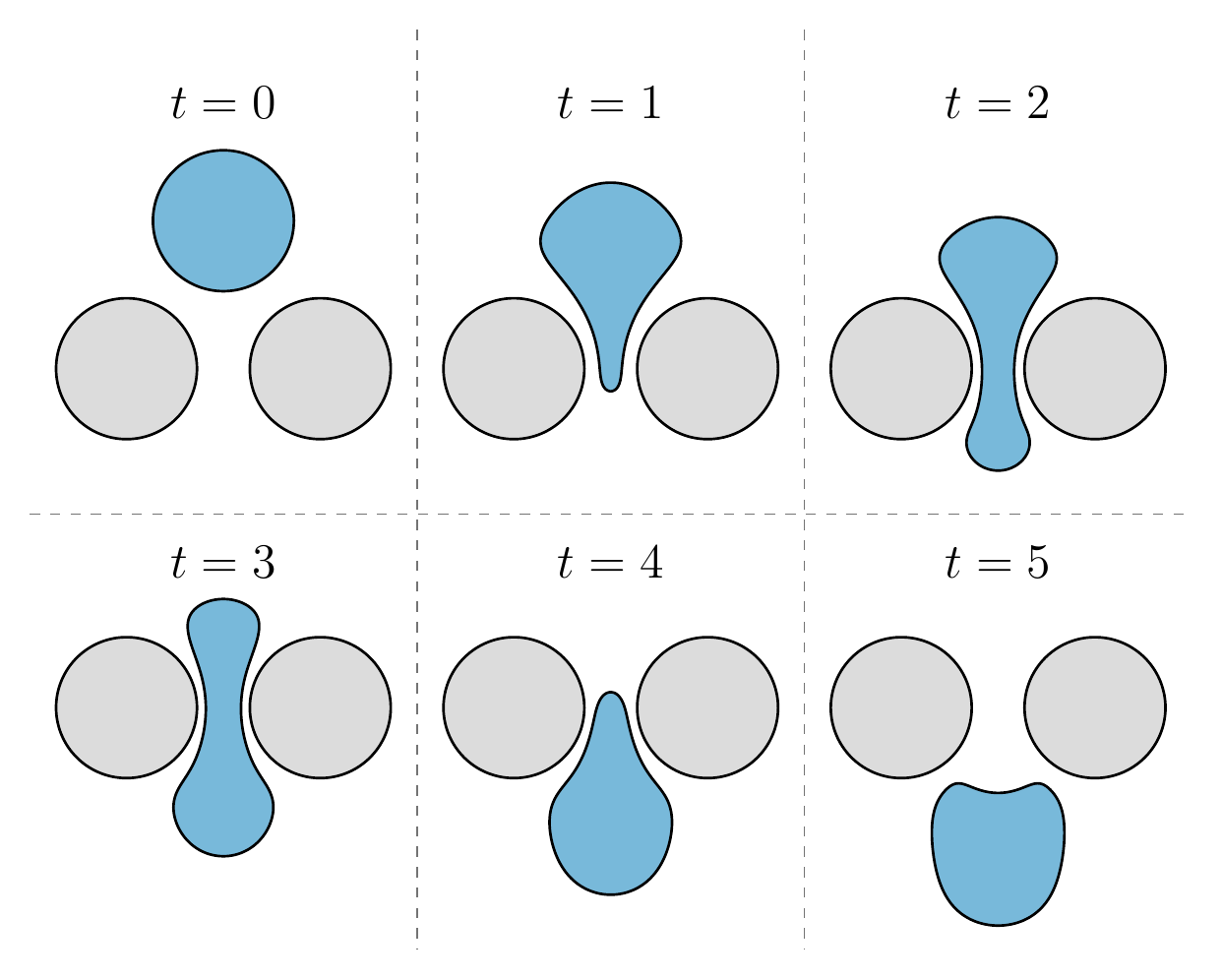}%
  \hspace{5mm}
  \includegraphics[width=0.4\textwidth]{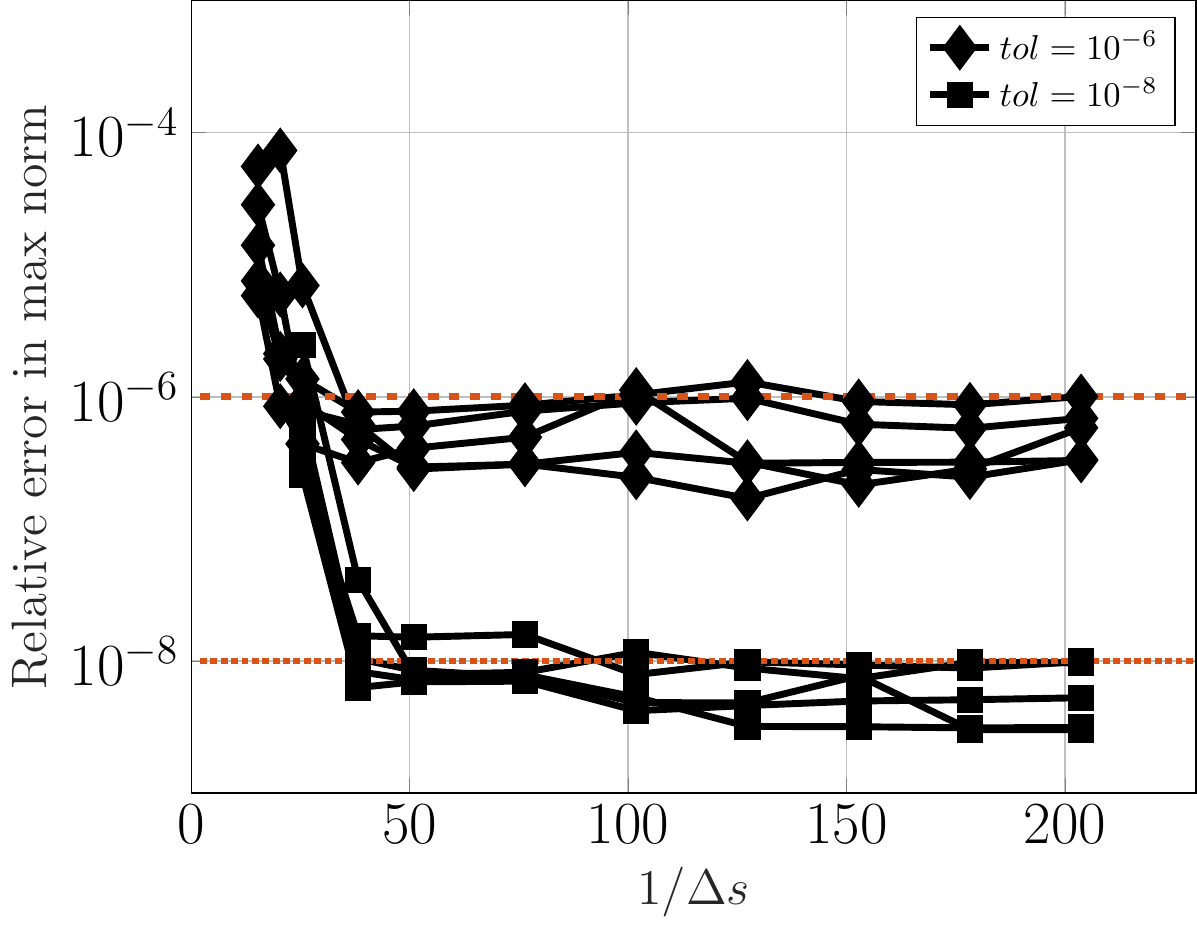}%
  \caption{Clean drop with viscosity ratio $\lambda=0.5$ squeezing through constriction. Left: Grey represents stationary solid discs, blue represents drop at times $t=0,1,2,3,4,5$. Right: relative error measured in max-norm as a function of $1/\Delta s$ for two time-stepping tolerances: $tol_1$ is marked with a red, dashed line, and corresponding errors black lines with diamonds ($ \Diamond$), $tol_2$ is marked with a red, dotted line, and corresponding errors black lines with squares ($\Box$).}
  \label{fig:P2_T1_L0p5}
\end{figure}

\begin{figure}[h!]
  \centering
  \includegraphics[width=0.4\textwidth]{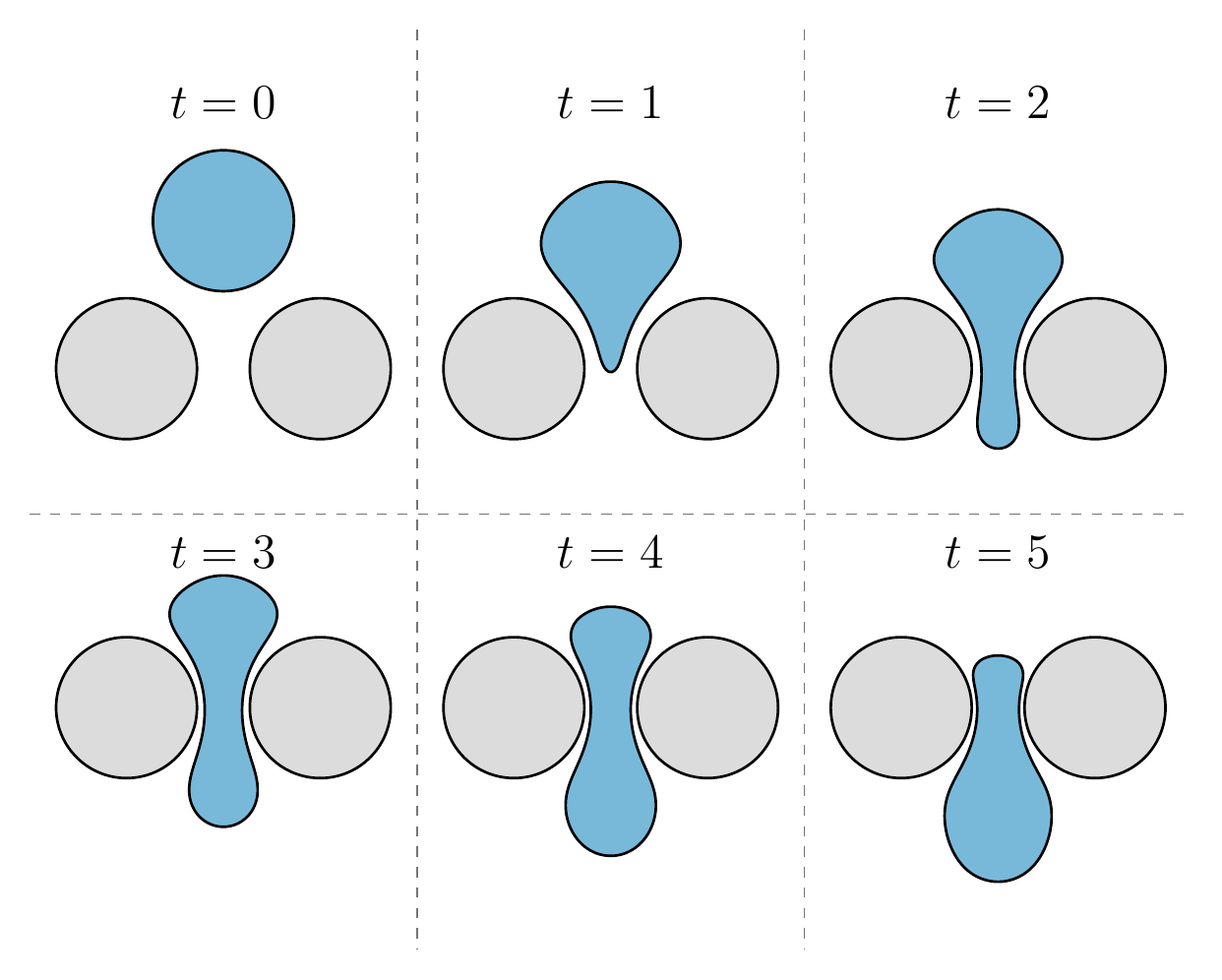}%
  \hspace{5mm}
  \includegraphics[width=0.4\textwidth]{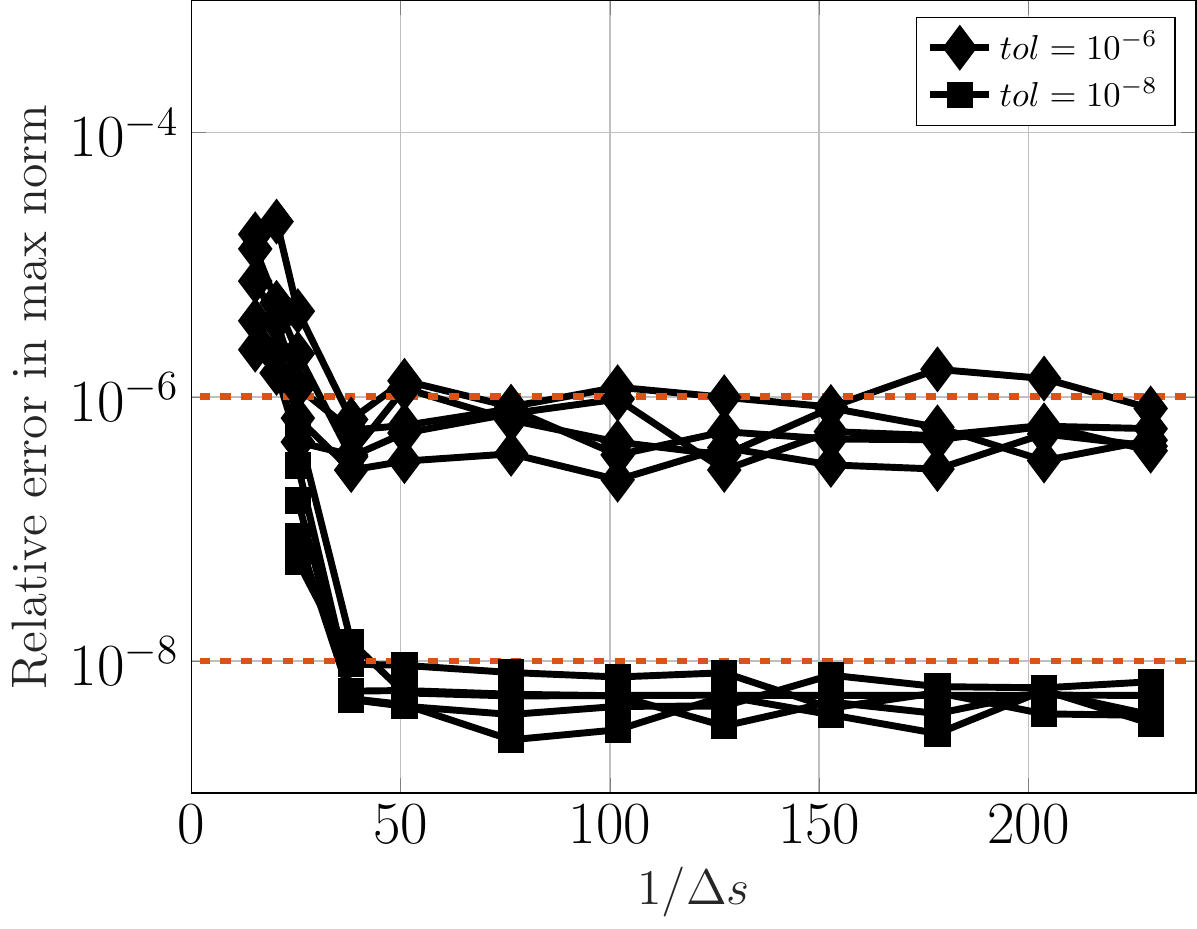}%
  \caption{Clean drop with viscosity ratio $\lambda=2$ squeezing through constriction. Left: grey represents stationary solid discs, blue represents drop at times $t=0,1,2,3,4,5$. Right: relative error measured in max-norm as a function of $1/\Delta s$ for two time-stepping tolerances: $tol_1$ is marked with a red, dashed line, and corresponding errors black lines with diamonds ($ \Diamond$), $tol_2$ is marked with a red, dotted line, and corresponding errors black lines with squares ($\Box$).}
  \label{fig:P2_T1_L2}
\end{figure}

In the case of no surfactants, the evolution of the domains and the error as a function of $1/\Delta s$ can be seen in Figure~\ref{fig:P2_T1_L0p5} for $\lambda=0.5$ and in Figure~\ref{fig:P2_T1_L2} for $\lambda=2$. The simulations were run to two time-stepping tolerances $tol_1=10^{-6}$ (marked with red dashes) and $tol_2=10^{-8}$ (marked with red dots). Throughout the simulations $\Delta s$ is approximately constant in time, due to the spatial adaptivity. Each time $t\in[1,2,3,4,5]$
is represented by a black line, with diamond markers for $tol_1$ and square markers for $tol_2$, showing the error as a function of $1/\Delta s$ at that particular time $t$. A comparison between the two viscosity ratios shows that for the same non-dimensional instance in time, the lower viscosity drop has deformed more, which is to be expected. For both cases the error decreases with an increase in $1/\Delta s$ until the time-stepping error dominates. The error is roughly the same for all times $t$. For both $\lambda=0.5$  and $\lambda=2$, the set tolerance is reached at an approximate $\Delta s \approx 0.04$ for $tol_1=10^{-6}$ and $\Delta s\approx 0.03$ for $tol_2=10^{-8}$.
\revOne{To reach $t=5$ with $tol_1$ and $\Delta s \approx 0.04$ takes in total 10 minutes in the case of $\lambda=0.5$ and 5 minutes in the case of $\lambda=2$. For these cases the average time for a time step is $3$ and $4.5$ seconds respectively. With the stricter tolerance $tol_2$ and $\Delta s\approx 0.03$ the time required becomes $16$ and $10$ minutes for $\lambda=0.5$ and $\lambda=2$ respectively, with an average time step taking $4.1$ and $4.6$ seconds. All timings are recorded on a standard desktop, with a 3.4GHz Intel Core i7 processor and 8 GB of RAM.}

\begin{figure}[h!]
  \centering
  \includegraphics[width=0.4\textwidth]{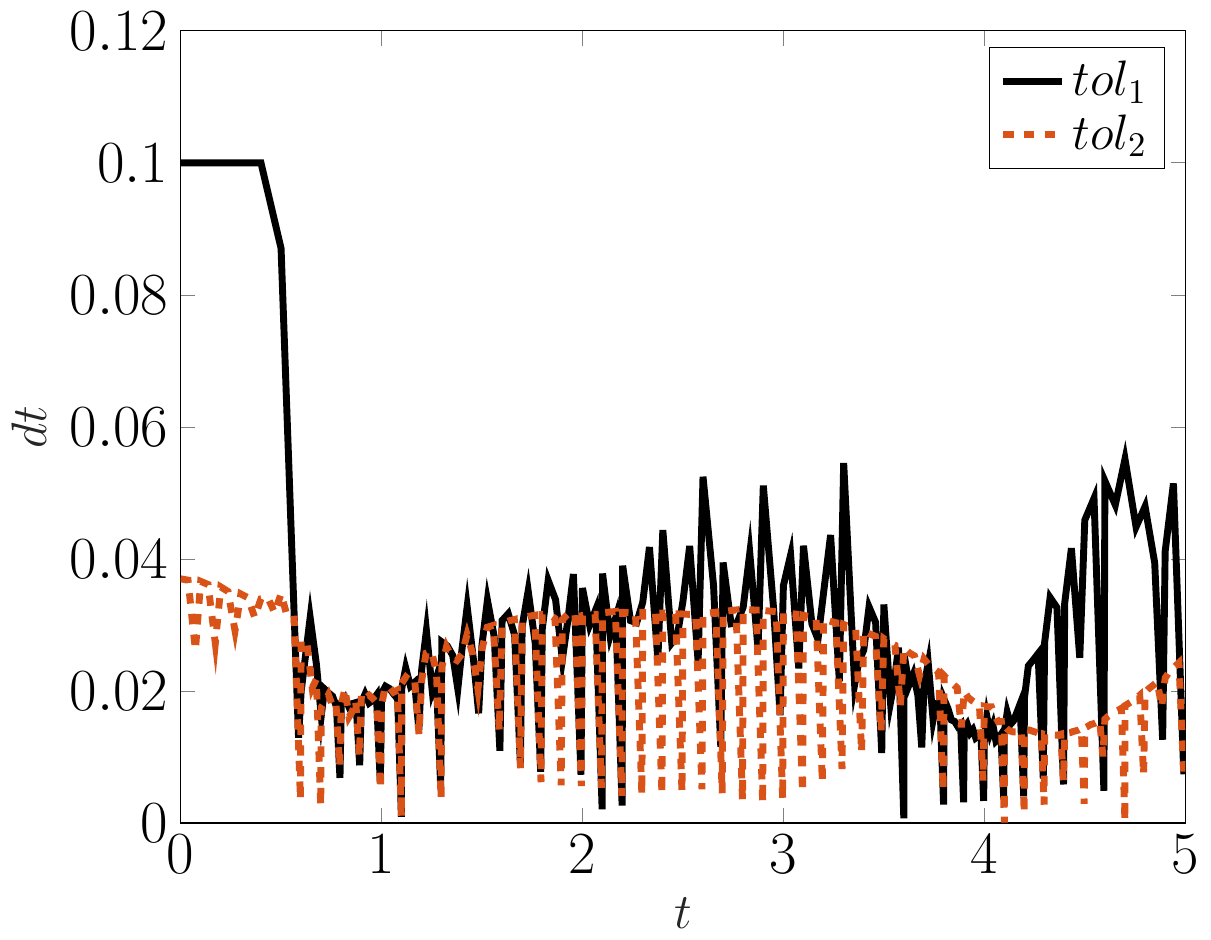}%
  \includegraphics[width=0.4\textwidth]{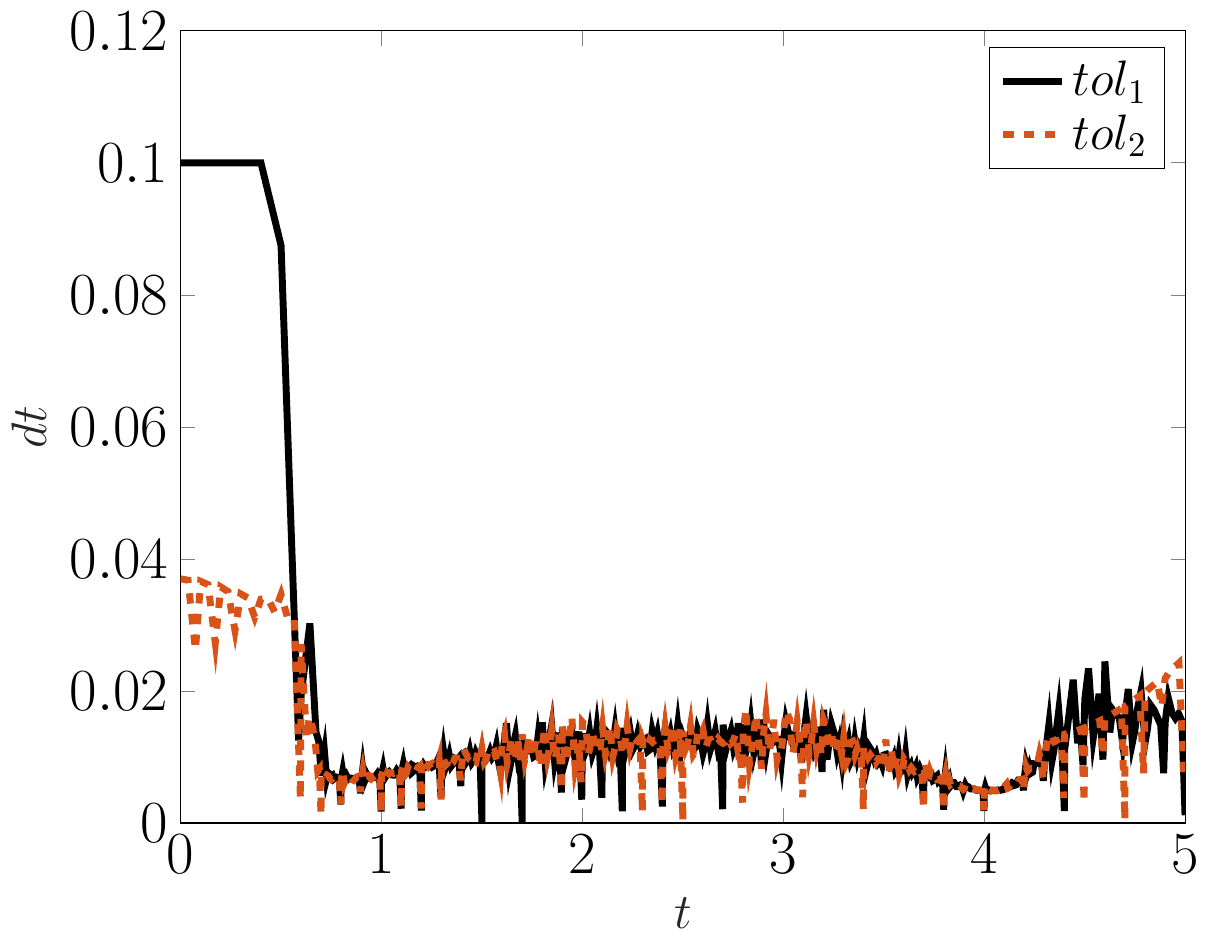}%
  \caption{Magnitude of successful time steps taken over time for a clean drop with viscosity ratio $\lambda=0.5$ squeezing through constriction. Solid, black lines for tolerance $tol_1=10^{-6}$ and dashed, red lines for tolerance $tol_2=10^{-8}$. Left: for $\Delta s=0.0131$, right: $\Delta s = 0.0049$, see Table~\ref{tab:res_timesteps_P2_T1_A}.}
  \label{fig:P2_T1_dt}
\end{figure}

\paragraph{The influence of time-stepping tolerance}
It is clear from Figure~\ref{fig:P2_T1_L0p5} and Figure~\ref{fig:P2_T1_L2} that both the correct tolerances can be reached. One would expect the less accurate solution to be cheaper to compute, however this is generally not the case.
In Table~\ref{tab:res_timesteps_P2_T1_A}, the number of failed and successful time steps up until time $t=5$ is shown for both tolerances, for four different values of $\Delta s$. \revOne{The first value is chosen such that only the larger tolerance $tol_1$ is reached and the second value is chosen such that also $tol_2$ is reached. They correspond to $1/\Delta s\approx 10$ and $15$ respectively. For the two final values of $\Delta s$ both tolerances are reached, and correspond to $1/\Delta s \approx 76$ and $1/\Delta s \approx 204$.} The cost of a simulation is defined as the number of velocity computations, i.e. the number of integral equation solves, since this is the most expensive part of the algorithm.
From Table~\ref{tab:res_timesteps_P2_T1_A} it is clear that the difference in cost between the two set tolerances is negligible, for \revOne{those values of $\Delta s$ where both tolerances can be reached}. \revOne{Note that this specific $\Delta s$ gives an under-resolved surface representation, the number of velocity computations becomes very high.}
Moreover, the cost for the larger tolerance $tol_1$ is even slightly higher than for $tol_2$, in the case of the smaller $\Delta s$. This is explained by Figure~\ref{fig:P2_T1_dt}, where the magnitude of each successful time step is shown over time, for both tolerances and $\Delta s$. For larger values of $\Delta s$, i.e. coarser discretisations, the larger tolerance $tol_1$ allows for larger time steps to be taken. However, as $\Delta s$ is decreased, the time steps become of equal size, except for at the very beginning.  The reason for this is discussed in the following paragraph, but one can conclude that using a stricter tolerance infers practically no additional cost.
\begin{table}[h!]
  \centering
  \begin{tabular}{c c | c c c }
    $\Delta s$ & Tolerance & \#failed $dt$ & \#successful $dt$ & \#velocity computations \\
    \hline\hline
    $0.0393$ & $10^{-6}$  & $40$ & $124$ & $820$ \\
    $0.0393$ & $10^{-8}$  & $2146$ & $4008$ & $30770$ \\
    \hline
    $0.0262$ & $10^{-6}$  & $21$ & $114$ & $675$ \\
    $0.0262$ & $10^{-8}$  & $6$ & $232$ & $1190$ \\
    \hline
    $0.0131$ & $10^{-6}$  & $54$ & $191$ & $1225$ \\
    $0.0131$ & $10^{-8}$  & $15$ & $239$ & $1270$ \\
    \hline
    $0.0049$ & $10^{-6}$  & $124$ & $464$ & $2940$ \\
    $0.0049$ & $10^{-8}$  & $105$ & $478$ & $2915$ \\
  \end{tabular}
  \caption{Number of failed and successful time steps up to time $t=5$ for a clean drop squeezing through a constriction, for $\lambda=0.5$, using the method and time-stepping scheme described in \S\ref{sec:meth}. \revOne{The choice of $\Delta s = 0.0393$ corresponds to the coarsest possible resolution to reach tolerance $10^{-6}$ and $\Delta s=0.0262$ the coarsest possible resolution to reach $10^{-8}$. }}
  \label{tab:res_timesteps_P2_T1_A}
\end{table}

\begin{figure}[h!]
  \centering
  \includegraphics[width=0.4\textwidth]{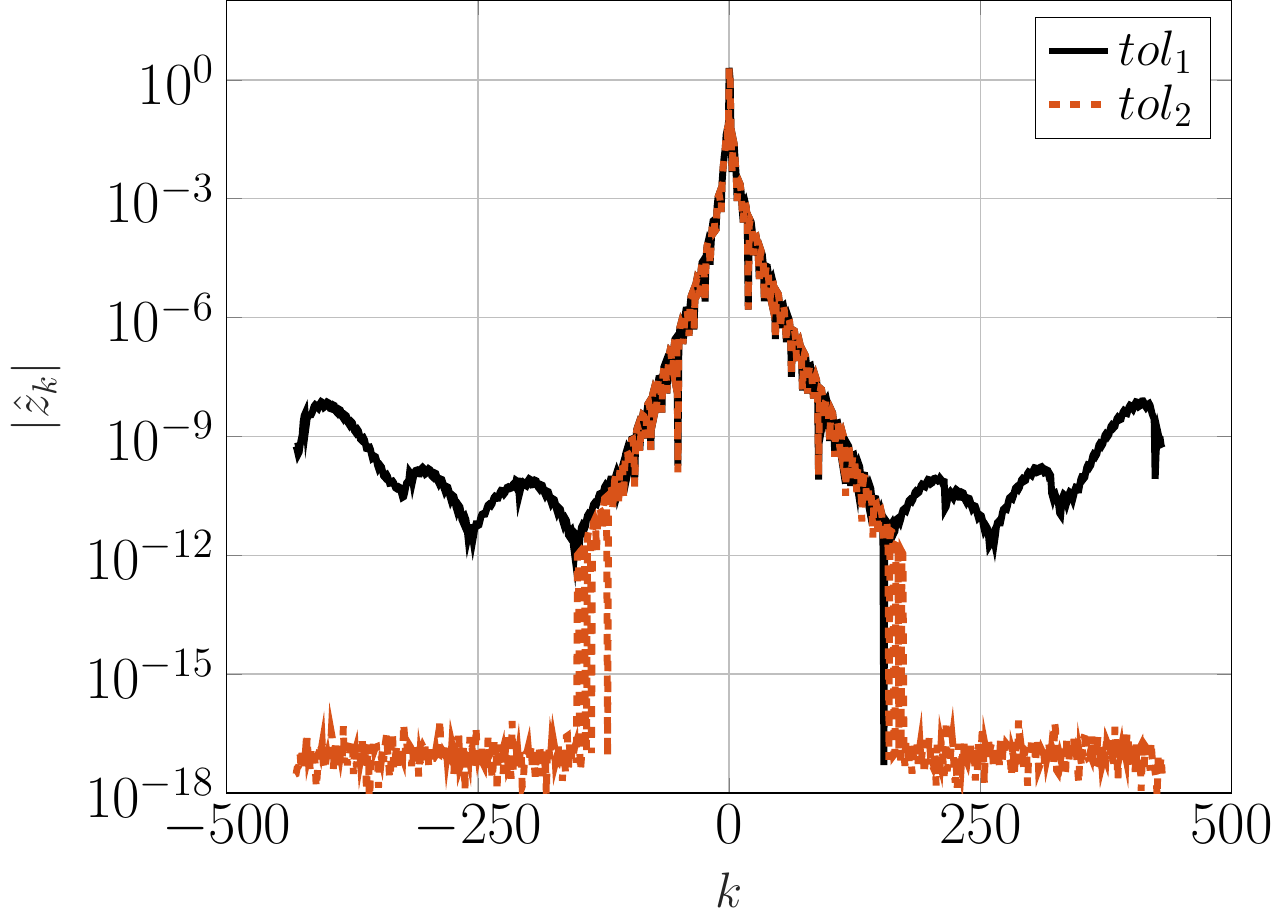}%
  \caption{Magnitude of Fourier coefficients obtained through an \texttt{FFT} of the drop shape $z(\alpha)$ for $\Delta s \approx 0.008$, for tolerance $tol_1$ (black, solid line) and $tol_2$ (red, dashed line), at time $t=5$.}
  \label{fig:P2_T1_spec}
\end{figure}

The fourth order time stepping scheme described in \S\ref{sec:meth} allows the method to take much larger time steps than the previous second order method in \cite{Palsson2017}. See Table~\ref{tab:res_timesteps_P2_T1_A_imex2} for an overview of the number of time steps taken with $tol_1$ and $1/\Delta s \approx 76$. The number of time steps taken increases from $191$ with the fourth order method to $3685$ with the second order method. This corresponds to approximately six times as many velocity evaluations. However, the required time steps are larger with the higher order method and they can come close to the stability limit. This is noticeable when regarding the equidistant spacing of the discretisation points. When taking very small time steps, such as is the case with the second order method, the points are held equidistant through time. With the larger time steps in the fourth order method, this only holds up to the time-stepping tolerance. This is the reason why computing the errors using \texttt{FFT}s for the coarser solutions and zero-padding is not viable, as it introduces additional errors. How the spectrum looks for the two tolerances is shown in Figure~\ref{fig:P2_T1_spec}. Thus, there is little computational gain when relaxing the time-stepping tolerance.

\begin{table}[h!]
  \centering
  \begin{tabular}{c c | c c c }
    $\Delta s$ & Tolerance & \#failed $dt$ & \#successful $dt$ & \#velocity computations \\
    \hline\hline
    $0.0131$ & $10^{-6}$  & 1 & 3685  & 7370  \\
  \end{tabular}
  \caption{Number of failed and successful time steps up to time $t=5$ for a clean drop squeezing through a constriction, for $\lambda=0.5$, using the second order time-stepping scheme described in \cite{Palsson2017}.}
  \label{tab:res_timesteps_P2_T1_A_imex2}
\end{table}

It is the recommendation of the authors, to in light of this information always run the simulations to the stricter time-stepping tolerance $tol_2$. This keeps the discretisation points equidistant with a clean Fourier spectrum (no ringing), and infers practically no additional cost.

\subsubsection{Surfactant-covered drops}
\begin{figure}[h!]
  \centering
  \includegraphics[width=0.8\textwidth]{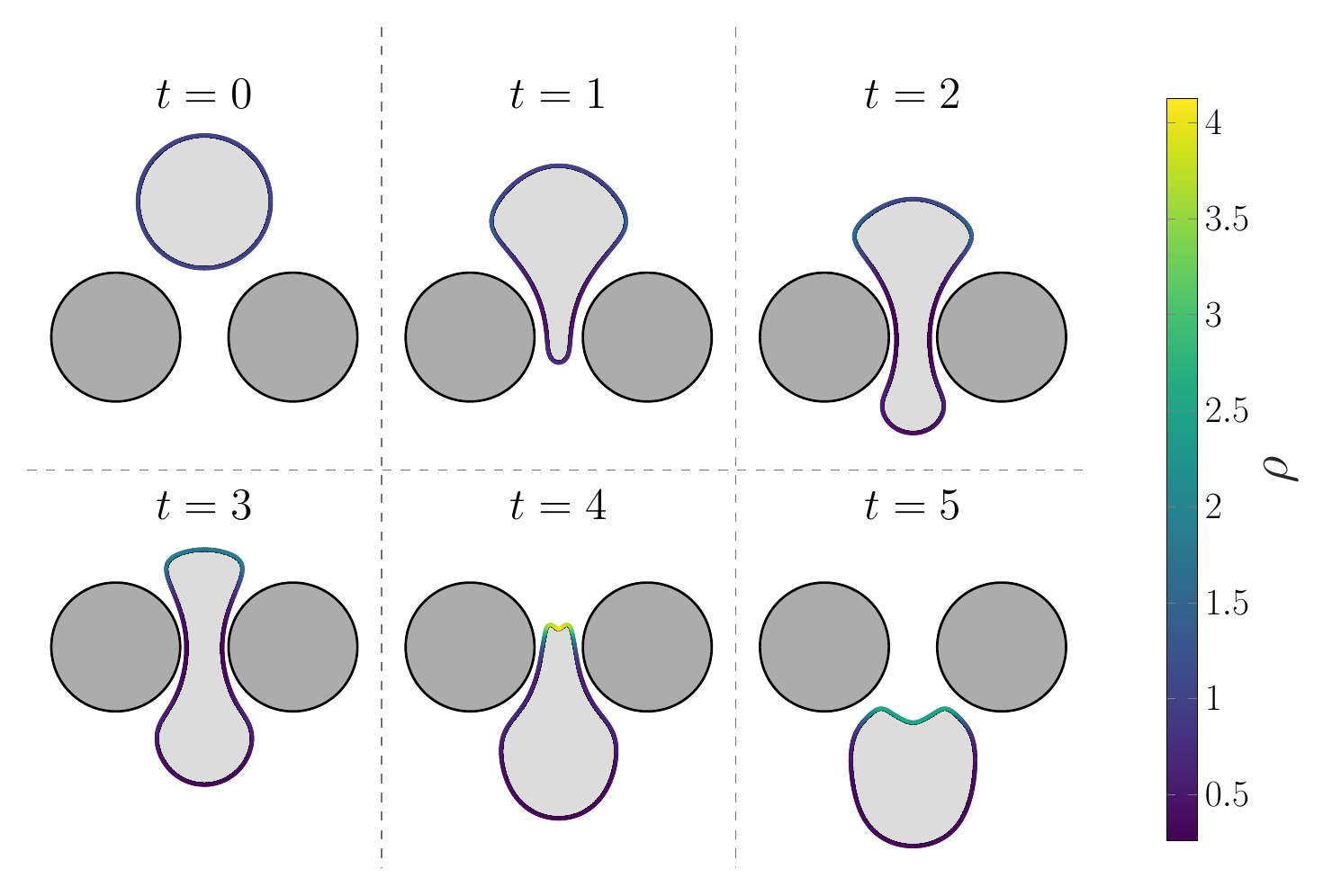}%
  \caption{Surfactant-covered drop with viscosity ratio $\lambda=0.5$ squeezing through constriction. Left: dark grey represents stationary solid discs, light grey represents drop at times $t=0,1,2,3,4,5$. Surfactant concentration on interface showed in colour.}
  \label{fig:P2_T1_B_L0p5}
\end{figure}

In this section, the same set-up as above is used with the addition of insoluble surfactants on the drop interface. The non-dimensional initial surfactant concentration is $\rho_0=1$. Furthermore, the elasticity number is set to $E=0.2$ and the P\'eclet number is set to $Pe=10$. How the drop squeezes through the constriction is shown in Figure~\ref{fig:P2_T1_B_L0p5}, together with the evolution of surfactant concentration on the interface. It is clear that the surfactant concentration affects the drop deformation, especially in places of high curvature.

\begin{figure}[h!]
  \centering
  \includegraphics[width=0.43\textwidth]{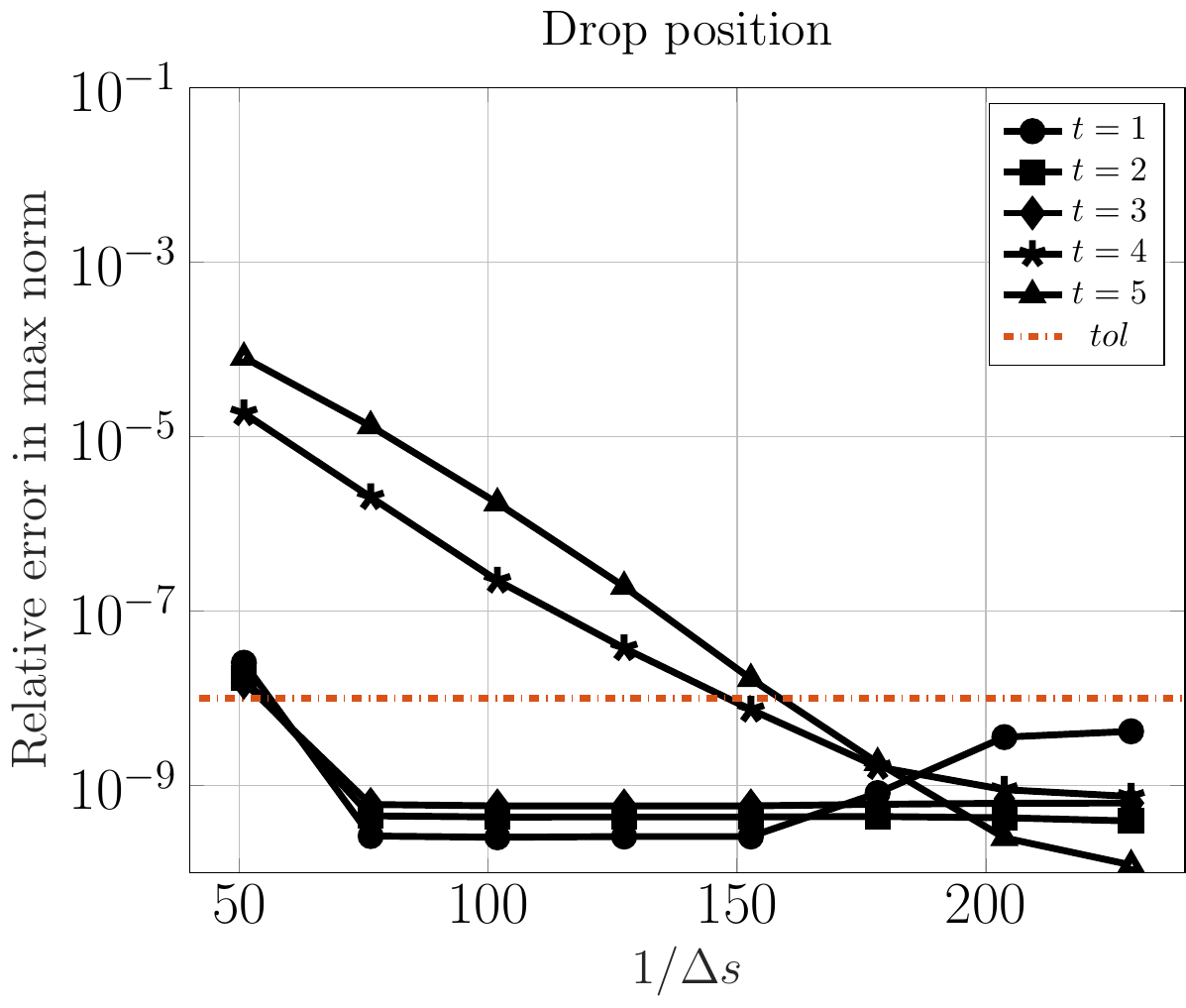}%
  \hspace{5mm}
  \includegraphics[width=0.4\textwidth]{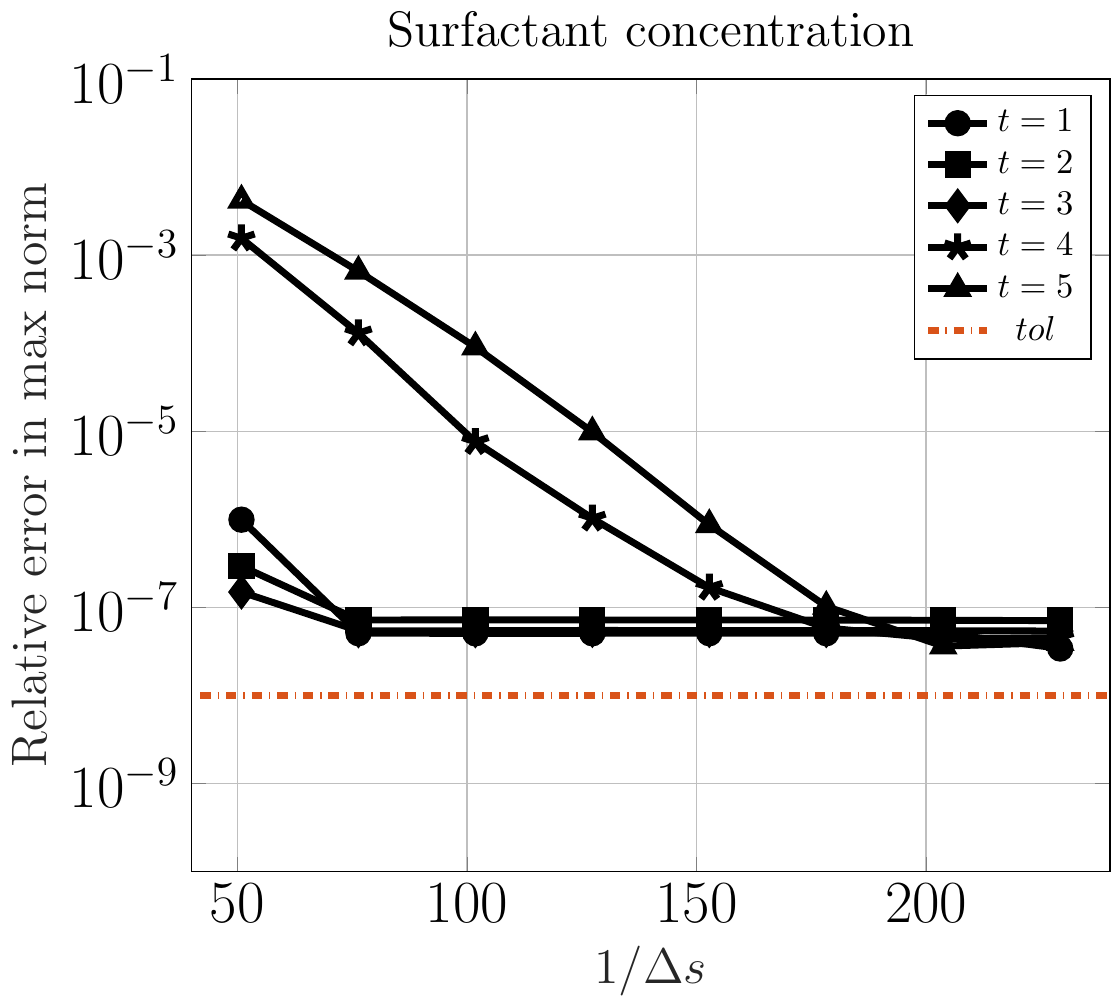}%
  \caption{Relative error measured in max-norm as a function of $1/\Delta s$ for drop position (left) and surfactant concentration (right). The different lines mark the error at time instances $t=1,2,3,4,5$.}
  \label{fig:P2_T1_B_L0p5_error}
\end{figure}

The relative error in max-norm compared to the reference solution is shown in Figure~\ref{fig:P2_T1_B_L0p5_error}. Several things should be noted with these errors. Firstly, the drop and surfactant errors are on different levels. I.e. for a set time-stepping tolerance of $10^{-8}$, the errors in position will be stable at around $10^{-9}$ whilst the surfactant concentration error level out at approximately $5\cdot 10^{-8}$. This could be easily controlled by using different time-stepping tolerances for the two quantities. Secondly, it is also clear that the error is larger for a set $\Delta s$ for the times $t=4$ and $t=5$ than for the earlier times, for both position and concentration. This is due to the increase in curvature of the drop shape, which can be seen in Figure~\ref{fig:P2_T1_B_L0p5} for time $t=4$.
This is a consequence of the fact that a smaller $\Delta s$ is needed to resolve interfaces with high curvature. Practically, this can be handled in a simulation by performing spatial adaptivity not only to keep $\Delta s$  constant, but also to decrease it as the curvature increases. This is currently not performed in the simulations.

\subsection{Multiple drops in a channel}
\begin{figure}[h!]
  \centering
  \includegraphics[width=0.4\textwidth]{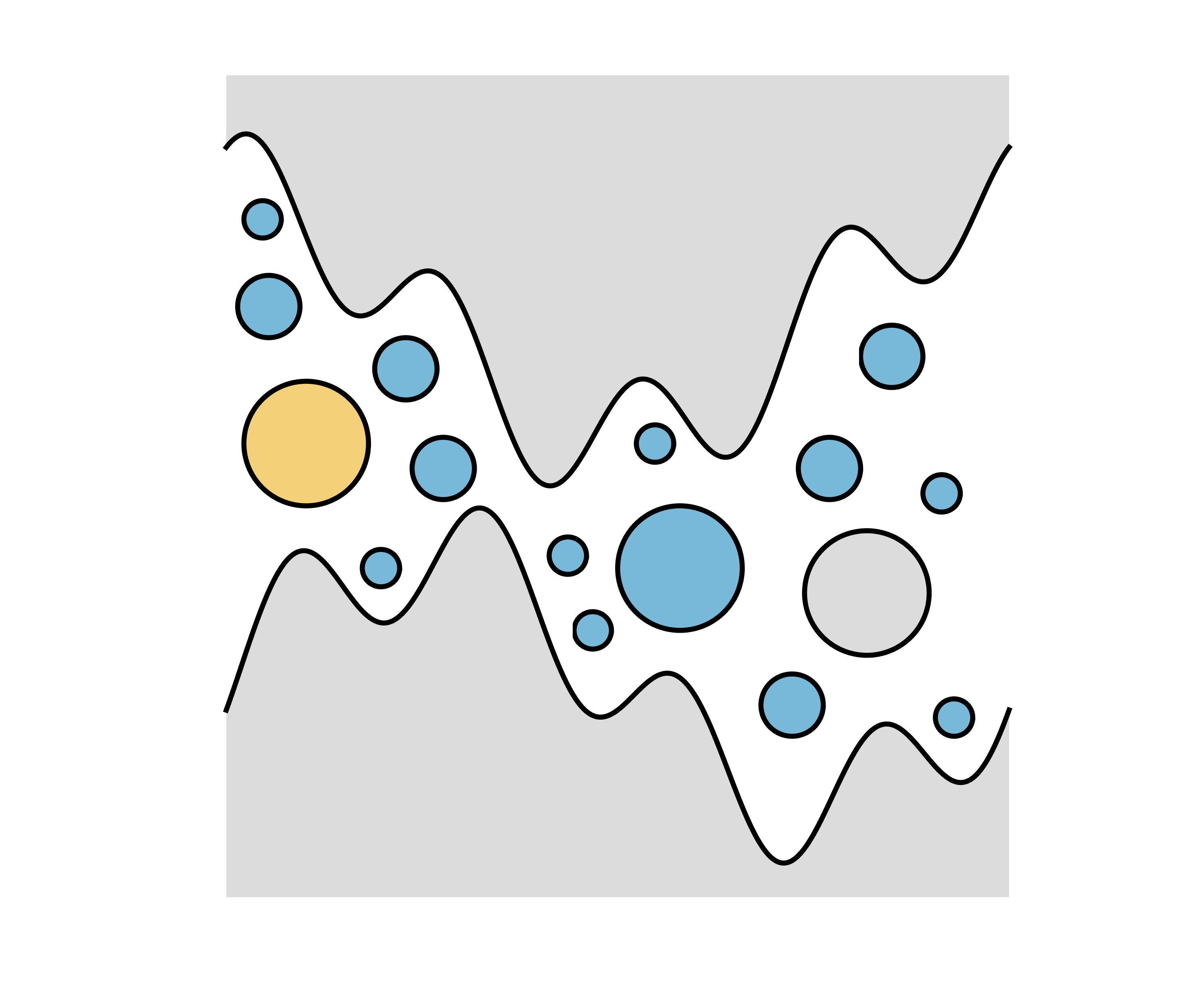}%
  \caption{Initial set-up of the channel flow. One drop is coloured yellow to facilitate the visualisation of time progressing. The periodic box is $2\pi\times 2\pi$ in size.}
  \label{fig:P2_T3_t0}
\end{figure}

\begin{figure}
  \centering
  \includegraphics[width=1\textwidth]{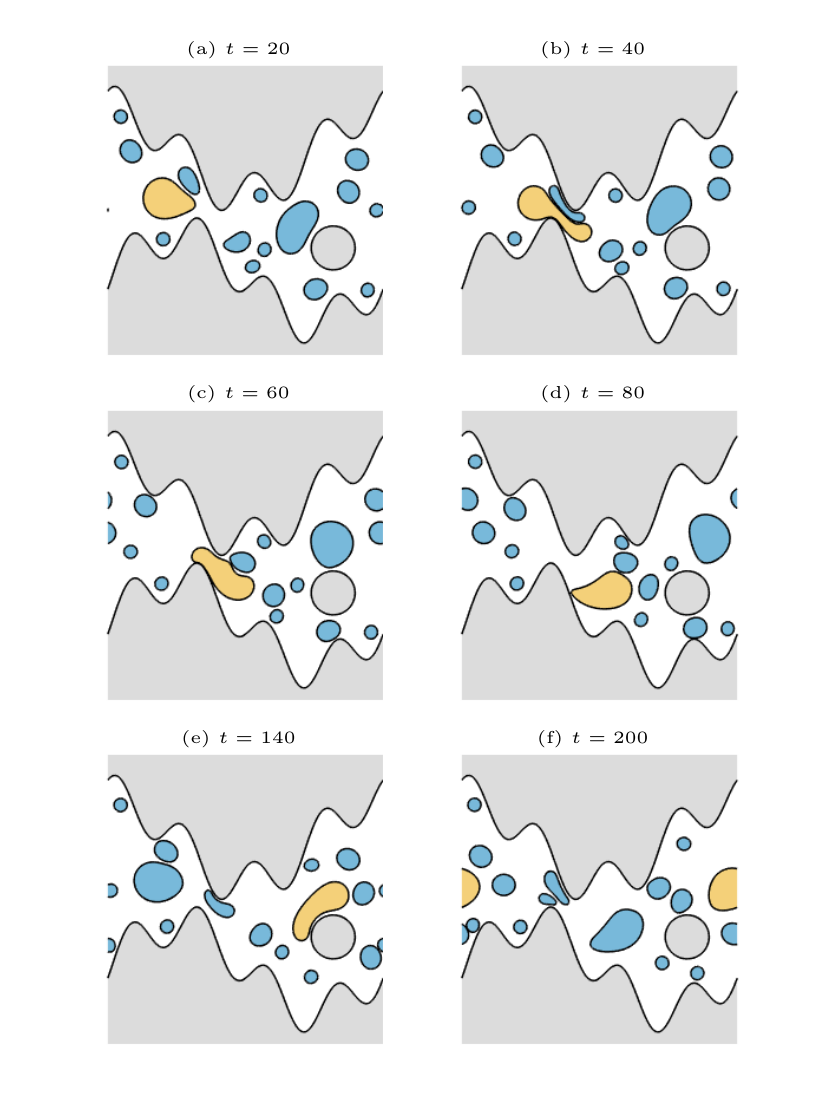}%
  \caption{Deformation over time for the channel flow. Viscosity ratio $\lambda=5$ for all drops, one drop is marked yellow to see passage of time.}
  \label{fig:P2_T3}%
\end{figure}

\begin{figure}[h!]
  \centering
  \includegraphics[width=0.4\textwidth]{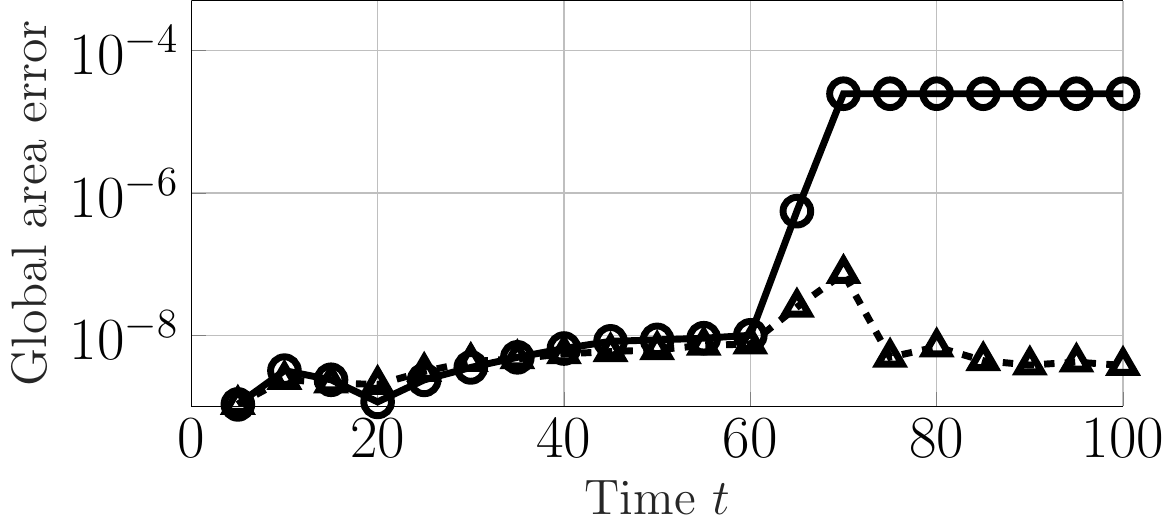}%
  \caption{Area error as a function of time for channel flow simulation, $\Delta s=0.01$ ($\circ$) and $\Delta s=0.005$ ($\triangle$).}
  \label{fig:P2_T3_ae}
\end{figure}

Here, the simulation of multiple drops in a periodic channel is shown. The set-up consists of $15$ drops of viscosity ratio $\lambda=5$: two with radius $0.5$, six with radius $0.25$ and seven with radius $0.15$. The walls are parametrised with $C^\infty$ curves and constructed through a superposition of sinus curves. Furthermore, a solid disc of radius $0.5$ is placed in the channel. The periodic length is $L=2\pi$. An added Poiseuille flow is driving the movement of the drops. The initial set-up is shown in Figure~\ref{fig:P2_T3_t0}. The minimum distance between the channel walls is $0.45$. The evolution of drops from time $t=0$ to $t=200$ is shown in Figure~\ref{fig:P2_T3}. The drops are initially discretised with 20, 10 and 6 panels for the three different drop sizes respectively, giving $\Delta s = 0.01$. The solids are discretised with a similar $\Delta s$.
Through the simulation the Capillary number is set to $\text{Ca}=5$. During the whole simulation (time $t=0$ to $t=200$), the minimum distance between two drops is $0.005$ and between a drop and a solid $0.008$. \revOne{The simulation time is approximately $58$ hours to final time on a standard desktop.} The time-stepping tolerance is set to $10^{-8}$ and the area error is less than $2.5\cdot10^{-5}$ for all times, and can be seen in Figure~\ref{fig:P2_T3_ae}. The increase in area error at time $t\approx 60$, is due to the increase in curvature in the yellow drop as seen in Figure~\ref{fig:P2_T3}(c). This higher curvature is due to the large Capillary number chosen for this simulation, and as a consequence more discretisation points to maintain a low error are needed. With twice as many points, i.e. $\Delta s=0.005$, the area error stays under $10^{-7}$
at all times, see Figure~\ref{fig:P2_T3_ae}.

\begin{figure}
  \centering
  \includegraphics[width=1\textwidth]{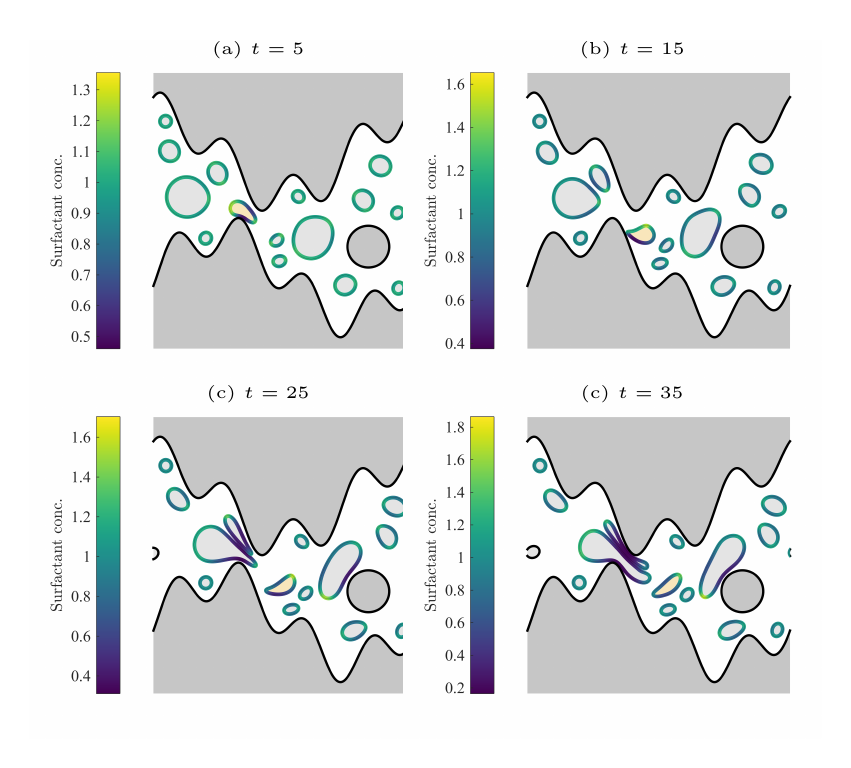}%
  \caption{Drop deformation and evolution of surfactant concentration for time instances $t=5,15,25,35$. Viscosity ratio $\lambda=5$ for all drops, one drop is marked yellow to see passage of time. Surfactant concentration of yellow drop is shown in Figure~\ref{fig:P2_T3_B_rho}.}
  \label{fig:P2_T3_B}%
\end{figure}

\begin{figure}
  \centering
  \includegraphics[width=1\textwidth]{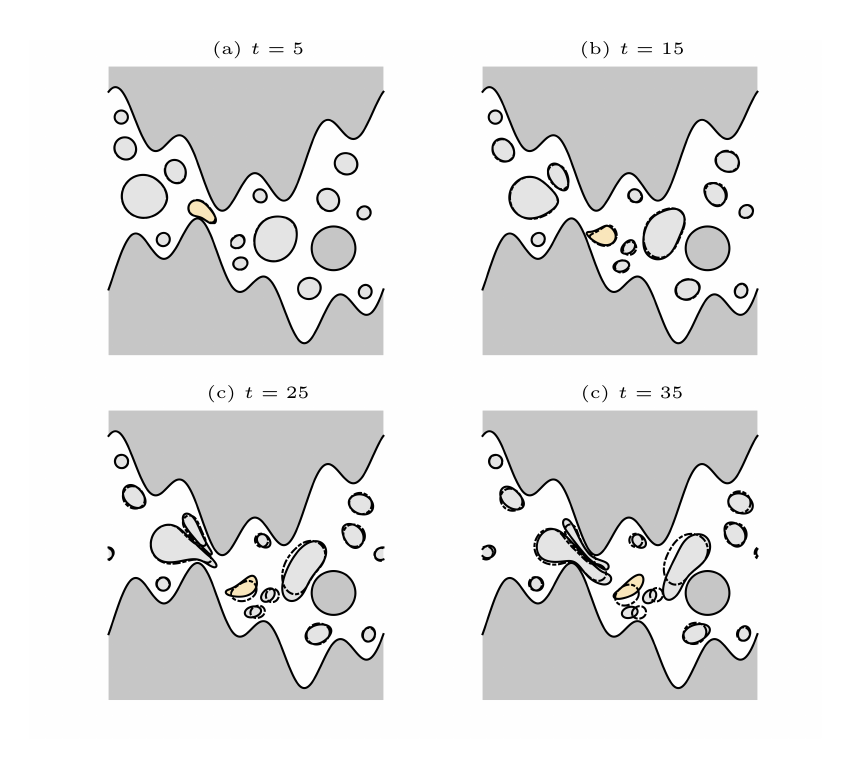}%
  \caption{Comparison between surfactant-covered case (solid lines) and clean case (dotted lines) for the channel flow.}
  \label{fig:P2_T3_Bcomp}%
\end{figure}

\begin{figure}[h!]
  \centering
  \includegraphics[width=0.8\textwidth]{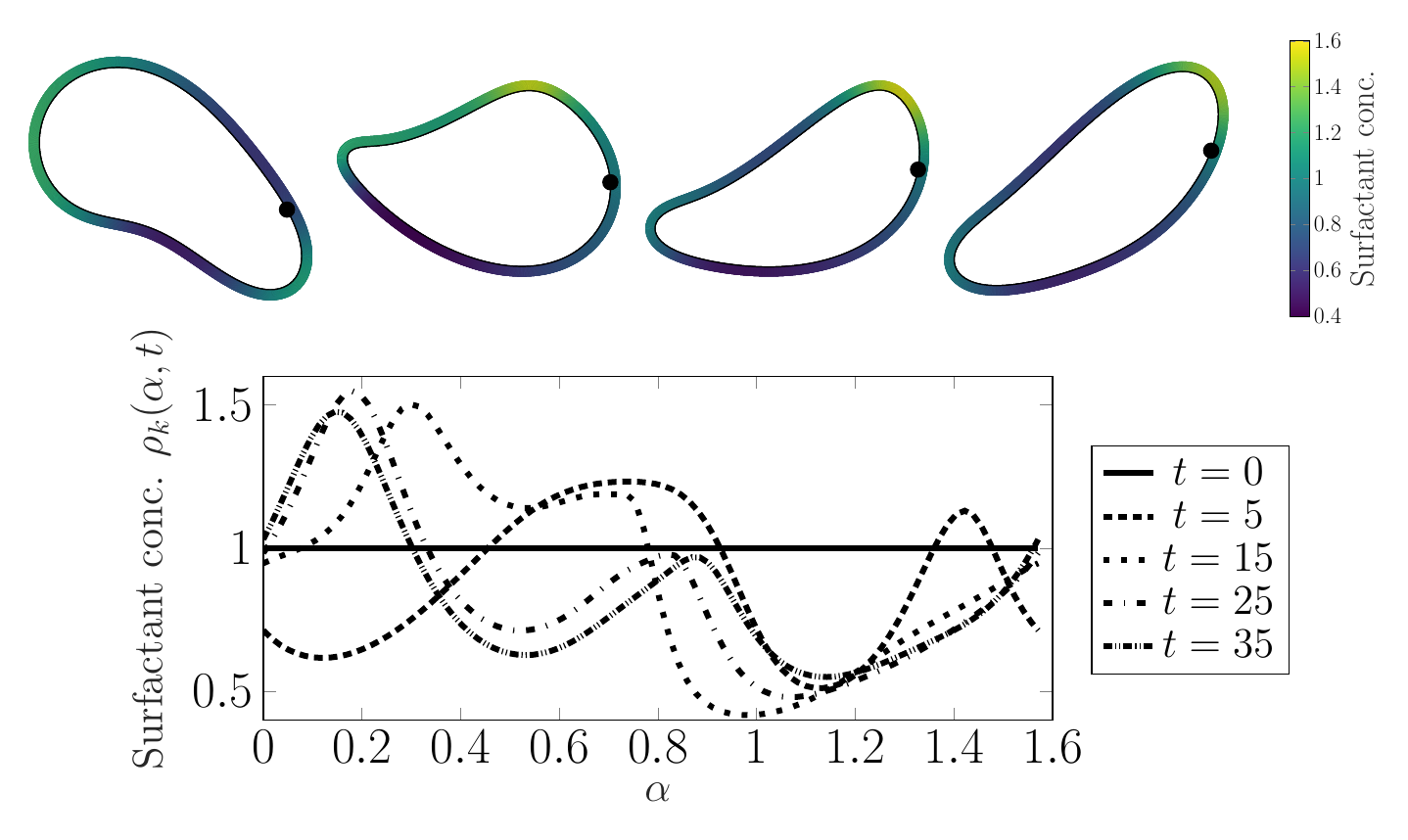}%
  \caption{Surfactant concentration for the yellow drop in Figure~\ref{fig:P2_T3_B} for times $t=5,15,25$ and $35$. Top: drop and surfactant concentration $\rho_k(\alpha,t)$,black dot marks $\alpha=0$ in bottom plot. Bottom: surfactant concentration as a function of the arc length around the drop.}
  \label{fig:P2_T3_B_rho}
\end{figure}

\begin{figure}[h!]
  \centering
  \includegraphics[width=0.4\textwidth]{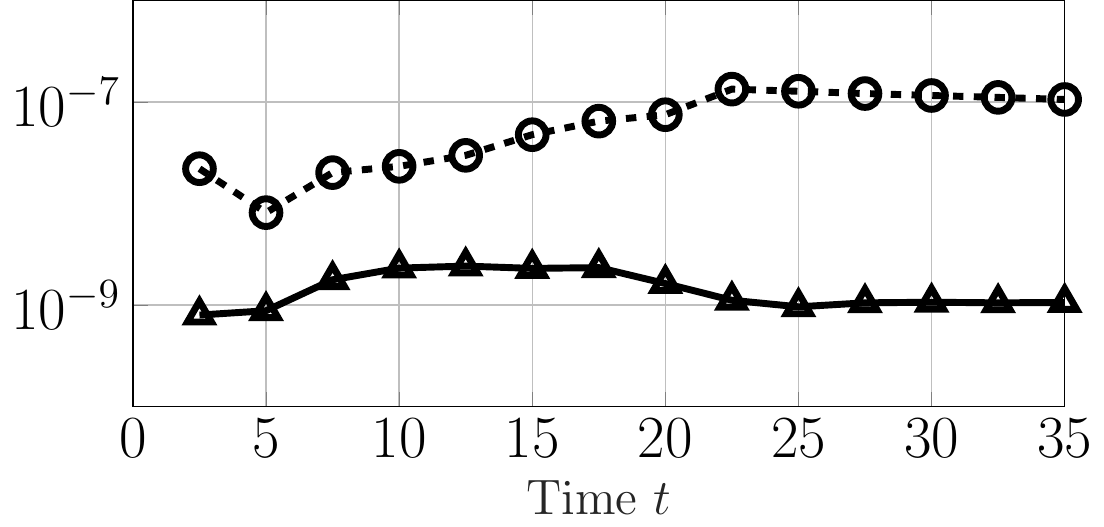}%
  \caption{Area error (\tiny{$\triangle$}\footnotesize) and surfactant concentration error ($\circ$) as a function of time for channel flow simulation.}
  \label{fig:P2_T3_B_ace}
\end{figure}

With the addition of surfactants, the surface tension of the drops is lowered and the drops therefore deform more. In Figure~\ref{fig:P2_T3_B} the deformation of the surfactant-covered drops can be seen for the case when the simulation in Figure~\ref{fig:P2_T3_t0} has been modified to include an initial surfactant concentration on all drops $\rho_0=1$, with elasticity number $E=0.5$ and P\'eclet number $Pe=1000$. A comparison of the deformation for drops with and without surfactants is shown in Figure~\ref{fig:P2_T3_Bcomp}. As can be seen, the addition of surfactants allows the drops to deform more. An example of the surfactant concentration on one drop
(the drop marked in yellow in Figure~\ref{fig:P2_T3_B}) can be seen in Figure~\ref{fig:P2_T3_B_rho} (left) for times $t=0,5,15,25,35$. The minimum distance between drops is $0.04$ and between drops and solids $0.03$.
The conservation error of the surfactants and area error of the drops can be seen in Figure~\ref{fig:P2_T3_B_ace}. It is clear that the surfactants suffer from errors greater than that of the drops position, which was already noted for the previous test case.

\section{Conclusions}
An accurate method for simulating droplets together with walls and solid stationary objects in two dimensional Stokes flow has been presented. The method allows for highly accurate solutions due to the boundary integral formulation together with the special quadrature scheme that allows for near-interaction to be well resolved. The method can handle both channel walls and solid constrictions for flow problems in a two-dimensional periodic setting. To match the high order accuracy in space, a fourth order adaptive time-stepping scheme is utilised.

In order to compute the periodic expressions efficiently, Ewald decompositions for both the Stokeslet and the stresslet have been derived and their truncation errors estimated. The decomposed expressions are then computed efficiently with the \myChange{spectral} Ewald method.

The accuracy of the method has been demonstrated through convergence tests both for clean and surfactant-covered drops. It's stability and robustness have been tested through challenging examples.

\revOne{This paper focuses solely on a boundary integral method for the two-dimensional case and demonstrates the applicability of such a method to the simulation of deforming droplets in Stokes flow. To extend this kind of method to three dimensions is a current topic of research, of which several different works have been mentioned in \S\ref{sec:intro}. One of the key challenges is then to maintain a high quality surface representation and ensure highly accurate quadrature also for drops in close proximity.}

\section{Acknowledgements}
We are very grateful to Dr. Rikard Ojala for his contributions at the initial stages of this work. This work is supported by the G\"{o}ran Gustafsson Foundation for Research in Nature and Medicine. A.K.T also gratefully acknowledges the support from the Swedish Research Council, Grant no 2015-04998.

\appendix
% !TEX root = ms.tex

% ******************************************************************************
\section{Ewald decompositions}
\label{sec:app_ewald}
% ******************************************************************************

% ******************************************************************************
\subsection{Decomposition of $\biharm(|\mb{r}|)$}
\label{sec:app_ewald_B}
% ******************************************************************************
In order to split the Green's function for the biharmonic equation, $\biharm(|\mb{r}|)$ as defined in \eqref{eq:ewald_B}, the following quantities need to be computed: $\widehat{\biharm}^F\left(|\mb{k}|,\xi\right)$ and $\biharm^R\left(\mb{r},\xi\right)$. The first correspond to the \kspace~and can be easily computed by
\begin{align}
  \widehat{\biharm}^F\left(|\mb{k}|,\xi\right) = \widehat{\gamma}(k,\xi)\widehat{\biharm}(|\mb{k}|) = \dfrac{-1}{k^4}\left(1+\dfrac{k^2}{4\xi^2}\right)e^{-k^2/4\xi^2},
  \label{eq:ewald_BF}
\end{align}
where $\mb{k}=(k_1,k_2)$, $k=|\mb{k}|$, $\widehat{\gamma}(k,\xi)$ is the Fourier transform of the Hasimoto screening function as defined in \eqref{eq:ewald_hasimoto} and  $\widehat{\biharm}(|\mb{k}|)=-1/k^4$. The \rspace~part, $\biharm^R$, is obtained through
\begin{align*}
  \biharm^R(\mb{r}) = \biharm(r) - \biharm(r)\ast\gamma(r,\xi),
\end{align*}
for the Hasimoto screening function $\gamma(r,\xi)$, where $r=|\mb{r}|$. Note that $\biharm^R$ is expected radial, but no assumption of this is made. Using that
\begin{align*}
  \widehat{\biharm}^R(k,\xi) = \widehat{\biharm}(k)-\widehat{\biharm}(k)\widehat{\gamma}(k,\xi) = \dfrac{1}{k^4}\left[ \left( 1+\dfrac{k^2}{4\xi^2}\right)e^{-k^2/4\xi^2}-1\right],
\end{align*}
it can be written
\begin{align*}
  \biharm^R(\mb{r},\xi) = \dfrac{1}{4\pi^2}\int\limits_{\mathbb{R}^2}\widehat{\biharm}^Re^{i\mb{k}\cdot\mb{r}}\dd\mb{k} = \dfrac{1}{4\pi^2}\int\limits_{\mathbb{R}^2} \dfrac{1}{k^4}\left[ \left( 1+\dfrac{k^2}{4\xi^2}\right)e^{-k^2/4\xi^2}-1\right] e^{i\mb{k}\cdot\mb{r}}\dd k_1 \dd k_2.
\end{align*}
To compute this integral, first switch to polar coordinates $(\kappa,\theta)$, where
\begin{align}
  \begin{cases}
  (k_1,k_2) &= \kappa\left(\cos(\theta),\sin(\theta)\right), \\
  \mb{r} &= r \left(\cos(\psi),\sin(\psi)\right).
\end{cases}
\label{eq:ewald_polar}
\end{align}
The integral to compute $\biharm^R$ can then be rewritten as
\begin{align*}
  \biharm^R(\mb{r},\xi) = \dfrac{1}{4\pi^2}\int\limits_{0}^{2\pi}\int\limits_{0}^{\infty}  \dfrac{1}{\kappa^4}\left[ \left( 1+\dfrac{\kappa^2}{4\xi^2}\right)e^{-\kappa^2/4\xi^2}-1\right] e^{i\kappa r(\cos(\theta-\psi))}\kappa \dd\kappa \dd\theta,
\end{align*}
which integrated over $\theta$ becomes
\begin{align*}
  \biharm^R(\mb{r},\xi) = \dfrac{1}{2\pi} \int\limits_0^\infty \dfrac{1}{\kappa^3}\left[ \left( 1+\dfrac{\kappa^2}{4\xi^2}\right)e^{-\kappa^2/4\xi^2}-1\right] J_0(\kappa r)\dd\kappa,
\end{align*}
where $J_0(x)$ is the Bessel function of the first kind of order $0$. This integral is difficult to compute, and is therefore differentiated w.r.t. $r$ according to a trick as described in \cite{Tornberg2016}. Note that $\frac{\partial J_0(\kappa r)}{\partial r} = -\kappa J_1(\kappa r)$. Differentiating $\biharm^R$ w.r.t. to $r$ then becomes
\begin{align}
  \dfrac{\partial \biharm^R}{\partial r} = \dfrac{1}{2\pi} \int\limits_{0}^\infty \dfrac{1}{\kappa^2}\left[ 1- \left( 1+\dfrac{\kappa^2}{4\xi^2}\right)e^{-\kappa^2/4\xi^2}\right]J_1(\kappa r) \dd\kappa.
  \label{eq:ewald_dBRdr}
\end{align}
To compute this integral, consider it in two steps:
\begin{align*}
  \dfrac{\partial \biharm^R}{\partial r} = \dfrac{1}{2\pi}  \int\limits_{0}^\infty  \dfrac{1-e^{-\kappa^2/4\xi^2}}{\kappa^2}J_1(\kappa r)\dd\kappa + \dfrac{1}{2\pi}\int\limits_0^\infty \dfrac{1}{4\xi^2}e^{-\kappa^2/4\xi^2}J_1(\kappa r)\dd\kappa = \dfrac{1}{8\pi}r E_1(\xi^2 r^2).
\end{align*}
Integrating this w.r.t. $r$ gives
\begin{align}
  \biharm^R(r,\xi) = \dfrac{1}{16\pi\xi^2}\left(\xi^2 r^2 E_1(\xi^2r^2)-e^{-\xi^2 r^2}\right),
  \label{eq:ewald_BR}
\end{align}
which indeed is radial. For this split to be independent of $\xi$, it is needed that $\sum_{n=1}^{M^\Lambda}f(\mb{x}_n)=0$.

% ******************************************************************************
\subsection{Decomposition of $\lapl(|\mb{r}|)$}
\label{sec:app_ewald_L}
% ******************************************************************************
In 2D, the Laplace Green's function is defined as
\begin{align*}
  \lapl(|\mb{r}|) = -\dfrac{1}{2\pi}\log(|\mb{r}|),
\end{align*}
which is the fundamental solution to $-\Delta \lapl(|\mb{r}|) = \delta(|\mb{r}|)$. When considering a sum
\begin{align*}
  u^{\lapl}(\mb{x}) = \sum\limits_{\mb{p}\in\mathbb{Z}^2}^{*}\sum\limits_{n=1}^{\pGen} \lapl(|\mb{x}-\mb{x}_n-\tau(\mb{p})|)f(\mb{x}_n),
\end{align*}
where the asterisk in the first sum corresponds to the exclusion of the term $\mb{x}-\mb{x}_k-\tau(\mb{p})=0$, the aim is to find a split into $\lapl^R(\mb{r},\xi)$ and $\widehat{\lapl}^F(|\mb{k}|,\xi)$ such that
\begin{align*}
  u^{\lapl}(\mb{x}) =& \sum\limits_{\mb{p}\in\mathbb{Z}^2}^{*}\sum\limits_{n=1}^{\pGen} \lapl^R(\mb{x}-\mb{x}_n-\tau(\mb{p}),\xi)f(\mb{x}_n) + \hdots \\
  &+ \dfrac{1}{V}\sum\limits_{\mb{k}\neq 0}\widehat{\lapl}^F(|\mb{k}|,\xi)\sum\limits_{n=1}^{\pGen} f(\mb{x}_n)e^{-i\mb{k}\cdot(\mb{x}-\mb{x}_n)} + \lim\limits_{|\mb{r}|\rightarrow 0} \left( \lapl^R(\mb{r},\xi)-\lapl(|\mb{r}|)\right)f(\mb{x}),
\end{align*}
where the last term is only included if $\mb{x}=\mb{x}_n$ for any $n\in[1,\pGen]$.
The Hasimoto split of $\lapl$ is obtained by convolving $\lapl$ with the Hasimoto screening function as defined in \eqref{eq:ewald_hasimoto}. This gives that
\begin{align*}
  \begin{cases}
  \lapl^R(\mb{r},\xi) &= \lapl(|\mb{r}|)-\lapl(|\mb{r}|)\ast \gamma(|\mb{r}|,\xi) \\
  \lapl^F(\mb{r},\xi) &=  \lapl(|\mb{r}|)\ast\gamma(|\mb{r}|,\xi).
\end{cases}
\end{align*}
Using that $\widehat{\lapl}(\mb{k}) = 1/k^2$, it follows that the \kspace~part corresponds to
\begin{align}
  \widehat{\lapl}^F(|\mb{k}|,\xi) = \widehat{\gamma}(k,\xi)\widehat{\lapl}(|\mb{k}|,\xi) = \dfrac{1}{k^2}\left(1+\dfrac{k^2}{4\xi^2}\right)e^{-k^2/4\xi^2}.
  \label{eq:ewald_LF}
\end{align}
Similarly as for $\biharm^R$ in \ref{sec:app_ewald_B}, $\widehat{\lapl}^R$ can be written as
\begin{align*}
  \widehat{\lapl}^R(|\mb{k}|,\xi) = \dfrac{1}{k^2} - \dfrac{1}{k^2}\left(1+\dfrac{k^2}{4\xi^2}\right)e^{-k^2/4\xi^2}.
\end{align*}
The inverse Fourier transform of this is
\begin{align*}
  \lapl^R(\mb{r},\xi) = \dfrac{1}{4\pi^2}\int\limits_{\mathbb{R}^2} \dfrac{1}{k^2}\left[1 - \left(1+\dfrac{k^2}{4\xi^2}\right)e^{-k^2/4\xi^2}\right]e^{i\mb{k}\cdot(\mb{x}-\mb{y})}\dd\mb{k}.
\end{align*}
Now, switching to polar coordinates \eqref{eq:ewald_polar} and integrating over $\theta$, this reads
\begin{align*}
  \lapl^R(\mb{r},\xi) &= \dfrac{1}{2\pi} \int\limits_0^\infty \dfrac{1}{\kappa} \left[ 1 - \left(1+\dfrac{\kappa^2}{4\xi^2}\right)e^{-\kappa^2/4\xi^2}\right]J_0(\kappa r)\dd\kappa =\hdots =  \dfrac{1}{4\pi}\left(-e^{-\xi^2 r^2}+E_1(\xi^2 r^2)\right).
\end{align*}
using the same trick as in \eqref{eq:ewald_dBRdr}.

% ******************************************************************************
\section{Truncation errors}
\label{sec:app_trunc}
% ******************************************************************************
For both the Stokeslet and the stresslet, there will be truncation errors when the infinite sums are truncated for computation. The \rspaces~is only evaluated for point pairs at a distance smaller than some cut-off radius $r_c$, and the \kspaces~is evaluated for all $k_1,k_2\in[-k_\infty,k_\infty]$ for some $k_\infty$. In general, denoting the truncated computation of one sum $\tilde{\mb{u}}(\mb{x})$ and the full solution $\mb{u}(\mb{x})$, the RMS truncation error is defined as
\begin{align*}
  \delta \mb{u} = \sqrt{\dfrac{1}{\pGen}\sum_{n=1}^{\pGen} |\mb{u}(\mb{x}_n)-\tilde{\mb{u}}(\mb{x}_n)|^2}.
\end{align*}

% ******************************************************************************
\subsection{Real space sum}
\label{sec:app_truncR}
% ******************************************************************************
The real space sums to compute are
\begin{align*}
  \begin{cases}
  u_j^{G,R}(\mb{x},\xi) &= \sum\limits_{\mb{p}\in\mathbb{Z}^2}^{*}\sum\limits_{n=1}^{\pGen}G^R_{jl}(\mb{x}-\mb{x}_n-\tau(\mb{p}),\xi)f_l(\mb{x}_n), \\
  u_j^{T,R}(\mb{x},\xi) &= \sum\limits_{\mb{p}\in\mathbb{Z}^2}^{*}\sum\limits_{n=1}^{\pGen} T^R_{jlm}(\mb{x}-\mb{x}_n-\tau(\mb{p}),\xi) f_l(\mb{x}_n) n_m(\mb{x}_n)
\end{cases}
\end{align*}
for the Stokeslet and stresslet respectively. The truncated sums are denoted by $\tilde{\mb{u}}_j^{G,R}$ and $\tilde{\mb{u}}_j^{T,R}$.

\paragraph{Stokeslet}
Starting with the Stokeslet, note that the error of computing such an infinite sum is given by
\begin{align*}
    \mb{u}_j^{G,R}(\mb{x}) - \tilde{\mb{u}}_j^{G,R}(\mb{x}) = \sum\limits_{s\in FL(\mb{x})} G^R_{jl}(\mb{x}-\mb{x}_s-\tau(\mb{p}),\xi)f_l(\mb{x}_s),
\end{align*}
where $FL(\mb{x}) = \left\{ (\mb{x}_s,\mb{p})\; : \; |\mb{x}-\mb{x}_s-\tau(\mb{p})|>r_c\right\}$ for a chosen cut-off radius $r_c$. The RMS error is given by
\begin{align*}
  \delta\mb{u}^{G,R} = \sqrt{\dfrac{1}{\pGen}\sum\limits_{n=1}^{\pGen}|\mb{u}^{G,R}(\mb{x}_n)-\tilde{\mb{u}}^{G,R}(\mb{x}_n)|^2}.
\end{align*}
Following \citeauthor{Kolafa1992} \cite{Kolafa1992}, this error can be approximated as
\begin{align*}
  (\delta \mb{u}^{G,R})^2 \approx \dfrac{1}{L^2} \sum_{n=1}^{\pGen} \sum_{j=1}^2 f_l^2(\mb{x}_n)\int\limits_{r>r_c}\left( G_{jl}^R\right)^2 \dd\mb{r}.
\end{align*}
The term $\sum_{j=1}^2 \left( G_{jl}^R \right)^2$ can be approximated as $2\overline{\left(G^R\right)^2}$, as follows,
\begin{align*}
  \sum_{j=1}^2 \left( G_{jl}^R \right)^2 \approx 2\overline{\left(G^R\right)^2} = \dfrac{2}{4}\sum_{j,l=1}^2 \left(G_{jl}^R\right)^2 = \dfrac{1}{2}\left[ e^{-2\xi^2 r^2}-e^{-\xi^2 r^2}E_1\left(\xi^2 r^2\right)+\dfrac{1}{2}E_1\left(\xi^2 r^2\right)^2\right].
\end{align*}
Computing $\delta\mb{u}^{G,R}$ thus reduces to
\begin{align*}
  \left( \delta \mb{u}^{G,R}\right)^2 \approx Q_GL^2\int\limits_{r>r_c} 2\overline{\left(G^R\right)^2}\dd\mb{r} = \dfrac{2\pi Q_G}{2L^2}\int\limits_{r>r_c} \left[ e^{-2\xi^2 r^2} -e^{-\xi^2 r^2}E_1\left(\xi^2 r^2\right) + \dfrac{1}{2} E_1\left(\xi^2 r^2\right)^2 \right]r\dd r.
\end{align*}
With the approximation that $E_1(x) \approx \frac{e^{-x}}{x}$ for large $x$, this can be computed and simplified as
\begin{align}
  \left(\delta \mb{u}^{G,R}\right)^2 \approx \dfrac{Q_G\pi}{4L^2} \dfrac{e^{-2\xi^2 r_c^2}}{\xi^4 r_c^2}\left(-1 + \xi^2 r_c^2\right) \approx \dfrac{Q_G\pi}{4L^2} \dfrac{e^{-2\xi^2 r_c^2}}{\xi^2},
  \label{eq:ewald_GRest2}
\end{align}
where in the last step only the leading order term in $r_c$ has been kept. Furthermore, $Q_G = \sum_{n=1}^{\pGen} f_l^2(\mb{x}_n)$. The truncation error and estimate can be seen in Figure~\ref{fig:ewald_realest}.

\paragraph{Stresslet} Similarly, for the stresslet the RMS error is given by
\begin{align*}
  \left(\delta\mb{u}^{T,R}\right)^2 = \dfrac{1}{\pGen}\sum\limits_{n=1}^{\pGen}|\mb{u}^{T,R}(\mb{x}_n)-\tilde{\mb{u}}^{T,R}(\mb{x}_n)|^2 \approx \dfrac{1}{L^2} \sum\limits_{n=1}^{\pGen} \sum\limits_{j=1}^2 f_l^2(\mb{x}_n)\normC{m}^2(\mb{x}_n) \int\limits_{r>r_c} \left(T^R_{jlm}\right)^2 \dd\mb{r}.
\end{align*}
Approximating
\begin{align*}
  \sum\limits_{j=1}^2\left(T^R_{jlm}\right)^2 \approx 2\overline{\left(T^R\right)^2} = \dfrac{2}{8}\sum\limits_{j,l,m=1}^2 \left(T^R_{jlm}\right)^2,
\end{align*}
it follows
\begin{align*}
  \int\limits_{r>r_c} \sum\limits_{j=1}^2 \left(T^R_{jlm}\right)^2 \dd\mb{r} \approx 2\int\limits_{r_c}^\infty \int\limits_0^{2\pi} \overline{\left(T^R\right)^2}  r\dd\theta \dd r = \pi\left[ e^{-2\xi^2 r_c^2}\left(5-6\xi^2+2\xi^2r_c^2\right) + 4\left(1-3\xi^2+3\xi^4\right)E_1\left(2\xi^2r_c^2\right) \right].
\end{align*}
Again, using a series expansion for $E_1(x)$ and approximating $E_1(x) \approx \frac{e^{-x}}{x}$ for large $x$, $\delta\mb{u}^{T,R}$ can be simplified into
\begin{align}
  \left( \delta\mb{u}^{T,R}\right)^2 \approx \dfrac{2\pi Q_T}{L^2} \xi^2 r_c^2 e^{-2\xi^2 r_c^2},
  \label{eq:ewald_TRest2}
\end{align}
and the error and estimate is shown in Figure~\ref{fig:ewald_realest} (right). Here, $Q_T= \sum_{n=1}^{\pGen} f_l^2(\mb{x}_n)\norm{m}^2(\mb{x}_n)$.

% ******************************************************************************
\subsection{Fourier space sum}
\label{sec:app_truncF}
% ******************************************************************************
The Fourier space sums to compute are
\begin{align*}
  \begin{cases}
      u_j^{G,F}(\mb{x},\xi)&= \dfrac{1}{V}\sum\limits_{\mb{k}\neq 0} \widehat{G}^F_{jl}\sum\limits_{n=1}^{\pGen}f_l(\mb{x}_n)e^{-i\mb{k}\cdot(\mb{x}-\mb{x}_n)}, \\
      u_j^{T,F}(\mb{x},\xi) &= \dfrac{1}{V}\sum\limits_{\mb{k}\neq 0} \widehat{T}^F_{jlm}\sum\limits_{n=1}^{\pGen}f_l(\mb{x}_n)\normC{m}(\mb{x}_n)e^{-i\mb{k}\cdot(\mb{x}-\mb{x}_n)}.
  \end{cases}
\end{align*}
If the error for a configuration of points $(\mb{x}_n,q_n)$ is defined as
\begin{align*}
  E(\mb{x}) = \sum_{n=1}^{\pGen} q_n\left(f(\mb{x}-\mb{x}_n)-\tilde{f}(\mb{x}-\mb{x}_n)\right),
\end{align*}
then the RMS can be approximated as
\begin{align*}
  \delta E^2 \approx \dfrac{1}{|\tilde{V}|} \sum\limits_n^{\pGen} q_n^2 \int\limits_{\tilde{V}}(f(\mb{r})-\tilde{f}(\mb{r}))^2\dd S,
\end{align*}
where $\tilde{V}$ is the volume enclosing all point-to-point vectors $r_{jl}=x_j-x_l$ \cite{AfKlinteberg2017a}.

\paragraph{Stokeslet} The truncation error comes from truncating the integral of the Fourier transform outside a maximum wave number, $k_\infty$, as
\begin{align*}
  \mb{u}^{G,F}(\mb{x})-\tilde{\mb{u}}^F(\mb{x}) = \dfrac{1}{(2\pi)^2}\int\limits_{|\mb{k}|>k_\infty}\widehat{\mb{G}}^F(\mb{k},\xi)\cdot\sum\limits_{n=1}^{\pGen} \mb{f}(\mb{x}_n)e^{i\mb{k}\cdot(\mb{x}-\mb{x}_n)}\dd\mb{k},
\end{align*}
for the Stokeslet. All the points $\mb{x}_n$ are contained within a periodic box of size $L\times L$, which means that $k_\infty=\frac{2\pi}{L}\frac{M}{2}$ when covering the box with $M^2$ points in a square grid. The RMS of the truncation error is computed as
\begin{align}
  \left(\delta \mb{u}^{G,F}\right)^2 = \dfrac{1}{\pGen}\sum\limits_{n=1}^{\pGen} |\mb{u}^{G,F}(\mb{x})-\tilde{\mb{u}}^{G,F}(\mb{x})|^2 \approx \sum\limits_{n=1}^{\pGen} \sum\limits_{j=1}^2\dfrac{1}{|\tilde{V}|}\int\limits_{\tilde{V}}\left(u_j^{G,F}-\tilde{u}_j^{G,F}\right)^2\dd\mb{r},
  \label{eq:ewald_duF}
\end{align}
where $\tilde{V}$ is a circle with radius $L/2$ containing all source and target points. Corresponding expressions hold for the stresslet.

To estimate the Fourier space truncation error for the Stokeslet, let
\begin{align*}
  \left(\mb{u}^{G,F}(\mb{x})-\tilde{\mb{u}}^{G,F}(\mb{x})\right)_j = e_{jl}f_l
\end{align*}
where
\begin{align*}
  e_{jl}(\mb{r}) = \dfrac{1}{L^2} \sum\limits_{\substack{\mb{k},\\ |\mb{k}|>k_\infty}} \widehat{\mb{G}}^F(\mb{k},\xi)e^{i\mb{k}\cdot\mb{r}} \approx  \dfrac{1}{L^2}\int\limits_{|\mb{k}|>k_\infty}\widehat{\mb{G}}^F(\mb{k},\xi)e^{i\mb{k}\cdot\mb{r}}\dd\mb{k},
\end{align*}
where the same approximation of the integral as in \cite{Kolafa1992} is used.
With $\widehat{G}^F_{jl}=\left(\delta_{jl}-\hat{k}_j\hat{k}_l\right)\left(1+\frac{k^2}{4\xi^2}\right)\frac{e^{-k^2/4\xi^2}}{k^2}$, for $\hat{k}_j = k_j/k$ where $k=|\mb{k}|$, and
\begin{align*}
  \sqrt{\dfrac{1}{4}\sum\limits_{j,l=1}^2 \left(\delta_{jl}-\hat{k}_j\hat{k}_l\right)^2} = \dfrac{1}{2},
\end{align*}
it holds that
\begin{align*}
  e_{jl}(\mb{r}) \approx \dfrac{1}{2 L^2} \int\limits_{k>k_\infty} \left(1+\frac{k^2}{4\xi^2}\right)\frac{e^{-k^2/4\xi^2}}{k^2} e^{i\mb{k}\cdot\mb{r}}\dd\mb{k} =
   \dfrac{1}{2L^2} \int\limits_0^{2\pi} \int\limits_{\kappa>k_\infty} \left(1+\frac{\kappa^2}{4\xi^2}\right)\frac{e^{-\kappa^2/4\xi^2}}{\kappa^2} e^{i\kappa\left(x\cos(\theta)+y\sin(\theta)\right)} \kappa \,\dd\kappa \dd\theta,
\end{align*}
where polar coordinates have been used in the last step. Integrating with respect to $\theta$ gives that
\begin{align*}
  e_{jl}(\mb{r}) \approx \dfrac{\pi}{L^2} \int\limits_{\kappa>k_\infty} \left(1+\frac{\kappa^2}{4\xi^2}\right)\frac{e^{-\kappa^2/4\xi^2}}{\kappa^2} J_0(\kappa r) \kappa \,\dd\kappa,
\end{align*}
for $r=\sqrt{x^2+y^2}$. Using $J_0(x) \approx \frac{\sqrt{2}}{\sqrt{\pi x}}$ for large $x$, the integral above is approximated as
\begin{align*}
  e_{jl}(\mb{r}) \approx \dfrac{\sqrt{2\pi}}{L^2\sqrt{r}} \int\limits_{\kappa>k_\infty} \left(1+\frac{\kappa^2}{4\xi^2}\right)\frac{e^{-\kappa^2/4\xi^2}}{\kappa^2} \dfrac{\kappa}{\sqrt{\kappa}}\, \dd\kappa =
  \dfrac{\sqrt{\pi}}{L^2\sqrt{r}2\xi^2} \Gamma\left(\dfrac{3}{4},\dfrac{k_\infty^2}{4\xi^2}\right) \approx \dfrac{\sqrt{\pi}}{L^2\sqrt{2r}}\dfrac{e^{-k_\infty^2/4\xi^2}}{\sqrt{k_\infty}},
\end{align*}
where $\Gamma(\nu,x)$ is the incomplete Gamma function, and $\Gamma(\frac{3}{4},x)\approx e^{-x}/x^{1/4}$ for large $x$.
Inserting this into \eqref{eq:ewald_duF},
\begin{align*}
  \left(\delta\mb{u}^{G,F}\right)^2 \approx \dfrac{8Q_G}{L^2 \pi} \int\limits_0^{2\pi}\int_0^{L/2} e_{jl}^2(r)r\,\dd r\,\dd\theta \approx \dfrac{8Q_G}{L^2\pi} \dfrac{\pi}{2L^4} \dfrac{e^{-2k_\infty^2/4\xi^2}}{k_\infty} \int\limits_0^{2\pi} \int\limits_0^{L/2} \dfrac{1}{r}r\,\dd r\,\dd\theta.
\end{align*}
This gives that
\begin{align}
  \left(\delta\mb{u}^{G,F}\right)^2 \approx \dfrac{4Q_G\pi}{L^5k_\infty}e^{-2k^2_\infty/4\xi^2}.
  \label{eq:ewald_GFest2}
\end{align}
The truncation error and estimate can be seen in Figure~\ref{fig:ewald_Fest} (left).

\paragraph{Stresslet} The approach to derive a truncation error estimate for the Fourier space sum of the stresslet is similar to that of the Stokeslet above. First, let
\begin{align*}
  \left(\mb{u}^{T,F}(\mb{x})-\tilde{\mb{u}}^{T,F}(\mb{x})\right)_j = e_{jlm}f_l\normC{m}
\end{align*}
where
\begin{align*}
  e_{jlm}(\mb{r}) = \dfrac{1}{L^2}\sum\limits_{\substack{\mb{k},\\ |\mb{k}|>k_\infty}}
\widehat{\mb{T}}^F(\mb{k},\xi)e^{i\mb{k}\cdot\mb{r}} \approx \dfrac{1}{L^2}
  \int\limits_{|\mb{k}|>k_\infty}\widehat{\mb{T}}^F(\mb{k},\xi)e^{i\mb{k}\cdot\mb{r}}\,\dd\mb{k}.
\end{align*}
Using
\[
\hat{T}_{jlm}^F = \left(1+\frac{k^2}{4\xi^2}\right)\left(\delta_{jl}\hat{k}_m + \delta_{jm}\hat{k}_l + \delta_{lm}\hat{k}_j-2\hat{k}_j\hat{k}_l\hat{k}_m\right)\frac{e^{-k^2/4\xi^2}}{k},
\]
and similar steps as above, with
\begin{align*}
  \sqrt{\dfrac{1}{8}\sum\limits_{j,l,m=1}^2 \left(\delta_{jl}\hat{k}_m + \delta_{jm}\hat{k}_l + \delta_{lm}\hat{k}_j-2\hat{k}_j\hat{k}_l\hat{k}_m\right)^2} = \dfrac{1}{\sqrt{2}},
\end{align*}
$e_{jlm}$ can be approximated as
\begin{align*}
  e_{jlm}(\mb{r}) \approx \dfrac{\sqrt{\pi}}{4L^2\sqrt{r}\xi^2}k_\infty^{5/2}E_{-1/4}\left(\dfrac{k_\infty^2}{4\xi^2}\right) \approx \dfrac{\sqrt{\pi k_\infty}}{L^2\sqrt{r}}e^{-k_\infty^2/4/xi^2},
\end{align*}
with the use of $E_{-1/4}(x) \approx e^{-x}/x$ for large $x$.
The RMS error then becomes
\begin{align}
  \left(\partial \mb{u}^{T,F}\right)^2 \approx \sum\limits_{n=1}^{\pGen} \sum\limits_{j=1}^2 f_l^2(\mb{x}_n)\normC{m}(\mb{x}_n)\dfrac{1}{|\tilde{V}|}\int\limits_{\tilde{V}} e_{jlm}^2r \,\dd\mb{r} \approx \dfrac{8\pi Q_T}{L^5}k_\infty e^{-2k_\infty^2/4\xi^2}.
  \label{eq:ewald_TFest2}
\end{align}
In this last step the same simplifications as for the Stokeslet have been applied, as well as $\text{erfc}(x)\approx\frac{e^{-x^2}}{\sqrt{\pi}x}$ for large $x$. The truncation error and estimate can be seen in Figure~\ref{fig:ewald_Fest} (right).

% !TEX root = ms.tex

% ******************************************************************************
\section{Special quadrature}
\label{sec:app_sq}
% ******************************************************************************
\noindent Here the recursion formulas for computing $p_\ell$, $q_\ell$ and $r_\ell$ in \eqref{eq:meth_pqr} are given. For $p_\ell$, it holds
\begin{align*}
  \begin{cases}
    p_0 &= \int\limits_{-1}^1 \dfrac{\dd\tau}{\tau-z} = \log(1-z) -\log(-1-z), \\
    p_\ell &= zp_{\ell-1} + \dfrac{1-(-1)^\ell}{\ell}, \; \ell=1,\hdots,n-1.
  \end{cases}
\end{align*}
Note that if $z$ is within the contour created by the transformed panel $\Lambda$ and the real axis from $-1$ to $1$, then a residue of $2\pi i$ must be added or subtracted from $p_0$ depending on if $z$ has a positive or negative imaginary part respectively. Similarly, to compute $q_l$, the recursion is
\begin{align*}
  \begin{cases}
    q_0 &= \int\limits_{-1}^1 \dfrac{\dd\tau}{(\tau-z)^2} = -\dfrac{1}{1+z} - \dfrac{1}{1-z}, \\
    q_\ell &= zq_{\ell-1} + p_\ell, \; \ell=1,\hdots,n-1.
  \end{cases}
\end{align*}
Finally, the recursion for $r_\ell$ is given by
\begin{align*}
  r_\ell = \dfrac{\log(1-z)-(-1)^{\ell+1}\log(-1-z)-p_{\ell+1}}{\ell+1} + \log(\gamma)\dfrac{1-(-1)^{\ell+1}}{\ell+1},
\end{align*}
for $\gamma = (\tau^{(2)}-\tau^{(1)})/2$ where $\tau^{(1)}$ and $\tau^{(2)}$ are the endpoints of the panel $\Lambda$ before it was scaled and rotated. This is needed as the $\log$-kernel is not scale and rotation invariant.

% \input{introduction/intro5}
% \input{introduction/formulation}
% \input{introduction/bie}
% \input{methods/meth_bie}
% \input{methods/meth_est}
% \input{methods/meth_spec}
% \input{methods/nummethods}
% \input{methods/validation4}
% \input{results/results}
% \input{conclusions/conclusions2}
% \input{conclusions/acknowledgement}
%
% \appendix
% \input{methods/appendix_est}
% \input{methods/appendix_spec}

\newpage

\section*{References}
%% References{}
%% Usage of \cite is as follows:{}
\bibliographystyle{plainnat}
% \bibliography{library.bib}
\bibliography{saraslib_abb.bib}

\end{document}